\documentclass[10pt]{article} 

\usepackage{lineno}

\usepackage{amssymb} 
\usepackage{etoolbox}
\usepackage{amsmath} 
\usepackage{amsthm} 
\usepackage{hyperref} 
\usepackage{amsfonts} 
\usepackage{color} 
\usepackage{epsfig} 
\usepackage{graphics,subfigure} 
\usepackage{graphicx} 
\usepackage{float} 
\usepackage{epsf} 
\usepackage{mathrsfs} 
\usepackage{epstopdf}

\newtheoremstyle{thmm}{1.5ex plus 1ex minus .2ex}{1.5ex plus 1ex minus
.2ex}{\rmfamily}{}{\bfseries}{}{1em}{} \theoremstyle{thmm} 
\newtheorem{theorem}{Theorem}[section] 
\newtheorem{lemma}{Lemma}[section] 
\newtheorem{remark}{Remark}[section] 
\newtheorem{algorithm}{Algorithm}[section] 

\allowdisplaybreaks 
\definecolor{myGreen}{rgb}{0.9, 0.99, 0.9}
\newcommand{\nn}{\nonumber} 
\def\refe#1{(\ref{#1})} 

\def\fr{\frac} 

\def\iy{\infty}

\def\O{\Omega}

\begin{document}
\date{\today}
\title{\bf A Crank-Nicolson leap-frog scheme for the unsteady incompressible magnetohydrodynamics equations\thanks{This work is supported by the National Natural Science Foundation of China (No. 11971152 \& 12126318) and the Fundamental
Research Funds for the Universities of Henan Province (No.NSFRF210437).}}
\author{Zhiyong Si
\footnote{School of Mathematics and Information Science, Henan Polytechnic University, 454003, Jiaozuo, P.R. China.  {\tt sizhiyong@hpu.edu.cn}
}\, ,\, Mingyi Wang\footnotemark[2], Yunxia Wang\footnotemark[2]}

\maketitle 

\begin{abstract}
This paper presents a Crank-Nicolson leap-frog (CNLF) scheme for the unsteady incompressible magnetohydrodynamics (MHD) equations. The spatial discretization adopts the Galerkin finite element method (FEM), and the temporal discretization employs the CNLF method for linear terms and the semi-implicit method for nonlinear terms. The first step uses Stokes style's scheme, the second step employs the Crank-Nicolson extrapolation scheme, and others apply the CNLF scheme. We testify that the fully discrete scheme is stable and convergent when the time step is less than or equal to a positive constant. The second order $L^{2}$ error estimates can be derived by a novel negative norm technique. The numerical results are consistent with our theoretical analysis, which indicates that the method has an optimal convergence order. Therefore, the scheme is effective for different parameters.

	\end{abstract}
	
	
	\section{Introduction}
	\setcounter{equation}{0}
	This work is considered for the following unsteady incompressible MHD system which coupled the incompressible Navier-Stokes equations and Maxwell's equations under the influence of external body forces. We focus on studying a CNLF scheme for the unsteady incompressible MHD equations.
	  The non-dimensional system of equations, which governs the flow under consideration is given as follows
	\begin{eqnarray}\label{1.1}
	\left\{\begin {array}{lll}
	\textbf{u}_t-\nu\Delta \textbf{u}+\textbf{u}\cdot\nabla\textbf{u}+\nabla p+\mu \textbf{H}\times\operatorname{curl}
	\textbf{H}=\textbf{f}, \quad \text{in }\Omega\times(0,T], \\
	\mu\textbf{H}_t+\sigma^{-1}\operatorname{curl}(\operatorname{curl}\textbf{H})-\mu\operatorname{curl}(\textbf{u}\times\textbf{H})=\textbf{g}, \quad \text{in }\Omega\times(0,T],\\
	\operatorname{div} \textbf{u}=0, \ \  \operatorname{div}\textbf{H}=0,  \quad \text{in }\Omega\times(0,T], \\
	\end{array}\right.
	\end{eqnarray}
	for $(x,t)\in Q_T$. Here $\Omega$ is an open bounded domain in $R^{d}$, $d=2, 3$ with a smooth boundary, $Q_T=\Omega\times (0,T]$ and $T\in(0,\iy)$. Where $\textbf{u}$, $p$ and \textbf{H} stand for  the velocity, pressure, and magnetic field, which are the unknowns and $\textbf{u}_t=\fr{\partial\textbf{u}}{\partial t}$, $\textbf{H}_t=\fr{\partial\textbf{H}}{\partial t}$.
	The function $\textbf{f}$ is the known body force and $\textbf{g}$ denotes the known applied current with $\operatorname{div}\textbf{g}=0$. Other parameters in the above equations are the viscosity coefficient $\nu$, the magnetic permeability $\mu$, the electric conductivity $\sigma$, and the density of the fluid is assumed to be 1. For the sake of simplicity, we assume that these constants are independent in the fluid. Additional discussion and derivation of these equations see \cite{1,2}. Typically, the system \refe{1.1} is complemented with the following initial and boundary conditions
	\begin{align}\label{1.2}
	\textbf{u}(x,0)=\textbf{u}_0(x),\quad \textbf{H}(x,0)=\textbf{H}_0(x), \ \ x\in\O,\\
    \label{1.3}
	 \textbf{u}=\textbf{0},\quad \textbf{H}\cdot\textbf{n}=0,\quad\operatorname{curl}\textbf{H}\times\textbf{n}=\textbf{0}, \quad \text{on  }\partial\Omega\times[0,T],
	\end{align}
	 with $\operatorname{div}{\textbf{u}_{0}=0}$, $\operatorname{div}{\textbf{H}_{0}=0}$, where $\partial\O$ is the boundary of $\O$, $\textbf{n}$ is the unit outer normal vector, and the physical interpretations can be seen in \cite{1,2,3}.
	
	 MHD is a subject of researching the interaction of electromagnetic fields and conducting fluids combining hydrodynamic and electrodynamic methods. The MHD modeling often plays a significant role in industrial MHD flows, astronomy, geophysics, and engineering problems, such as metallurgical engineering, electromagnetic pumping, stirring of liquid metals, liquid metal cooling of nuclear reactors, the electromagnetic casting of metals and measuring flow quantities based on induction (see \cite{4,5}).
	 There is usually no analytical solution for the unsteady incompressible MHD system. The fully implicit time-stepping schemes are (almost) unconditionally stable. However, they need to solve a nonlinear system at each step. A fully explicit scheme is simple and easy to carry out on the computer, but it has more restrictions on time steps about stability.
	 There are many efficient semi-implicit time-stepping schemes, which treat implicit schemes for the linear terms and semi-implicit schemes or explicit schemes for the nonlinear terms. The numerical methods can maintain stability with fewer restrictions of the time step, save memory space for large-scale problems and be easy for practical implementation. Many researchers have devoted themselves to nonlinear partial differential equations.
	In \cite{5}, Gunzburger et al. discussed the existence and uniqueness of the solution of a weak formulation of the equations of stationary, incompressible magnetohydrodynamics and proved an optimal error estimate for the approximation solution of the finite element discretization. In \cite{WD}, Wang established the existence and uniqueness of a global solution to the initial-boundary value problem with general large initial data in $H^1$.
	In \cite{9} and \cite{10}, He provided a fully discrete stabilized finite-element method for the two-dimensional time-dependent Navier-Stokes problem and an optimal error order for the velocity and pressure as well as presented stability and error analysis for a spectral Galerkin method for the Navier-Stokes equations with $H^2$ or $H^{1}$ initial data, respectively.
	In \cite{8}, Diening et al. studied a semi-implicit Euler scheme for generalized Newtonian fluids. In \cite{6}, He and Sun showed the stability and convergence of the Crank-Nicolson/Adams-Bashforth scheme for the two-dimensional nonstationary Navier-Stokes equations. In \cite{12}, Prohl verified convergence of iterates of different coupling and decoupling fully discrete schemes for MHD equations and Sovinec as well as King analyzed a mixed semi-implicit/implicit algorithm for low-frequency two-fluid plasma modeling in \cite{13}.
	In \cite{14}, error estimates of a linear decoupled Euler-FEM scheme for a mass diffusion model have been studied. In \cite{15}, Dong and He investigated two-level Newton iterative method for the 2D/3D stationary incompressible magnetohydrodynamics.
	In \cite{3}, He studied unconditional convergence of the Euler semi-implicit scheme for the three-dimensional incompressible MHD equations and Zhang as well as He analyzed unconditional convergence of the Euler semi-implicit scheme for the 3D incompressible MHD equations in \cite{ZH}.
	In \cite{16} and \cite{17}, Si et al. confirmed unconditional stability and error estimates of modified characteristics FEMs for the Navier-Stokes equations and MHD equations, respectively. There are also several applications of other finite element methods for different equations, which can be found in\cite{18,20,19,SWW}.

	In recent years, the implicit-explicit combination of Crank-Nicolson and leapfrog is a popular second-order accurate difference format and is widely used in atmosphere, ocean, and climate codes. In \cite{21}, Li developed a leapfrog mixed finite element method for Maxwell's equations resulting from metamaterials and proved its convergence. In \cite{22}, Huang et al. presented superconvergence analysis for time-dependent Maxwell's equations in metamaterials.
    In \cite{23}, Layton and Trenchea considered a time-stepping CNLF scheme for uncoupling systems of evolution equations and testified that it is conditionally stable. In \cite{24}, Kubacki investigated uncoupling evolutionary groundwater-surface water flows using the CNLF method. In \cite{25}, Jiang et al. promoted a CNLF scheme to solve the geophysical flow, uncoupling groundwater-surface water flow and Stokes flow plus a Coriolis term. What's more, there are also numerous works to dedicate second-order schemes for the non-stationary Navier-Stokes equations.
	In \cite{26}, Heywood and Rannacher deduced $L^{2}$-almost unconditional convergence for the fully implicit Crank-Nicolson scheme. In \cite{27}, Shen presented the error estimates of the projection methods for the Navier-Stokes equations. He and Sun studied the $L^{2}$-almost unconditional convergence for the semi-implicit Crank-Nicolson extrapolation scheme for the Navier-Stokes equations in \cite{28, 29}. Tang and Huang proposed a CNLF scheme for the unsteady incompressible Navier-Stokes equations and proved stability and convergence of the CNLF scheme in \cite{30}. For the unsteady incompressible MHD equations, Yuksel and Ingram investigated the fully implicit Crank-Nicolson discretization for MHD flow at small magnetic Reynolds numbers in \cite{31}. Zhang et al. derived the almost unconditional convergence of the Crank-Nicolson extrapolation time-stepping scheme in \cite{32} and Dong as well as He gave the optimal convergence analysis of the Crank-Nicolson extrapolation scheme for the three-dimensional incompressible MHD equations in \cite{4}. Zhang et al. studied second-order unconditional linear energy stable, rotational velocity correction method for unsteady incompressible magnetohydrodynamics equations in \cite{ZSF22}.

	In this paper, we introduce a fully discrete CNLF semi-implicit scheme based on FEM for numerically solving the unsteady incompressible MHD equations, which deals with spatial discretization by Galerkin finite element approximation.
	For the temporal discretization, we use the Euler-backward implicit/explicit scheme in the first step, apply the Crank-Nicolson extrapolation scheme in the second step and employ the CNLF scheme for the linear terms and the semi-implicit method for nonlinear terms in others. The error estimates between exact solution $(\textbf{u}(t_{n}), p(t_{n}), \textbf{H}(t_{n}))$ and the fully discrete finite element solution $(\textbf{u}_{h}^{n}, p_{h}^{n}, \textbf{H}_{h}^{n})$ of the CNLF scheme can be divided into $(\textbf{u}(t_{n}), p(t_{n}), \textbf{H}(t_{n}))$ and the spatially discrete solution $(\textbf{u}_{h}(t_{n}), p_{h}(t_{n}), \textbf{H}_{h}(t_{n}))$ as well as $(\textbf{u}_{h}(t_{n}), p_{h}(t_{n}), \textbf{H}_{h}(t_{n}))$ and $(\textbf{u}_{h}^{n}, p_{h}^{n}, \textbf{H}_{h}^{n})$.  Meanwhile, the almost unconditional stability and convergence will be presented, and the optimal error estimates will be derived by a new negative norm technique. Other sections of this article are organized as follows. In Section 2, we give the functional setting of the problem \refe{1.1}\--\refe{1.3}, some basic assumptions, subsequent Galerkin FEM, and some smooth properties of finite element solution. In Section 3, we provide the fully discrete CNLF scheme and devote it to the almost unconditional stability of the scheme. Section 4 considers $L^{2}$-optimal error estimates of the method. Section 5 shows some numerical experiments to confirm the validity of the theoretical analysis.
	\section{The functional setting and Galerkin finite element formula}
	\setcounter{equation}{0}

	In order to provide the variation form of the problem \refe{1.1}\--\refe{1.3}, we employ the Sobolev spaces
	\begin{align*}
	&X=H_0^1(\Omega)^d,\ \ H_{0}^{1}(\Omega)=\{v\in H^{1}(\Omega),v|_{\partial\Omega}=0\}, \ \ X_{0}=\{\textbf{v}\in X, \operatorname{div}\textbf{v}=0 \},\\
	&W=H_n^1(\Omega)=\{\textbf{B}\in H^1(\O)^d,\textbf{B}\cdot\textbf{n}|_{\partial\O}=0\},\ \ V_n=\{\textbf{B}\in W,\operatorname{div}\textbf{B}=0\},\\
	&M=L_0^2(\Omega)=\{q\in L^2(\Omega), \int_\Omega qdx=0\},\ \ H=\{\textbf{v}\in L^{2}(\O)^{d}, \operatorname{div}\textbf{v}=0, \textbf{v}\cdot \textbf{n}|_{\partial\Omega}=0\}.
	\end{align*}
	The space $L^{2}$ is equipped with the $L^{2}$-scalar product $(\cdot,\cdot)$ and $L^{2}$-norm $\|\cdot\|_{0}$.
Others are equipped with the norms
	\begin{align*}
	\|\textbf{v}\|_{X}=\|\nabla\textbf{v}\|_{0}, \quad \|q\|_{M}=\|q\|_0,\quad \|\textbf{B}\|_{W}=\|\textbf{B}\|_{1}.
	\end{align*}
	We apply the standard scalar Sobolev space $H^{k}(\Omega)=W^{k,2}(\Omega)$ for nonnegative integer $k$ with norm $\|v\|_{k}=(\sum_{|r|=0}^{k}\|\partial^{r}v\|_{0}^{2})^{1/2}$, where  $\partial^{r}=\frac{\partial^{|r|}}{\partial x_{1}^{r_{1}} \ldots \partial x_{d}^{r_{d}}}$, $\partial$ is the differential operator with respect to $x$ for the multi-index $r=(r_{1}, \ldots,r_{d})$ and $|r|=r_{1}+\ldots+r_{d}$ with $r_{1}, \ldots, r_{d}\geq0$. For vector-value functions, the norm in Sobolev space $\mathbf{H}^{k}(\Omega)=H^{k}(\Omega)^{d}$ is defined by $\|\textbf{v}\|_{k}=(\sum_{i=1}^{d}\|v_{i}\|_{k}^{2})^{1/2}$ (refer to \cite{1,34} for more details). In this article, we assume that $X\times M$ satisfies the following inf-sup condition
	\begin{align*}
		\exists \beta_{0}>0, \ \  \forall q \in M, \ \ \beta_{0}\|q\|_{0} \leq \sup _{\textbf{0}\neq\textbf{v} \in X} \frac{d(\textbf{v}, q)}{\|\nabla \textbf{v}\|_{0}} .
	\end{align*}	
	Denote
	\begin{align*}
	a(\textbf{u}, \textbf{v})&=\nu(\nabla \textbf{u}, \nabla \textbf{v}), \ \  \forall \textbf{u}, \textbf{v} \in X ,\\
	d(\textbf{v}, p)&=(\operatorname{div} \textbf{v},p), \ \ \forall\textbf{v} \in X,\ \ p \in M ,
\end{align*}
and the explicitly skew-symmetric convection term
\begin{align*}
	b(\textbf{u}, \textbf{v}, \textbf{w})&=\int_{\Omega}(\textbf{u} \cdot \nabla) \textbf{v} \cdot \textbf{w} dx+\frac{1}{2} \int_{\Omega}(\operatorname{div} \textbf{u}) \textbf{v} \cdot \textbf{w} d x \nn\\
	&=\frac{1}{2} \int_{\Omega}(\textbf{u} \cdot \nabla) \textbf{v} \cdot \textbf{w} d x-\frac{1}{2}\int_{\Omega}(\textbf{u} \cdot \nabla) \textbf{w} \cdot \textbf{v} d x,  \quad\forall \textbf{u} \in X_{0}, \textbf{v}, \textbf{w} \in X .
	\end{align*}
	The variational formulation of the problem \refe{1.1}\--\refe{1.3} is given as follows: Find
	 $(\textbf{u}, p, \textbf{H}) \in L^{2}(0, T ; X) \times L^{2}(0, T ; M) \times L^{2}(0, T ; W)$ for all $(\textbf{v},q,\textbf{B}) \in X \times M \times V_{n}$ and for all almost $t\in(0,T)$
\begin{align}\label{2.1}
	&(\textbf{u}_{t}, \textbf{v})+ a(\textbf{u}, \textbf{v})-d(\textbf{v}, p)+d(\textbf{u}, q)+ b(\textbf{u}, \textbf{u}, \textbf{v})+\mu(\textbf{H} \times \operatorname{curl} \textbf{H}, \textbf{v})=(\textbf{f}, \textbf{v}),\\
\label{2.2}
	&\mu(\textbf{H}_{t},\textbf{B})+\sigma^{-1}(\operatorname{curl}\textbf{H}, \operatorname{curl}\textbf{B})-\mu(\textbf{u}\times \textbf{H},\operatorname{curl}\textbf{B})=(\textbf{g},\textbf{B}),\\
\label{2.3}
	&\textbf{u}(0)=\textbf{u}_{0}, \quad \textbf{H}(0)=\textbf{H}_{0}.
	\end{align}
Like in \cite{3}, we display the following three assumptions.

\textbf{Assumption(A0)}: The initial data $\textbf{u}_{0}\in X_{0}\cap H^{2}(\Omega)^{d}$, $\textbf{H}_{0}\in V_{n}\cap H^{2}(\Omega)^{d}$, the force $\textbf{f}$ and the applied current $\textbf{g}$ satisfy the following bound
\begin{align*}
\sup_{0\leq t\leq T}\{\|\textbf{f}(t)\|_{0}^{2}+\|\textbf{f}_{t}(t)\|_{0}^{2}+\|\textbf{g}(t)\|_{0}^{2}+\|\textbf{g}_{t}(t)\|_{0}^{2}+\|\textbf{u}_{0}\|_{2}^{2}+\|\textbf{H}_{0}\|_{2}^{2}+\|\textbf{f}_{tt}\|_{1}^{2}+\|\textbf{g}_{tt}\|_{1}^{2}\}
\leq\kappa_{0}.
\end{align*}

\textbf{Assumption(A1)}: The unique solution $(\textbf{u}(t),p(t),\textbf{H}(t))$ of the problem \refe{2.1}\--\refe{2.3} satisfies $\textbf{u}(t)\in L^{2}(0,T;X_{0})$, $\textbf{H}(t)\in L^{2}(0,T;V_{n})$ and $p(t)\in L^{2}(0,T;M)$ such that
\begin{align*}
\int_{0}^{T}(\|\nabla\textbf{u}(t)\|_{0}^{4}+\|\operatorname{curl}\textbf{H}(t)\|_{0}^{4})dt\leq\kappa_{1}.
\end{align*}
 From here to later, $\kappa_{i}$ ($i=0$, $1$, $2$, $\ldots$) is a positive constant depending on \{$\Omega$, $\nu$, $\sigma$, $\mu$, $T$, $\textbf{u}_{0}$, $\textbf{H}_{0}$, $\textbf{f}$, $\textbf{g}$\}.

\textbf{Assumption(A2)}: Assume that the boundary of $\Omega$ is smooth so that the steady Stokes problem
\begin{equation}\label{2.4}
\begin{aligned}
&-\Delta \textbf{v}+\nabla q=\textbf{f}, \ \operatorname{div} \textbf{v}=0 \text { in } \Omega,   \left.\textbf{v}\right|_{\partial \Omega}=\textbf{0},
\end{aligned}
\end{equation}
for prescribed $\textbf{f} \in L^{2}(\Omega)^{d}$ has a unique solution ($\textbf{v},q$) which satisfies
\begin{equation}\label{2.5}
\|\textbf{v}\|_{2}+\|q\|_{1} \leq c_{0}\|\textbf{f}\|_{0};
\end{equation}
	and Maxwell's equations
	\begin{equation}\label{2.6}
	\operatorname{curl}\operatorname{curl}\textbf{H}=\textbf{g}, \   \operatorname{div}\textbf{H}=0 \text { in } \Omega,  \  \textbf{n} \times \operatorname{curl}\textbf{H}=\textbf{0},  \  \textbf{H} \cdot \textbf{n}=0 \text { on } \partial \Omega,
	\end{equation}
	for the prescribed $\textbf{g} \in L^{2}(\Omega)^{d}$ admit a unique solution $\textbf{H} \in V_{n}$ which satisfies
	\begin{equation}\label{2.7}
	\|\textbf{H}\|_{2} \leq c_{0}\|\textbf{g}\|_{0}.
	\end{equation}
Here and after, $c_{i}$ ($i=0$, $1$, $2$, $\ldots$) presents a generic positive constant depending only on $\Omega$.
\begin{remark}\label{r1}
Let $\textbf{B}=\nabla\phi$ with $\phi\in H_{0}^{2}(\Omega)$ in \refe{2.2}, using $\operatorname{curl}\nabla\phi=\textbf{0}$, we can conclude that
\begin{align}\label{2.70}
\mu(\textbf{H}_{t},\nabla \phi)=(\textbf{g},\nabla \phi).
\end{align}
\refe{2.70} employs divergence theorem, Green formula and $\phi\in H_{0}^{2}(\Omega)$. Then, there holds
\begin{align}\label{2.71}
\mu(\operatorname{div}\textbf{H}_{t},\phi)=(\operatorname{div} \textbf{g},\phi).
\end{align}
Combining \refe{2.71} with $\operatorname{div} \textbf{g}=0$ and $\operatorname{div}\textbf{H}_{0}=0$, it yields
$$\operatorname{div}\textbf{H}(t)=0,$$
for all $0\leq t\leq T$.
\end{remark}

\begin{remark}\label{r2}In \cite{3,35}, we can see the validity of \refe{2.5} if $\partial\O$ is of $C^2$, or if $\Omega$ is a convex polyhedron. In \cite{36,37}, the \refe{2.7} is known if $\partial\O$ is of $C^2$. Assumption(A0) guarantees the existence of a unique strong solution to the problem \refe{1.1}\--\refe{1.3} on some interval $[0,T)$ with $\textbf{u}\in C([0,T];X)\cap L^{2}(0,T;H^{2}(\Omega)^{d})$, $\textbf{H}\in C([0,T];W)\cap L^{2}(0,T;H^{2}(\Omega)^{d})$, $p\in L^{2}(0,T;M\cap H^{1}(\Omega))$, and Assumption(A1) ensures the uniqueness of a weak solution to the problem \refe{2.1}\--\refe{2.3} (refer to \cite{36}).
\end{remark}

Let $T_{h}$ be a family of regular and quasi-uniform triangulation partition of $\Omega$, which consists of tetrahedral elements $K$, where $h=\max_{K\in T_{h}}h_{K}$ and $h_{K}$ is the diameter of the element $K$. Each tetrahedron $K$ is supposed to be the image of a reference tetrahedron $\hat{K}$ under an affine map $F_{K}$ and we introduce the following finite element spaces \cite{3,4}.
\begin{align*}
&Y_{h}^{k}=\{\textbf{v}_{h}\in C^{0}(\Omega),\ \ \textbf{v}_{h}(F_{K}(\hat{x}))\in P_{k}(\hat{K}) ,\ \ \forall K\in T_{h} \}.\\
&P_{1,h}^{b}=\{\textbf{v}_{h}\in C^{0}(\Omega),\ \ \textbf{v}_{h}(F_{K}(\hat{x}))\in P_{1}(\hat{K})\oplus \text{span}\{\hat{b}\},\ \ \forall K\in T_{h}\}.\\
&X_{h}=(P_{1,h}^{b})^{d}\cap X, \ \ M_{h}=Y_{h}^{1}\cap M, \ \ W_{h}=(Y_{h}^{1})^{d}\cap W.\\
&X_{0h}=\{\textbf{v}_{h}\in X_{h},\ \  d(\textbf{v}_{h},q_{h})=0,\ \ \forall q_{h}\in M_{h}\},
\end{align*}
    where $P_{k}(\hat{K})$ is the space of polynomials of degree $k$ defined on $\hat{K}$, $Y_{h}^{k}$ is the classical Lagrange finite element space for $k\geq1$, and $\hat{b}\in H_{0}^{1}(\hat{K})$ is called a cubic bubble function, which takes the value 1 at the barycenter of $\hat{K}$ and satisfies $0\leq \hat{b}(\hat{x})\leq1$.
    The mixed finite element method was proposed by Gunzburger et al. in \cite{5}, which admits the classical $H^{1}$-conforming finite elements. The mixed finite element space pair $(X_{h},M_{h})$ use the Mini-element \cite{38} to approximate the velocity and pressure, which satisfies the discrete inf-sup condition. Likely, $W_{h}$ approximates the magnetic field and $X_{0h}$ is the subspace of $X_{h}$.

    The mixed finite element space pair $(X_{h}, M_{h})$ satisfies the so-called inf-sup condition
    \begin{align}\label{2.8}
    \sup_{\textbf{0}\neq \textbf{v}_{h}\in X_{h}}\frac{d(\textbf{v}_{h},q_{h})}{\|\nabla \textbf{v}_{h}\|_{0}}\geq\beta_{1}\|q_{h}\|_{0}, \ \ \forall q_{h}\in M_{h},
    \end{align}
    where $\beta_{1}$ is a positive constant depending on $\Omega$.
    Denote by $P_{h}$: $L^{2}(\Omega)^{d}\rightarrow X_{0h}$, $R_{0h}$: $L^{2}(\Omega)^{d}\rightarrow W_{h}$ the $L^{2}$-orthogonal projection. With the above statements, the standard Galerkin finite element approximation of problem \refe{2.1}\--\refe{2.3} is to find $(\textbf{u}_{h},p_{h},\textbf{H}_{h})\in
    X_{h}\times M_{h} \times W_{h}$, such that for all $(\textbf{v}_{h},q_{h},\textbf{B}_{h})\in X_{h}\times M_{h} \times W_{h}$ and $t\in(0,T)$
    \begin{align}\label{2.9}
    &(\textbf{u}_{ht},\textbf{v}_{h})+a(\textbf{u}_{h},\textbf{v}_{h})-d(\textbf{v}_{h},p_{h})
     +d(\textbf{u}_{h},q_{h})+b(\textbf{u}_{h},\textbf{u}_{h},\textbf{v}_{h})\nn\\
     &\ \ +\mu(\textbf{H}_{h}\times\operatorname{curl}\textbf{H}_{h},
    \textbf{v}_{h})=(\textbf{f},\textbf{v}_{h}),\\
    \label{2.10}
    &\mu(\textbf{H}_{ht},\textbf{B}_{h})+\sigma^{-1}(\operatorname{curl}\textbf{H}_{h},\operatorname{curl}\textbf{B}_{h})+\sigma^{-1}(\operatorname{div}\textbf{H}_{h},\operatorname{div}\textbf{B}_{h})\nn\\
    &\ \ -\mu(\textbf{u}_{h}\times \textbf{H}_{h},\operatorname{curl}\textbf{B}_{h})=(\textbf{g},\textbf{B}_{h}),\\
    \label{2.11}
    &\textbf{u}_{h}(0)=P_{h}\textbf{u}_{0},\ \ \textbf{H}_{h}(0)=R_{0h}\textbf{H}_{0}.
    \end{align}

    Here, the discrete Stokes operator is defined by $A_{1h}=-P_{h}\Delta_{h}$, such that $\Delta_{h}$ (refer to \cite{3,4,28,35,36})
    $$-( \Delta_{h}\textbf{u}_{h},\textbf{v}_{h})=(\nabla \textbf{u}_{h}, \nabla \textbf{v}_{h}), \ \ \forall \textbf{u}_{h}, \textbf{v}_{h}\in X_{h},$$
    and the discrete norm $\|\textbf{v}_{h}\|_{\alpha}=\|A^{\frac{\alpha}{2}}_{1h}\textbf{v}_{h}\|_{0}$ of the $\alpha\in R$ order can be denoted, where
    $$\|\textbf{v}_{h}\|_{1}=\|\nabla \textbf{v}_{h}\|_{0},\ \ \|\textbf{v}_{h}\|_{2}=\|A_{1h}\textbf{v}_{h}\|_{0}, \ \ \|\textbf{v}_{h}\|_{-1}=\|A^{-\frac{1}{2}}_{1h}\textbf{v}_{h}\|_{0}, \ \ \forall \textbf{v}_{h}\in X_{0h}.$$
    Likewise, another discrete operator $A_{2h}\textbf{B}_{h}=R_{0h}(\nabla_{h}\times\operatorname{curl}\mathbf{B}_{h}+\nabla_{h}\operatorname{div}\textbf{B}_{h})\in W_{h}$ (see \cite{3,4}) is defined by
    $$(A_{2h}\textbf{B}_{h},\textbf{C}_{h})=(A^{\frac{1}{2}}_{2h}\textbf{B}_{h},A^{\frac{1}{2}}_{2h}\textbf{C}_{h})=(\operatorname{curl}\textbf{B}_{h},\operatorname{curl}\textbf{C}_{h})+(\operatorname{div}\textbf{B}_{h},\operatorname{div}\textbf{C}_{h}),\ \ \forall\textbf{B}_{h}, \textbf{C}_{h} \in W_{h},$$
    and the discrete norm $\|\textbf{B}_{h}\|_{\alpha}=\|A^{\frac{\alpha}{2}}_{2h}\textbf{B}_{h}\|_{0}$ of the $\alpha\in R$ order can be denoted, where
    $$\|\textbf{B}_{h}\|^{2}_{1}=\|A^{\frac{1}{2}}_{2h}\textbf{B}_{h}\|^{2}_{0}=\|\operatorname{curl}\textbf{B}_{h}\|^{2}_{0}+\|\operatorname{div}\textbf{B}_{h}\|^{2}_{0},$$ $$\|\textbf{B}_{h}\|_{2}=\|\nabla_{h}\times\operatorname{curl}\textbf{B}_{h}+\nabla_{h}\operatorname{div}\textbf{B}_{h}\|_{0},\ \  \|\textbf{B}_{h}\|_{-1}=\sup_{\textbf{0}\neq\textbf{C}_{h}\in W_{h}}\frac{(\textbf{B}_{h},\textbf{C}_{h})}{\|A^{\frac{1}{2}}_{2h}\textbf{C}_{h}\|_{0}},$$
    which can be found in \cite{3,4}.

    Furthermore, we need some discrete estimates \cite{3,4,35,36} and some estimates of the trilinear form $b(\cdot,\cdot,\cdot)$ \cite{3,4,6,sj}
    \begin{align}\label{2.12}
    &\|\nabla\textbf{v}_{h}\|_{L^{3}}+\|\textbf{v}_{h}\|_{L^{\infty}}\leq c_{1}\|\nabla\textbf{v}_{h}\|^{\frac{1}{2}}_{0}\|A_{1h}\textbf{v}_{h}\|^{\frac{1}{2}}_{0},\ \ \|\nabla\textbf{v}_{h}\|_{0}+\|\nabla\textbf{v}_{h}\|_{L^{6}}\leq c_{1}\|A_{1h}\textbf{v}_{h}\|_{0},\\
    \label{2.13}
    &\|\nabla\textbf{B}_{h}\|_{L^{3}}+\|\textbf{B}_{h}\|_{L^{\infty}}\leq c_{1}\|\nabla\textbf{B}_{h}\|^{\frac{1}{2}}_{0}\|A_{2h}\textbf{B}_{h}\|^{\frac{1}{2}}_{0},\ \ \|\nabla\textbf{B}_{h}\|_{0}+\|\nabla\textbf{B}_{h}\|_{L^{6}}\leq c_{1}\|A_{2h}\textbf{B}_{h}\|_{0},\\
    \label{2.14'}
    &b(\textbf{u}_{h},\textbf{v}_{h},\textbf{w}_{h})=-b(\textbf{u}_{h},\textbf{w}_{h},\textbf{v}_{h}),\\
    \label{2.14}
    &|b(\textbf{u}_{h},\textbf{v}_{h},\textbf{w}_{h})|+|b(\textbf{w}_{h},\textbf{u}_{h},\textbf{v}_{h})|+|b(\textbf{v}_{h},\textbf{u}_{h},\textbf{w}_{h})|\leq c_{2}\|\textbf{u}_{h}\|^{\frac{1}{2}}_{1}\|A_{1h}\textbf{u}_{h}\|^{\frac{1}{2}}_{0}\|\textbf{v}_{h}\|_{1}\|\textbf{w}_{h}\|_{0},\\
    \label{2.15}
    &|(\textbf{u}_{h}\times \textbf{H}_{h},\operatorname{curl}\textbf{B}_{h})|+|(\textbf{u}_{h}\times \textbf{B}_{h},\operatorname{curl}\textbf{H}_{h})|\leq c_{3}\|\textbf{B}_{h}\|^{\frac{1}{2}}_{1}\|A_{2h}\textbf{B}_{h}\|^{\frac{1}{2}}_{0}\|\textbf{H}_{h}\|_{1}\|\textbf{u}_{h}\|_{0},\\
    \label{2.16}
    &|(\textbf{u}_{h}\times \textbf{H}_{h},\operatorname{curl}\textbf{B}_{h})|+|(\textbf{u}_{h}\times \textbf{B}_{h},\operatorname{curl}\textbf{H}_{h})|\leq c_{3}\|\textbf{u}_{h}\|_{1}^{\frac{1}{2}}\|A_{1h}\textbf{u}_{h}\|_{0}^{\frac{1}{2}}\|\textbf{B}_{h}\|_{1}\|\textbf{H}_{h}\|_{0},\\
    \label{2.17}
    &|b(\textbf{u}_{h},\textbf{v}_{h},\textbf{w}_{h})|+|b(\textbf{w}_{h},\textbf{u}_{h},\textbf{v}_{h})|+|b(\textbf{v}_{h},\textbf{u}_{h},\textbf{w}_{h})|\nn\\
    &\ \ \leq\frac{c_{4}}{2}(\|A_{1h}\textbf{u}_{h}\|_{0}\|\textbf{v}_{h}\|^{\frac{1}{2}}_{1}\|A_{1h}\textbf{v}_{h}\|_{0}^{\frac{1}{2}}+\|A_{1h}\textbf{v}_{h}\|_{0}
    \|\textbf{u}_{h}\|_{1}^{\frac{1}{2}}\|A_{1h}\textbf{u}_{h}\|^{\frac{1}{2}}_{0})\|\textbf{w}_{h}\|_{-1},\\
    \label{2.18}
    &|(\textbf{u}_{h}\times \textbf{H}_{h},\operatorname{curl}\textbf{B}_{h})|+|(\textbf{u}_{h}\times \textbf{B}_{h},\operatorname{curl}\textbf{H}_{h})|\nn\\
    &\ \ \leq\frac{c_{5}}{2}(\|A_{1h}\textbf{u}_{h}\|_{0}\|\textbf{H}_{h}\|^{\frac{1}{2}}_{1}\|A_{2h}\textbf{H}_{h}\|^{\frac{1}{2}}_{0}
    +\|A_{2h}\textbf{H}_{h}\|_{0}\|\textbf{u}_{h}\|^{\frac{1}{2}}_{1}\|A_{1h}\textbf{u}_{h}\|^{\frac{1}{2}}_{0})\|\textbf{B}_{h}\|_{-1},\\
    \label{2.19}
    &|(\textbf{u}_{h}\times \textbf{H}_{h},\operatorname{curl}\textbf{B}_{h})|+|(\textbf{u}_{h}\times \textbf{B}_{h},\operatorname{curl}\textbf{H}_{h})|\nn\\
    &\ \ \leq\frac{c_{5}}{2}(\|A_{2h}\textbf{H}_{h}\|_{0}\|\textbf{B}_{h}\|^{\frac{1}{2}}_{1}\|A_{2h}\textbf{B}_{h}\|^{\frac{1}{2}}_{0}
    +\|A_{2h}\textbf{B}_{h}\|_{0}\|\textbf{H}_{h}\|^{\frac{1}{2}}_{1}\|A_{2h}\textbf{H}_{h}\|^{\frac{1}{2}}_{0})\|\textbf{u}_{h}\|_{-1},
    \end{align}
    for all $\textbf{u}_{h}, \textbf{v}_{h}, \textbf{w}_{h} \in X_{0h}$, $\textbf{B}_{h}, \textbf{H}_{h} \in W_{h}$.

    Moreover, we recall the following smooth properties of the solution ($\textbf{u}_{h}(t),\textbf{p}_{h}(t),\textbf{B}_{h}(t)$) of scheme \refe{2.9}\--\refe{2.11}, which are deduced in \cite{41,40}.
    \begin{lemma} \label{L2.1}
	Assume that Assumptions(A0-A2) are satisfied. Then the following smooth properties hold
\begin{align}\label{2.20}
&\|\textbf{u}_{h}(t)\|^{2}_{2}+\|\textbf{H}_{h}(t)\|^{2}_{2}\leq \kappa_{2},\ \  \tau^{r}(t)(\|\textbf{u}_{ht}(t)\|^{2}_{r}+\|\textbf{H}_{ht}(t)\|^{2}_{r})\leq \kappa_{2},\ \ r=0,\ 1,\ 2,\\
\label{2.21}
&\tau^{r+2}(t)(\|\textbf{u}_{htt}(t)\|^{2}_{r}+\|\textbf{H}_{htt}(t)\|^{2}_{r})\leq \kappa_{2},\ \ r=-1,\ 0,\ 1,\ \ \tau^{2}(t)\|p_{h}(t)\|_{0}^{2}\leq\kappa_2.\\
\label{2.22}
&\int_{0}^{t}{\tau^{r}(s)(\|\textbf{u}_{ht}(s)\|_{r+1}^{2}+\|\textbf{H}_{ht}(s)\|_{r+1}^{2})}ds\leq \kappa_{2},\ \ r=0,\ 1,\\
\label{2.23}
&\int_{0}^{t}{\tau^{r+1}(s)(\|\textbf{u}_{htt}(s)\|_{r}^{2}+\|\textbf{H}_{htt}(s)\|_{r}^{2})}ds\leq \kappa_{2},\ \ r=-1,\ 0,\ 1,\\
\label{2.24}
&\int_{0}^{t}{\tau^{r+2}(s)(\|\textbf{u}_{httt}(s)\|_{r-1}^{2}+\|\textbf{H}_{httt}(s)\|_{r-1}^{2})}ds\leq \kappa_{2},\ \ r=-1,\ 0,\ 1,\\
\label{2.25}
&\int_{0}^{t}({\tau^{3}(s)(\|A_{1h}\textbf{u}_{htt}(s)\|_{0}^{2}+\|A_{2h}\textbf{H}_{htt}(s)\|_{0}^{2}+\|p_{htt}(s)\|_{0}^{2})}+\tau(s)\|p_{ht}(s)\|_{0}^{2})ds\leq \kappa_{2},
\end{align}
for all $t\in[0,T]$, where $\tau(t)=\min\{1,t\}$.
	\end{lemma}
Next, we need to derive the following Theorem \ref{T2.1}.
\begin{theorem}\label{T2.1}
Suppose that Assumptions (A0-A2) hold, $\textbf{u}_h$ and $\textbf{H}_h$ are the solutions of  \refe{2.9} and \refe{2.10}, respectively. Then there holds that
\begin{align}\label{2.26}
\tau^{4}(t)(\|A_{1h}\textbf{u}_{htt}\|_{0}^{2}+\|A_{2h}\textbf{H}_{htt}\|_{0}^{2})+\int_{0}^{t}{\tau^{4}(s)(\|A_{1h}^{\frac{3}{2}}\textbf{u}_{htt}(s)\|_{0}^{2}+\|A_{2h}^{\frac{3}{2}}\textbf{H}_{htt}(s)\|_{0}^{2})}ds\leq\kappa_{3},
\end{align}
for all $t\in(0,T]$.

\begin{proof}
Calculating the second-order derivative of \refe{2.9} and \refe{2.10} with respect to $t$, it yields
\begin{align}\label{2.27}
&(\textbf{u}_{httt},\textbf{v}_{h})+a(\textbf{u}_{htt},\textbf{v}_{h})-d(\textbf{v}_{h},p_{htt})+d(\textbf{u}_{htt},q_{h})=(\textbf{f}_{tt},\textbf{v}_{h})
-b_{tt}(\textbf{u}_{h},\textbf{u}_{h},\textbf{v}_{h})\nn\\
&\ \  -\mu(\textbf{H}_{h}\times\operatorname{curl}\textbf{H}_{h},\textbf{v}_{h})_{tt},\\
\label{2.28}
&\mu(\textbf{H}_{httt},\textbf{B}_{h})+\sigma^{-1}(\operatorname{curl}\textbf{H}_{htt},\operatorname{curl}\textbf{B}_{h})+\sigma^{-1}(\operatorname{div}\textbf{H}_{htt},\operatorname{div}\textbf{B}_{h})=
(\textbf{g}_{tt},\textbf{B}_{h})+\mu(\textbf{u}_{h}\times\textbf{H}_{h},\operatorname{curl}\textbf{B}_{h})_{tt},
\end{align}
where
\begin{align}\label{2.29}
&b_{tt}(\textbf{u}_{h},\textbf{u}_{h},\textbf{v}_{h})=b(\textbf{u}_{htt},\textbf{u}_{h},\textbf{v}_{h})+2b(\textbf{u}_{ht},\textbf{u}_{ht},\textbf{v}_{h})+b(\textbf{u}_{h},\textbf{u}_{htt},\textbf{v}_{h}),\\
\label{2.30}
&(\textbf{H}_{h}\times\operatorname{curl}\textbf{H}_{h},\textbf{v}_{h})_{tt}=(\textbf{H}_{htt}\times\operatorname{curl}\textbf{H}_{h},\textbf{v}_{h})+2(\textbf{H}_{ht}\times\operatorname{curl}\textbf{H}_{ht},\textbf{v}_{h})+(\textbf{H}_{h}\times\operatorname{curl}\textbf{H}_{htt},\textbf{v}_{h}),\\
\label{2.31}
&(\textbf{u}_{h}\times\textbf{H}_{h},\operatorname{curl}\textbf{B}_{h})_{tt}=(\textbf{u}_{htt}\times\textbf{H}_{h},\operatorname{curl}\textbf{B}_{h})+2(\textbf{u}_{ht}\times\textbf{H}_{ht},\operatorname{curl}\textbf{B}_{h})
+(\textbf{u}_{h}\times\textbf{H}_{htt},\operatorname{curl}\textbf{B}_{h}).
\end{align}
Taking $\textbf{v}_{h}=A_{1h}^{2}\textbf{u}_{htt}\in X_{0h}$, $q_{h}=0$ in \refe{2.27} and $\textbf{B}_{h}=A_{2h}^{2}\textbf{H}_{htt}\in W_{h}$ in \refe{2.28}, we can acquire
\begin{align}\label{2.32}
&\frac{1}{2}\frac{d}{dt}\|A_{1h}\textbf{u}_{htt}\|_{0}^{2}+\nu\|A_{1h}^{\frac{3}{2}}\textbf{u}_{htt}\|_{0}^{2}=(\textbf{f}_{tt},A_{1h}^{2}\textbf{u}_{htt})-b_{tt}(\textbf{u}_{h},\textbf{u}_{h},A_{1h}^{2}\textbf{u}_{htt})-\mu(\textbf{H}_{h}\times\operatorname{curl}\textbf{H}_{h},A_{1h}^{2}\textbf{u}_{htt})_{tt},\\
\label{2.33}
&\frac{\mu}{2}\frac{d}{dt}\|A_{2h}\textbf{H}_{htt}\|_{0}^{2}+\sigma^{-1}\|A_{2h}^{\frac{3}{2}}\textbf{H}_{htt}\|_{0}^{2}=(\textbf{g}_{tt},A_{2h}^{2}\textbf{H}_{htt})+\mu(\textbf{u}_{h}\times\textbf{H}_{h},\operatorname{curl}A_{2h}^{2}\textbf{H}_{htt})_{tt}.
\end{align}
Using Young's inequality, the following estimate holds
\begin{align*}
&|(\textbf{f}_{tt},A_{1h}^{2}\textbf{u}_{htt})|\leq \|A_{1h}^{2}\textbf{u}_{htt}\|_{-1}\|\textbf{f}_{tt}\|_{1}\leq \frac{\nu}{10}\|A_{1h}^{\frac{3}{2}}\textbf{u}_{htt}\|_{0}^{2}+\frac{5}{2}\nu^{-1}\|\textbf{f}_{tt}\|_{1}^{2}.
\end{align*}
Taking advantage of \refe{2.12}, \refe{2.17}, \refe{2.29} and Young's inequality leads to
\begin{align*}
|b_{tt}(\textbf{u}_{h},\textbf{u}_{h},A_{1h}^{2}\textbf{u}_{htt})|
&= |b(\textbf{u}_{htt},\textbf{u}_{h},A_{1h}^{2}\textbf{u}_{htt})+2b(\textbf{u}_{ht},\textbf{u}_{ht},A_{1h}^{2}\textbf{u}_{htt})+b(\textbf{u}_{h},\textbf{u}_{htt},A_{1h}^{2}\textbf{u}_{htt})|\nn\\
&\leq 2c_{4}c_{1}^{\frac{1}{2}}\|A_{1h}^{\frac{3}{2}}\textbf{u}_{htt}\|_{0}\|A_{1h}\textbf{u}_{h}\|_{0}\|A_{1h}\textbf{u}_{htt}\|_{0}+2c_{4}c_{1}^{\frac{1}{2}}\|A_{1h}^{\frac{3}{2}}\textbf{u}_{htt}\|_{0}\|A_{1h}\textbf{u}_{ht}\|_{0}\|A_{1h}\textbf{u}_{ht}\|_{0}\nn\\
& \leq \frac{\nu}{5}\|A_{1h}^{\frac{3}{2}}\textbf{u}_{htt}\|_{0}^{2}+10\nu^{-1}c_{4}^{2}c_{1}(\|A_{1h}\textbf{u}_{h}\|_{0}^{2}\|A_{1h}\textbf{u}_{htt}\|_{0}^{2}+\|A_{1h}\textbf{u}_{ht}\|_{0}^{2}\|A_{1h}\textbf{u}_{ht}\|_{0}^{2}).
\end{align*}
Employing \refe{2.13}, \refe{2.19}, \refe{2.30} and Young's inequality yields
\begin{align*}
&\mu|(\textbf{H}_{h}\times\operatorname{curl}\textbf{H}_{h},A_{1h}^{2}\textbf{u}_{htt})_{tt}|\nn\\
&\leq\mu|(\textbf{H}_{htt}\times\operatorname{curl}\textbf{H}_{h},A_{1h}^{2}\textbf{u}_{htt})+(\textbf{H}_{h}\times\operatorname{curl}\textbf{H}_{htt},A_{1h}^{2}\textbf{u}_{htt})|
 +2\mu|(\textbf{u}_{ht}\times\operatorname{curl}\textbf{H}_{ht},A_{1h}^{2}\textbf{u}_{htt})|\\
& \leq 2c_{5}\mu c_{1}^{\frac{1}{2}}\|A_{1h}^{\frac{3}{2}}\textbf{u}_{htt}\|_{0}\|A_{2h}\textbf{H}_{htt}\|_{0}\|A_{2h}\textbf{H}_{h}\|_{0}+2c_{5}\mu c_{1}^{\frac{1}{2}}\|A_{1h}^{\frac{3}{2}}\textbf{u}_{htt}\|_{0}\|A_{2h}\textbf{H}_{ht}\|_{0}\|A_{2h}\textbf{H}_{ht}\|_{0}\nn\\
& \leq \frac{\nu}{5}\|A_{1h}^{\frac{3}{2}}\textbf{u}_{htt}\|_{0}^{2}+10\nu^{-1}\mu^{2} c_{5}^{2}c_{1}(\|A_{2h}\textbf{H}_{htt}\|_{0}^{2}\|A_{2h}\textbf{H}_{h}\|_{0}^{2}+\|A_{2h}\textbf{H}_{ht}\|_{0}^{2}\|A_{2h}\textbf{H}_{ht}\|_{0}^{2}).
\end{align*}
Similarly, making use of \refe{2.12}, \refe{2.13}, \refe{2.18}, \refe{2.31} and Young's inequality, the following estimates exist
\begin{align*}
&|(\textbf{g}_{tt},A_{2h}^{2}\textbf{H}_{htt})|\leq\|A_{2h}^{2}\textbf{H}_{htt}\|_{-1}\|\textbf{g}_{tt}\|_{1}\leq \frac{\sigma^{-1}}{6}\|A_{2h}^{\frac{3}{2}}\textbf{H}_{htt}\|_{0}^{2}+\frac{3}{2}\sigma\|\textbf{g}_{tt}\|_{1}^{2},\\
&\mu|(\textbf{u}_{h}\times\textbf{H}_{h},\operatorname{curl}A_{2h}^{2}\textbf{H}_{htt})_{tt}|\nn\\
&\leq\mu|(\textbf{u}_{htt}\times\textbf{H}_{h},\operatorname{curl}A_{2h}^{2}\textbf{H}_{htt})+(\textbf{u}_{h}\times\textbf{H}_{htt},\operatorname{curl}A_{2h}^{2}\textbf{H}_{htt})|
  +2\mu|(\textbf{u}_{ht}\times\textbf{H}_{ht},\operatorname{curl}A_{2h}^{2}\textbf{H}_{htt})|\nn\\
&\leq \mu c_{5}c_{1}^{\frac{1}{2}}\|A_{2h}^{\frac{3}{2}}\textbf{H}_{htt}\|_{0}\|A_{1h}\textbf{u}_{htt}\|_{0}\|A_{2h}\textbf{H}_{h}\|_{0}+\mu c_{5}c_{1}^{\frac{1}{2}}\|A_{2h}^{\frac{3}{2}}\textbf{H}_{htt}\|_{0}\|A_{1h}\textbf{u}_{h}\|_{0}\|A_{2h}\textbf{H}_{htt}\|_{0}\\
&\ \ \ +2\mu c_{5}c_{1}^{\frac{1}{2}}\|A_{2h}^{\frac{3}{2}}\textbf{H}_{htt}\|_{0}\|A_{1h}\textbf{u}_{ht}\|_{0}\|A_{2h}\textbf{H}_{ht}\|_{0}\\
& \leq\frac{\sigma^{-1}}{3}\|A_{2h}^{\frac{3}{2}}\textbf{H}_{htt}\|_{0}^{2}+3\mu^{2}\sigma c_{5}^{2}c_{1}(\|A_{1h}\textbf{u}_{htt}\|_{0}^{2}\|A_{2h}\textbf{H}_{h}\|_{0}^{2}+\|A_{1h}\textbf{u}_{h}\|_{0}^{2}\|A_{2h}\textbf{H}_{htt}\|_{0}^{2})\\
&\ \ \ +6\mu^{2}\sigma c_{5}^{2}c_{1}\|A_{1h}\textbf{u}_{ht}\|_{0}^{2}\|A_{2h}\textbf{H}_{ht}\|_{0}^{2}.
\end{align*}
Taking the sum of \refe{2.32} and \refe{2.33}, multiplying it by $\tau^{4}(t)$, and substituting the above estimates into the equality, we conclude that
\begin{align}\label{2.34}
&\tau^{4}(t)(\nu\|A_{1h}^{\frac{3}{2}}\textbf{u}_{htt}\|_{0}^{2}+\sigma^{-1}\|A_{2h}^{\frac{3}{2}}\textbf{H}_{htt}\|_{0}^{2})+\frac{d}{dt}(\tau^{4}(t)(\|A_{1h}\textbf{u}_{htt}\|_{0}^{2}+\mu\|A_{2h}\textbf{H}_{htt}\|_{0}^{2}))\nn\\
&\ \ -4\tau^{3}(t)(\|A_{1h}\textbf{u}_{htt}\|_{0}^{2} +\mu\|A_{2h}\textbf{H}_{htt}\|_{0}^{2})\nn\\
&\ \ \leq \tau^{4}(t) [20\nu^{-1}c_{4}^{2}c_{1}(\|A_{1h}\textbf{u}_{h}\|_{0}^{2}\|A_{1h}\textbf{u}_{htt}\|_{0}^{2}+\|A_{1h}\textbf{u}_{ht}\|_{0}^{2}\|A_{1h}\textbf{u}_{ht}\|_{0}^{2}) \nn\\
&\ \ \ \ \
+20\nu^{-1}\mu^{2} c_{5}^{2}c_{1}(\|A_{2h}\textbf{H}_{htt}\|_{0}^{2}\|A_{2h}\textbf{H}_{h}\|_{0}^{2}+\|A_{2h}\textbf{H}_{ht}\|_{0}^{2}\|A_{2h}\textbf{H}_{ht}\|_{0}^{2})\nn\\
&\ \ \ \ \ +6\mu^{2}\sigma c_{5}^{2}c_{1}(\|A_{1h}\textbf{u}_{htt}\|_{0}^{2}\|A_{2h}\textbf{H}_{h}\|_{0}^{2}+\|A_{1h}\textbf{u}_{h}\|_{0}^{2}\|A_{2h}\textbf{H}_{htt}\|_{0}^{2})\nn\\
&\ \ \ \ \ +12\mu^{2}\sigma c_{5}^{2}c_{1}\|A_{1h}\textbf{u}_{ht}\|_{0}^{2}\|A_{2h}\textbf{H}_{ht}\|_{0}^{2}+5\nu^{-1}\|\textbf{f}_{tt}\|_{1}^{2} +3\sigma\|\textbf{g}_{tt}\|_{1}^{2}].
\end{align}
Integrating \refe{2.34} from $0$ to $t$, and taking full advantage of Lemma \ref{L2.1}, \refe{2.26} can be derived.
\end{proof}
\end{theorem}

The optimal error estimates of finite element solution to \refe{2.9}\--\refe{2.11} are given as follows in \cite{40}.
\begin{theorem}\label{T2.2}
Suppose that Assumptions (A0-A2) hold, the solution $(\textbf{u}_{h}(t),p_{h}(t),\textbf{H}_{h}(t))$ of \refe{2.9}\--\refe{2.11} satisfies
\begin{align}\label{2.35}
&\tau(t)(\|\nabla(\textbf{u}(t)-\textbf{u}_{h}(t))\|_{0}^{2}+\|\nabla(\textbf{H}(t)-\textbf{H}_{h}(t))\|_{0}^{2})\leq \kappa_{4}h^{4},\\
\label{2.36}
&\tau(t)(\|\textbf{u}(t)-\textbf{u}_{h}(t)\|_{0}^{2}+\|\textbf{H}(t)-\textbf{H}_{h}(t)\|_{0}^{2})+\tau^2(t)\|p(t)-p_{h}(t)\|_{0}^{2}\leq \kappa_{4}h^{4},
\end{align}
for all $t\in[0,T]$.
\end{theorem}

Then, we provide the discrete Gronwall's inequality in \cite{10, 26, 30}.
	\begin{lemma} \label{L2.3}
	Let $C_{0}, a_n,b_n, d_{n}$ for integers $n_{0}\leq n\leq m$ be non-negative numbers such that
	\begin{align}\label{2.39}
	a_m+\Delta t\sum_{n=k}^{m} b_n\leq
	\Delta t\sum_{n=k}^{m}d_{n} a_n+C_{0}, \ \ n_{0}\leq k \leq m.
	\end{align}
	Suppose that $\Delta t d_n<1$ and set
	$\gamma_n=(1-\Delta t d_n)^{-1}$ for all $n_{0}\leq n\leq m$. Then
    \begin{align}\label{2.40}
	a_m+\Delta t\sum_{n=k}^{m}b_n\leq
	C_{0}\exp(\Delta t\sum_{n=k}^{m}\gamma_n d_n), \ \ n_{0}\leq k\leq m.
	\end{align}
	\end{lemma}
	\section{The fully discrete CNLF scheme and stability of the scheme}
	\setcounter{equation}{0}
	In this section, we consider a fully discrete CNLF scheme for solving the problem \refe{1.1}\--\refe{1.3}. Let $0=t_0<t_1<\ldots <t_N=T$ be uniform partition of the time interval $[0,T]$ with time step $\Delta t=\frac{T}{N}$ and $t_n=n\Delta t$, we use $(\textbf{u}_{h}^{n}, p_{h}^{n}, \textbf{H}_{h}^{n})$ to approximate $(\textbf{u}_{h}(t_{n}),p_{h}(t_{n}),\textbf{H}_{h}(t_{n}))$ at $t=t_{n}$ for $n=0, 1, \ldots, N$.
	The error estimates between exact solution $(\textbf{u}(t_{n}), p(t_{n}), \textbf{H}(t_{n}))$ to \refe{2.1}\--\refe{2.3} and the fully discrete finite element solution $(\textbf{u}_{h}^{n}, p_{h}^{n}, \textbf{H}_{h}^{n})$ of the CNLF scheme can be divided into two parts. The first part is  $(\textbf{u}(t_{n}), p(t_{n}), \textbf{H}(t_{n}))$ and the Galerkin finite element approximation solution $(\textbf{u}_{h}(t_{n}), p_{h}(t_{n}), \textbf{H}_{h}(t_{n}))$ to \refe{2.9}\--\refe{2.11}. The second part is $(\textbf{u}_{h}(t_{n}), p_{h}(t_{n}), \textbf{H}_{h}(t_{n}))$ and $(\textbf{u}_{h}^{n}, p_{h}^{n}, \textbf{H}_{h}^{n})$.
	Thus, we mainly analyze the error estimates of $(\textbf{u}_{h}(t_{n}), p_{h}(t_{n}), \textbf{H}_{h}(t_{n}))$ and $(\textbf{u}_{h}^{n}, p_{h}^{n}, \textbf{H}_{h}^{n})$ due to Theorem \ref{T2.2}.
We adopt the following notations for any function sequence $\{\phi_{h}^n\}$
\begin{align*}
&d_{t}\phi_{h}^{n}=\frac{\phi_{h}^{n+1}-\phi_{h}^{n-1}}{2\Delta t},\ \ \bar{\phi}_{h}^{n}=\frac{1}{2\Delta t}\int_{t_{n-1}}^{t_{n+1}}{\phi_{h}(\cdot,t)}dt,\\ &\tilde{\phi}_{h}^{n}=\frac{\phi_{h}^{n+1}+\phi_{h}^{n-1}}{2},\ \ \tilde{\phi}_{h}(t_{n})=\frac{\phi_{h}(t_{n+1})+\phi_{h}(t_{n-1})}{2}, \ \ \forall n\geq1.
\end{align*}

\begin{algorithm}(The fully discrete CNLF scheme for the MHD equations)
	
	Find ${(\textbf{u}_{h}^{n+1}, p_{h}^{n}, \textbf{H}_{h}^{n+1})}\in X_{h}\times M_{h}\times W_{h}, n=2, 3, \ldots, N-1$ such that, for all $\textbf{v}_{h}\in X_{h}, q_{h}\in M_{h}, \textbf{B}_{h}\in W_{h},$
\begin{align}\label{3.1}
&(d_{t}\textbf{u}_{h}^{n},\textbf{v}_{h})+a(\tilde{\textbf{u}}_{h}^{n},\textbf{v}_{h})-d(\textbf{v}_{h},p_{h}^{n})+d(\tilde{\textbf{u}}_{h}^{n},q_{h})
+b(\textbf{u}_{h}^{n},\tilde{\textbf{u}}_{h}^{n},\textbf{v}_{h})\nn\\
&\ \ +\mu(\textbf{H}_{h}^{n}\times\operatorname{curl}\tilde{\textbf{H}}_{h}^{n},\textbf{v}_{h})=(\textbf{f}(t_{n}),\textbf{v}_{h}),\\
\label{3.2}
&\mu(d_{t}\textbf{H}_{h}^{n},\textbf{B}_{h})+\sigma^{-1}(\operatorname{curl}\tilde{\textbf{H}}_{h}^{n},\operatorname{curl}\textbf{B}_{h})+\sigma^{-1}(\operatorname{div}\tilde{\textbf{H}}_{h}^{n},\operatorname{div}\textbf{B}_{h})\nn\\
&\ \ -\mu(\tilde{\textbf{u}}_{h}^{n}\times\textbf{H}_{h}^{n},\operatorname{curl}\textbf{B}_{h})=(\textbf{g}(t_{n}),\textbf{B}_{h}).
\end{align}
	\end{algorithm}
Here, we define the initial approximating data $(\textbf{u}_{h}^{0},\textbf{H}_{h}^{0})=(P_{h}\textbf{u}_{0},R_{0h}\textbf{H}_{0})$, $(\textbf{u}_{h}^{1},p_{h}^{1},\textbf{H}_{h}^{1})$ is defined by the Euler-backward implicit/explicit scheme and $(\textbf{u}_{h}^{2},p_{h}^{2},\textbf{H}_{h}^{2})$ is denoted by the Crank-Nicolson extrapolation scheme
\begin{align}\label{3.3}
&\left(\frac{\textbf{u}_{h}^{1}-\textbf{u}_{h}^{0}}{\Delta t},\textbf{v}_{h}\right)+a(\textbf{u}_{h}^{1},\textbf{v}_{h})-d(\textbf{v}_{h},p_{h}^{1})+d(\textbf{u}_{h}^{1},q_{h})+b(\textbf{u}_{h}^{0},\textbf{u}_{h}^{0},\textbf{v}_{h})\nn\\
&\ \ +\mu(\textbf{H}_{h}^{0}\times\operatorname{curl}\textbf{H}_{h}^{0},\textbf{v}_{h})=(\textbf{f}(t_{1}),\textbf{v}_{h}), \ \ \forall \textbf{v}_{h}\in X_{h}, q_{h}\in M_{h},\\
\label{3.4}
&\mu\left(\frac{\textbf{H}_{h}^{1}-\textbf{H}_{h}^{0}}{\Delta t},\textbf{B}_{h}\right)+\sigma^{-1}(\operatorname{curl}\textbf{H}_{h}^{1},\operatorname{curl}\textbf{B}_{h})+\sigma^{-1}(\operatorname{div}\textbf{H}_{h}^{1},\operatorname{div}\textbf{B}_{h})\nn\\
&\ \ -\mu(\textbf{u}_{h}^{0}\times \textbf{H}_{h}^{0},\operatorname{curl}\textbf{B}_{h})=(\textbf{g}(t_{1}),\textbf{B}_{h}), \ \ \forall \textbf{B}_{h}\in W_{h}.\\
\label{3.5}
&\left(\frac{\textbf{u}_{h}^{2}-\textbf{u}_{h}^{1}}{\Delta t},\textbf{v}_{h}\right)+a\left(\frac{\textbf{u}_{h}^{2}+\textbf{u}_{h}^{1}}{2},\textbf{v}_{h}\right)-d(\textbf{v}_{h},p_{h}^{2})+d(\textbf{u}_{h}^{2},q_{h})+b\left(\frac{3}{2}\textbf{u}_{h}^{1}-\frac{1}{2}\textbf{u}_{h}^{0}, \frac{\textbf{u}_{h}^{2}+\textbf{u}_{h}^{1}}{2}, \textbf{v}_{h}\right)\nn\\
&\ \ +\mu\left(\left(\frac{3}{2}\textbf{H}_{h}^{1}-\frac{1}{2}\textbf{H}_{h}^{0}\right)\times\operatorname{curl}\frac{\textbf{H}_{h}^{2}+\textbf{H}_{h}^{1}}{2},\textbf{v}_{h}\right)=\left(\frac{\textbf{f}(t_{2})+\textbf{f}(t_{1})}{2}, \textbf{v}_{h}\right),
\\
\label{3.6}
&\mu\left(\frac{\textbf{H}_{h}^{2}-\textbf{H}_{h}^{1}}{\Delta t},\textbf{B}_{h}\right)+\sigma^{-1}\left(\operatorname{curl}\frac{\textbf{H}_{h}^{2}+\textbf{H}_{h}^{1}}{2},\operatorname{curl}\textbf{B}_{h}\right)
+\sigma^{-1}\left(\operatorname{div}\frac{\textbf{H}_{h}^{2}+\textbf{H}_{h}^{1}}{2},\operatorname{div}\textbf{B}_{h}\right)\nn\\
&\ \ -\mu\left(\frac{\textbf{u}_{h}^{2}+\textbf{u}_{h}^{1}}{2}\times\left(\frac{3}{2}\textbf{H}_{h}^{1}-\frac{1}{2}\textbf{H}_{h}^{0}\right),\operatorname{curl}\textbf{B}_{h}\right)=\left(\frac{\textbf{g}(t_{2})+\textbf{g}(t_{1})}{2}, \textbf{B}_{h}\right).
\end{align}

The following stability and convergence have been given in \cite{4, 41}.
\begin{lemma}\label{L3.1}
Assume that Assumptions (A0-A2) hold and $\Delta t\leq \kappa_{5}$, then
\begin{align*}
&\|\textbf{u}_{h}^{j}\|_{r}^{2}+\|\textbf{H}_{h}^{j}\|_{r}^{2}\leq\kappa_{6},\ \ r=0, 1, 2,\\
 &\|\textbf{u}_{h}(t_{j})-\textbf{u}_{h}^{j}\|_{r}^{2}+\|\textbf{H}_{h}(t_{j})-\textbf{H}_{h}^{j}\|_{r}^{2}\leq \kappa_{6}(\Delta t)^{2-r},\ \ r=-2, -1, 0, 1, 2,\ \ j=1,2.
 \end{align*}
\end{lemma}

For simplicity, define $e^{n}=\textbf{u}_{h}(t_{n})-\textbf{u}_{h}^{n}$, $\eta^{n}=\bar{p}_{h}^{n}-p_{h}^{n}$, $\xi^{n}=\textbf{H}_{h}(t_{n})-\textbf{H}_{h}^{n},$
\begin{align*}
&d_t\textbf{u}_{h}(t_{n})=\frac{\textbf{u}_{h}(t_{n+1})-\textbf{u}_{h}(t_{n-1})}{2\Delta t}, \ \ \bar{b}^{n}(\textbf{u}_{h},\textbf{u}_{h},\textbf{v}_{h})=\frac{1}{2\Delta t}\int_{t_{n-1}}^{t_{n+1}}{b(\textbf{u}_{h}(t),\textbf{u}_{h}(t),\textbf{v}_{h})}dt,\\
&\tilde{b}^{n}(\textbf{u}_{h},\textbf{u}_{h},\textbf{v}_{h})=\frac{1}{2}(b(\textbf{u}_{h}(t_{n+1}),\textbf{u}_{h}(t_{n+1}),\textbf{v}_{h})+b(\textbf{u}_{h}(t_{n-1}),\textbf{u}_{h}(t_{n-1}),\textbf{v}_{h})).
\end{align*}
Integrating \refe{2.9} and \refe{2.10} from $t_{n-1}$ to $t_{n+1}$ yields
\begin{align}\label{3.7}
&(d_{t}\textbf{u}_{h}(t_{n}),\textbf{v}_{h})+a(\bar{\textbf{u}}_{h}^{n},\textbf{v}_{h})-d(\textbf{v}_{h},\bar{p}_{h}^{n})+d(\bar{\textbf{u}}_{h}^{n},q_{h})+\bar{b}^{n}(\textbf{u}_{h},\textbf{u}_{h},\textbf{v}_{h})\nn\\
&\ \ +\frac{\mu}{2\Delta t}\int_{t_{n-1}}^{t_{n+1}}(\textbf{H}_{h}(t)\times\operatorname{curl}\textbf{H}_{h}(t),\textbf{v}_{h})dt=(\overline{\textbf{f}}^{n},\textbf{v}_{h}),\\
\label{3.8}
&\mu(d_{t}\textbf{H}_{h}(t_{n}),\textbf{B}_{h})+\sigma^{-1}(\operatorname{curl}\bar{\textbf{H}}_{h}^{n},\operatorname{curl}\textbf{B}_{h})
+\sigma^{-1}(\operatorname{div}\bar{\textbf{H}}_{h}^{n},\operatorname{div}\textbf{B}_{h})\nn\\
&\ \ -\frac{\mu}{2\Delta t}\int_{t_{n-1}}^{t_{n+1}}(\textbf{u}_{h}(t)\times \textbf{H}_{h}(t),\operatorname{curl}\textbf{B}_{h})dt=(\bar{\textbf{g}}^{n}, \textbf{B}_{h}).
\end{align}

	Subtracting \refe{3.1} and \refe{3.2} from \refe{3.7} and \refe{3.8}, respectively, we have
\begin{align}\label{3.9}
&(d_{t}e^{n},\textbf{v}_{h})+a(\tilde{e}^{n},\textbf{v}_{h})-d(\textbf{v}_{h},\eta^{n})+d(\tilde{e}^{n},q_{h})=-d(\bar{\textbf{u}}_{h}^{n}-\tilde{\textbf{u}}_{h}^{n},q_{h})-b(\textbf{u}_{h}^{n},\tilde{e}^{n},\textbf{v}_{h})\nn\\
&\ \ -b(e^{n},\tilde{\textbf{u}}_{h}(t_{n}),\textbf{v}_{h})-\mu(\textbf{H}_{h}^{n}\times\operatorname{curl}\tilde{\xi}^{n},\textbf{v}_{h})-\mu(\xi^{n}\times\operatorname{curl}\tilde{\textbf{H}}_{h}(t_{n}),\textbf{v}_{h})+(\textbf{E}^{n},\textbf{v}_{h}),\\
\label{3.10}
&\mu(d_{t}\xi^{n},\textbf{B}_{h})+\sigma^{-1}(\operatorname{curl}\tilde{\xi}^{n},\operatorname{curl}\textbf{B}_{h})+\sigma^{-1}(\operatorname{div}\tilde{\xi}^{n},\operatorname{div}\textbf{B}_{h})=\mu(\tilde{\textbf{u}}_{h}(t_{n})\times\xi^{n},\operatorname{curl}\textbf{B}_{h})\nn\\
&\ \ +\mu(\tilde{e}^{n}\times\textbf{H}_{h}^{n},\operatorname{curl}\textbf{B}_{h})+(\textbf{F}^{n},\textbf{B}_{h}),
\end{align}
where
\begin{align}\label{3.11}
&(\textbf{E}^{n},\textbf{v}_{h})=-a(\bar{\textbf{u}}_{h}^{n}-\tilde{\textbf{u}}_{h}(t_{n}),\textbf{v}_{h})+(\overline{\textbf{f}}^{n}-\textbf{f}(t_{n}),\textbf{v}_{h})-(\bar{b}^{n}(\textbf{u}_{h},\textbf{u}_{h},\textbf{v}_{h})-b(\textbf{u}_{h}(t_{n}),
\tilde{\textbf{u}}_{h}(t_{n}),\textbf{v}_{h}))\nn\\
&\ \ \ \ \ \ \ \ \ \ \ \ \ \ \ -\mu(\frac{1}{2\Delta t}\int_{t_{n-1}}^{t_{n+1}}(\textbf{H}_{h}(t)\times\operatorname{curl}{\textbf{H}_{h}(t)},\textbf{v}_{h})dt-(\textbf{H}_{h}(t_{n})\times\operatorname{curl}\tilde{\textbf{H}}_{h}(t_{n}),\textbf{v}_{h})),
\\ \label{3.12}
&(\textbf{F}^{n},\textbf{B}_{h})=-\sigma^{-1}((\operatorname{curl}(\bar{\textbf{H}}_{h}^{n}-\tilde{\textbf{H}}_{h}(t_{n})),\operatorname{curl}\textbf{B}_{h})+(\operatorname{div}(\bar{\textbf{H}}_{h}^{n}-\tilde{\textbf{H}}_{h}(t_{n})),\operatorname{div}\textbf{B}_{h}))+(\bar{\textbf{g}}^{n}
-\textbf{g}(t_{n}),\textbf{B}_{h})\nn\\
&\ \ \ \ \ \ \ \ \ \ \ \ \ \ \  +\mu(\frac{1}{2\Delta t}\int_{t_{n-1}}^{t_{n+1}}(\textbf{u}_{h}(t)\times\textbf{H}_{h}(t),\operatorname{curl}\textbf{B}_{h})dt-(\tilde{\textbf{u}}_{h}(t_{n})\times\textbf{H}_{h}(t_{n}),\operatorname{curl}\textbf{B}_{h})).
\end{align}
Using Taylor's formula and combining \refe{3.11} with \refe{3.12}, we deduce that
\begin{align}\label{3.13}
(\textbf{E}^{n},\textbf{v}_{h})&=a(\frac{1}{4\Delta t}\int_{t_{n-1}}^{t_{n+1}}(t_{n+1}-t)(t-t_{n-1})\textbf{u}_{htt}(t)dt,\textbf{v}_{h})\nn\\
&\ \ \ +\frac{1}{4\Delta t}(\int_{t_{n-1}}^{t_{n+1}}(t-t_{n+1})(t-t_{n-1})\textbf{f}_{tt}(t)dt
,\textbf{v}_{h})\nn\\
&\ \ \ -\frac{1}{2}(\int_{t_{n}}^{t_{n+1}}(t-t_{n+1})\textbf{f}_{tt}(t)dt-\int_{t_{n-1}}^{t_{n}}(t-t_{n-1})\textbf{f}_{tt}(t)dt,\textbf{v}_{h})\nn\\
&\ \ \ +\frac{1}{4\Delta t}
\int_{t_{n-1}}^{t_{n+1}}(t_{n+1}-t)(t-t_{n-1})b_{tt}(\textbf{u}_{h}(t),\textbf{u}_{h}(t),\textbf{v}_{h})dt\nn\\
&\ \ \ -b(\textbf{u}_{h}(t_{n}),\frac{1}{2}\int_{t_{n}}^{t_{n+1}}(t-t_{n+1})\textbf{u}_{htt}(t)dt-\frac{1}{2}\int_{t_{n-1}}^{t_{n}}(t-t_{n-1})\textbf{u}_{htt}(t)dt,\textbf{v}_{h})\nn\\
&\ \ \ +\frac{1}{2}\int_{t_{n}}^{t_{n+1}}(t-t_{n+1})b_{tt}(\textbf{u}_{h}(t),\textbf{u}_{h}(t),\textbf{v}_{h})dt-\frac{1}{2}\int_{t_{n-1}}^{t_{n}}(t-t_{n-1})b_{tt}(\textbf{u}_{h}(t),\textbf{u}_{h}(t),\textbf{v}_{h})dt\nn\\
&\ \ \ +\frac{\mu}{4\Delta t}\int_{t_{n-1}}^{t_{n+1}}(t_{n+1}-t)(t-t_{n-1})(\textbf{H}_{h}(t)\times \operatorname{curl}\textbf{H}_{h}(t),\textbf{v}_{h})_{tt}dt\nn\\
&\ \ \ -\mu(\textbf{H}_{h}(t_{n})\times\operatorname{curl}(\frac{1}{2}\int_{t_{n}}^{t_{n+1}}(t-t_{n+1})\textbf{H}_{htt}(t)dt-\frac{1}{2}\int_{t_{n-1}}^{t_{n}}(t-t_{n-1})\textbf{H}_{htt}(t)dt),\textbf{v}_{h})\nn\\
&\ \ \ + \frac{\mu}{2}\int_{t_{n}}^{t_{n+1}}(t-t_{n+1})(\textbf{H}_{h}(t)\times \operatorname{curl}\textbf{H}_{h}(t),\textbf{v}_{h})_{tt}dt\nn\\
&\ \ \
-\frac{\mu}{2}\int_{t_{n-1}}^{t_{n}}(t-t_{n-1})(\textbf{H}_{h}(t)\times \operatorname{curl}\textbf{H}_{h}(t),\textbf{v}_{h})_{tt}dt,
\\
\label{3.14}
(\textbf{F}^{n},\textbf{B}_{h})&=\frac{\sigma^{-1}}{4\Delta t}(\int_{t_{n-1}}^{t_{n+1}}(t_{n+1}-t)(t-t_{n-1})\operatorname{curl}\textbf{H}_{htt}(t)dt,\operatorname{curl}\textbf{B}_{h})\nn\\
&\ \ \ +\frac{\sigma^{-1}}{4\Delta t}(\int_{t_{n-1}}^{t_{n+1}}(t_{n+1}-t)(t-t_{n-1})\operatorname{div}\textbf{H}_{htt}(t)dt,\operatorname{div}\textbf{B}_{h})\nn\\
&\ \ \  +\frac{1}{4\Delta t}(\int_{t_{n-1}}^{t_{n+1}}(t-t_{n+1})(t-t_{n-1})\textbf{g}_{tt}(t)dt,\textbf{B}_{h})\nn\\
&\ \ \  -\frac{1}{2}(\int_{t_{n}}^{t_{n+1}}(t-t_{n+1})\textbf{g}_{tt}(t)dt-\int_{t_{n-1}}^{t_{n}}(t-t_{n-1})\textbf{g}_{tt}(t)dt,\textbf{B}_{h})\nn\\
&\ \ \  -\frac{\mu}{4\Delta t}\int_{t_{n-1}}^{t_{n+1}}(t_{n+1}-t)(t-t_{n-1})(\textbf{u}_{h}(t)\times\textbf{H}_{h}(t),\operatorname{curl}\textbf{B}_{h})_{tt}dt\nn\\
&\ \ \  +\mu((\frac{1}{2}\int_{t_{n}}^{t_{n+1}}(t-t_{n+1})\textbf{u}_{htt}(t)dt-\frac{1}{2}\int_{t_{n-1}}^{t_{n}}(t-t_{n-1})\textbf{u}_{htt}(t)dt)\times\textbf{H}_{h}(t_{n}),\operatorname{curl}\textbf{B}_{h})\nn\\
&\ \ \ -\frac{\mu}{2}\int_{t_{n}}^{t_{n+1}}(t-t_{n+1})(\textbf{u}_{h}(t)\times\textbf{H}_{h}(t),\operatorname{curl}\textbf{B}_{h})_{tt}dt\nn\\
&\ \ \ +\frac{\mu}{2}\int_{t_{n-1}}^{t_{n}}(t-t_{n-1})(\textbf{u}_{h}(t)\times\textbf{H}_{h}(t),\operatorname{curl}\textbf{B}_{h})_{tt}dt,
\end{align}
where we have used the following integral identities
\begin{align}\label{3.15}
&u(s)=u(t_{k})+(s-t_{k})u_{t}(t_{k})+\int_{s}^{t_{k}}(t-s)u_{htt}(t)dt,\\
\label{3.16}
&\frac{1}{2\Delta t}\int_{t_{k-1}}^{t_{k+1}}u(t)dt-\frac{u(t_{k+1})+u(t_{k-1})}{2}=\frac{1}{4\Delta t}\int_{t_{k-1}}^{t_{k+1}}(t-t_{k+1})(t-t_{k-1})u_{tt}(t)dt.
\end{align}
\begin{lemma}\label{T3.1}
Suppose that Assumptions (A0-A2) and $\Delta t\leq \frac{1}{2}$ are satisfied, then it is established that
\begin{align}\label{3.17}
&\Delta t\sum_{n=2}^{m}\|A_{1h}^{-\frac{3}{2}}P_{h}\textbf{E}^{n}\|_{0}^{2}\leq \lambda_{0}(\Delta t)^{4},\  \ \ \ \ \ \ \ \ \Delta t\sum_{n=2}^{m}\|A_{2h}^{-\frac{3}{2}}R_{0h}\textbf{F}^{n}\|_{0}^{2}\leq \lambda_{0}(\Delta t)^{4},\\
\label{3.18}
&\Delta t\sum_{n=2}^{m}\tau^{i}(t_{n+1})\|A_{1h}^{-1}P_{h}\textbf{E}^{n}\|_{0}^{2}\leq\lambda_{0}(\Delta t)^{3+i}, \nn\\
 &\Delta t\sum_{n=2}^{m}\tau^{i}(t_{n+1})\|A_{2h}^{-1}R_{0h}\textbf{F}^{n}\|_{0}^{2}\leq\lambda_{0}(\Delta t)^{3+i},  i=0,1,\\
\label{3.19}
&\Delta t\sum_{n=2}^{m}\tau^{i}(t_{n+1})\|A_{1h}^{-\frac{1}{2}}P_{h}\textbf{E}^{n}\|_{0}^{2}\leq\lambda_{0}(\Delta t)^{2+i},  \nn\\
 &\Delta t\sum_{n=2}^{m}\tau^{i}(t_{n+1})\|A_{2h}^{-\frac{1}{2}}R_{0h}\textbf{F}^{n}\|_{0}^{2}\leq\lambda_{0}(\Delta t)^{2+i},   i=0,1, 2,\\
\label{3.20}
&\Delta t\sum_{n=2}^{m}\tau^{i}(t_{n+1})\|P_{h}\textbf{E}^{n}\|_{0}^{2}\leq\lambda_{0}(\Delta t)^{1+i},    \  \Delta t\sum_{n=2}^{m}\tau^{i}(t_{n+1})\|R_{0h}\textbf{F}^{n}\|_{0}^{2}\leq\lambda_{0}(\Delta t)^{1+i},  i=0,1, 2,3,\\
\label{3.21}
&\Delta t\sum_{n=2}^{m}\tau^{i}(t_{n+1})\|P_{h}\textbf{E}^{n}\|_{1}^{2}\leq\lambda_{0}(\Delta t)^{i},  \ \ \ \ \Delta t\sum_{n=2}^{m}\tau^{i}(t_{n+1})\|R_{0h}\textbf{F}^{n}\|_{1}^{2}\leq\lambda_{0}(\Delta t)^{i},  i=0,1, 2,3,4,
\end{align}
for all $2\leq m \leq N-1$. Here and after, $\lambda_{i}$ ($i=0$, $1$, $2$, $\ldots$) is a positive constant depending on $\Omega$, $\nu$, $\sigma$, $\mu$, $T$, $\textbf{u}_{0}$, $\textbf{H}_{0}$, $\textbf{f}$ and $\textbf{g}$.
\end{lemma}
\begin{proof}
Applying \refe{2.14}\--\refe{2.19} and Cauchy-Schwarz inequality to \refe{3.13} and \refe{3.14}, we can get
\begin{align}\label{3.22}
\|A_{1h}^{-\frac{3}{2}}P_{h}\textbf{E}^{n}\|_{0}&=\sup_{\textbf{v}_{h}\in X_{0h}}\frac{(\textbf{E}^{n},\textbf{v}_{h})}{\|A_{1h}^{\frac{3}{2}}\textbf{v}_{h}\|_{0}}\leq c(\Delta t)^{-\frac{1}{2}}(\int_{t_{n-1}}^{t_{n+1}}(t_{n+1}-t)^{2}(t-t_{n-1})^{2}\|\textbf{u}_{htt}(t)\|_{-1}^{2}dt)^{\frac{1}{2}}\nn\\
&\ \ \ \  +c(\Delta t)^{-\frac{1}{2}}(\int_{t_{n-1}}^{t_{n+1}}(t_{n+1}-t)^{2}(t-t_{n-1})^{2}\|\textbf{f}_{tt}(t)\|_{0}^{2}dt)^{\frac{1}{2}}\nn\\
&\ \ \ \  +c(\Delta t)^{\frac{1}{2}}[(\int_{t_{n}}^{t_{n+1}}(t_{n+1}-t)^{2}\|\textbf{f}_{tt}(t)\|_{0}^{2}dt)^{\frac{1}{2}}+(\int_{t_{n-1}}^{t_{n}}(t-t_{n-1})^{2}\|\textbf{f}_{tt}(t)\|_{0}^{2}dt)^{\frac{1}{2}}]\nn\\
&\ \ \ \  +c(\Delta t)^{-\frac{1}{2}}(\int_{t_{n-1}}^{t_{n+1}}(t_{n+1}-t)^{2}(t-t_{n-1})^{2}\|\textbf{u}_{htt}(t)\|_{-1}^{2}\|\textbf{u}_{h}(t)\|_{2}^{2}dt)^{\frac{1}{2}}\nn\\
&\ \ \ \ +c(\Delta t)^{-\frac{1}{2}}(\int_{t_{n-1}}^{t_{n+1}}(t_{n+1}-t)^{2}(t-t_{n-1})^{2}\|\textbf{u}_{ht}(t)\|_{0}^{2}\|\textbf{u}_{ht}(t)\|_{1}^{2}dt)^{\frac{1}{2}}\nn\\
&\ \ \ \ +c(\Delta t)^{\frac{1}{2}}\|\textbf{u}_{h}(t_{n})\|_{2}[(\int_{t_{n}}^{t_{n+1}}(t_{n+1}-t)^{2}\|\textbf{u}_{htt}(t)\|_{-1}^{2}dt)^{\frac{1}{2}}\nn\\&\ \ \ \ +(\int_{t_{n-1}}^{t_{n}}(t-t_{n-1})^{2}\|\textbf{u}_{htt}(t)\|_{-1}^{2}dt)^{\frac{1}{2}}]\nn\\
&\ \ \ \  +c(\Delta t)^{\frac{1}{2}}[(\int_{t_{n}}^{t_{n+1}}(t_{n+1}-t)^{2}\|\textbf{u}_{htt}(t)\|_{-1}^{2}\|\textbf{u}_{h}(t)\|_{2}^{2}dt)^{\frac{1}{2}}\nn\\&\ \ \ \ +(\int_{t_{n-1}}^{t_{n}}(t-t_{n-1})^{2}\|\textbf{u}_{htt}(t)\|_{-1}^{2}\|\textbf{u}_{h}(t)\|_{2}^{2}dt)^{\frac{1}{2}}]\nn\\
&\ \ \ \  +c(\Delta t)^{\frac{1}{2}}[(\int_{t_{n}}^{t_{n+1}}(t_{n+1}-t)^{2}\|\textbf{u}_{ht}(t)\|_{0}^{2}\|\textbf{u}_{ht}(t)\|_{1}^{2}dt)^{\frac{1}{2}}\nn\\ & \ \ \ \ +(\int_{t_{n-1}}^{t_{n}}(t-t_{n-1})^{2}\|\textbf{u}_{ht}(t)\|_{0}^{2}\|\textbf{u}_{ht}(t)\|_{1}^{2}dt)^{\frac{1}{2}}]\nn\\
&\ \ \ \  +c(\Delta t)^{-\frac{1}{2}}(\int_{t_{n-1}}^{t_{n+1}}(t_{n+1}-t)^{2}(t-t_{n-1})^{2}\|\textbf{H}_{htt}(t)\|_{-1}^{2}\|A_{2h}\textbf{H}_{h}(t)\|_{0}^{2}dt)^{\frac{1}{2}}\nn\\
&\ \ \ \  +c(\Delta t)^{-\frac{1}{2}}(\int_{t_{n-1}}^{t_{n+1}}(t_{n+1}-t)^{2}(t-t_{n-1})^{2}\|\textbf{H}_{ht}(t)\|_{0}^{2}\|\textbf{H}_{ht}(t)\|_{1}^{2}dt)^{\frac{1}{2}}\nn\\
&\ \ \ \  +c(\Delta t)^{\frac{1}{2}}\|A_{2h}\textbf{H}_{h}(t_{n})\|_{0}[(\int_{t_{n}}^{t_{n+1}}(t_{n+1}-t)^{2}\|\textbf{H}_{htt}(t)\|_{-1}^{2}dt)^{\frac{1}{2}}\nn\\
&\ \ \ \ +(\int_{t_{n-1}}^{t_{n}}(t-t_{n-1})^{2}\|\textbf{H}_{htt}(t)\|_{-1}^{2}dt)^{\frac{1}{2}}]\nn\\
&\ \ \ \  +c(\Delta t)^{\frac{1}{2}}[(\int_{t_{n}}^{t_{n+1}}(t_{n+1}-t)^{2}\|\textbf{H}_{htt}(t)\|_{-1}^{2}\|A_{2h}\textbf{H}_{h}(t)\|_{0}^{2}dt)^{\frac{1}{2}}\nn\\&\ \ \ \ +(\int_{t_{n-1}}^{t_{n}}(t-t_{n-1})^{2}\|\textbf{H}_{htt}(t)\|_{-1}^{2}\|A_{2h}\textbf{H}_{h}(t)\|_{0}^{2}dt)^{\frac{1}{2}}]\nn\\
&\ \ \ \  +c(\Delta t)^{\frac{1}{2}}[(\int_{t_{n}}^{t_{n+1}}(t_{n+1}-t)^{2}\|\textbf{H}_{ht}(t)\|_{0}^{2}\|\textbf{H}_{ht}(t)\|_{1}^{2}dt)^{\frac{1}{2}}\nn\\ & \ \ \ \ +(\int_{t_{n-1}}^{t_{n}}(t-t_{n-1})^{2}\|\textbf{H}_{ht}(t)\|_{0}^{2}\|\textbf{H}_{ht}(t)\|_{1}^{2}dt)^{\frac{1}{2}}],\\
\label{3.23}
\|A_{2h}^{-\frac{3}{2}}R_{0h}\textbf{F}^{n}\|_{0}&=\sup_{\textbf{B}_{h}\in W_{h}}\frac{(\textbf{F}^{n},\textbf{B}_{h})}{\|A_{2h}^{\frac{3}{2}}\textbf{B}_{h}\|_{0}}\leq
c(\Delta t)^{-\frac{1}{2}}(\int_{t_{n-1}}^{t_{n+1}}(t_{n+1}-t)^{2}(t-t_{n-1})^{2}\|\textbf{H}_{htt}(t)\|_{-1}^{2}dt)^{\frac{1}{2}}\nn\\
&\ \ \ \  +c(\Delta t)^{-\frac{1}{2}}(\int_{t_{n-1}}^{t_{n+1}}(t_{n+1}-t)^{2}(t-t_{n-1})^{2}\|\textbf{g}_{tt}(t)\|_{0}^{2}dt)^{\frac{1}{2}}\nn\\
&\ \ \ \  +c(\Delta t)^{\frac{1}{2}}[(\int_{t_{n}}^{t_{n+1}}(t_{n+1}-t)^{2}\|\textbf{g}_{tt}(t)\|_{0}^{2}dt)^{\frac{1}{2}}+(\int_{t_{n-1}}^{t_{n}}(t-t_{n-1})^{2}\|\textbf{g}_{tt}(t)\|_{0}^{2}dt)^{\frac{1}{2}}]\nn\\
&\ \ \ \  +c(\Delta t)^{-\frac{1}{2}}[(\int_{t_{n-1}}^{t_{n+1}}(t_{n+1}-t)^{2}(t-t_{n-1})^{2}\|\textbf{u}_{htt}(t)\|_{-1}^{2}\|A_{2h}\textbf{H}_{h}(t)\|_{0}^{2}dt)^{\frac{1}{2}}\nn\\
&\ \ \ \ +(\int_{t_{n-1}}^{t_{n+1}}(t_{n+1}-t)^{2}(t-t_{n-1})^{2}\|A_{1h}\textbf{u}_{h}(t)\|_{0}^{2}\|\textbf{H}_{htt}(t)\|_{-1}^{2}dt)^{\frac{1}{2}}\nn\\
&\ \ \ \  +(\int_{t_{n-1}}^{t_{n+1}}(t_{n+1}-t)^{2}(t-t_{n-1})^{2}\|\textbf{u}_{ht}(t)\|_{0}^{2}\|\textbf{H}_{ht}(t)\|_{1}^{2}dt)^{\frac{1}{2}}]\nn\\
&\ \ \ \  +c(\Delta t)^{\frac{1}{2}}\|A_{2h}\textbf{H}_{h}(t_{n})\|_{0}[(\int_{t_{n}}^{t_{n+1}}(t_{n+1}-t)^{2}\|\textbf{u}_{htt}(t)\|_{-1}^{2}dt)^{\frac{1}{2}}\nn\\
&\ \ \ \  +(\int_{t_{n-1}}^{t_{n}}(t-t_{n-1})^{2}\|\textbf{u}_{htt}(t)\|_{-1}^{2}dt)^{\frac{1}{2}}]
\nn\\
&\ \ \ \  +c(\Delta t)^{\frac{1}{2}}[(\int_{t_{n}}^{t_{n+1}}(t_{n+1}-t)^{2}\|\textbf{u}_{htt}(t)\|_{-1}^{2}\|A_{2h}\textbf{H}_{h}(t)\|_{0}^{2}dt)^{\frac{1}{2}}\nn\\
&\ \ \ \  +(\int_{t_{n}}^{t_{n+1}}(t_{n+1}-t)^{2}\|A_{1h}\textbf{u}_{h}(t)\|_{0}^{2}\|\textbf{H}_{htt}(t)\|_{-1}^{2}dt)^{\frac{1}{2}}]
\nn\\
&\ \ \ \  +c(\Delta t)^{\frac{1}{2}}[(\int_{t_{n-1}}^{t_{n}}(t-t_{n-1})^{2}\|\textbf{u}_{htt}(t)\|_{-1}^{2}\|A_{2h}\textbf{H}_{h}(t)\|_{0}^{2}dt)^{\frac{1}{2}}\nn\\
&\ \ \ \  +(\int_{t_{n-1}}^{t_{n}}(t-t_{n-1})^{2}\|A_{1h}\textbf{u}_{h}(t)\|_{0}^{2}\|\textbf{H}_{htt}(t)\|_{-1}^{2}dt)^{\frac{1}{2}}]
\nn\\
&\ \ \ \  +c(\Delta t)^{\frac{1}{2}}[(\int_{t_{n}}^{t_{n+1}}(t_{n+1}-t)^{2}\|\textbf{u}_{ht}(t)\|_{0}^{2}\|\textbf{H}_{ht}(t)\|_{1}^{2}dt)^{\frac{1}{2}}\nn\\
&\ \ \ \  +(\int_{t_{n-1}}^{t_{n}}(t-t_{n-1})^{2}\|\textbf{u}_{ht}(t)\|_{0}^{2}\|\textbf{H}_{ht}(t)\|_{1}^{2}dt)^{\frac{1}{2}}].
\end{align}
In the same way, there holds that\begin{align}
\label{3.24}
\|A_{1h}^{\frac{\alpha-1}{2}}P_{h}\textbf{E}^{n}\|_{0}&=\sup_{\textbf{v}_{h}\in X_{0h}}\frac{(\textbf{E}^{n},\textbf{v}_{h})}{\|A_{1h}^{\frac{1-\alpha}{2}}\textbf{v}_{h}\|_{0}}\leq c(\Delta t)^{-\frac{1}{2}}(\int_{t_{n-1}}^{t_{n+1}}(t_{n+1}-t)^{2}(t-t_{n-1})^{2}\|A_{1h}^{\frac{\alpha+1}{2}}\textbf{u}_{htt}(t)\|_{0}^{2}dt)^{\frac{1}{2}}\nn\\
&\ \ \ \ +c(\Delta t)^{-\frac{1}{2}}(\int_{t_{n-1}}^{t_{n+1}}(t_{n+1}-t)^{2}(t-t_{n-1})^{2}\|\textbf{f}_{tt}(t)\|_{0}^{2}dt)^{\frac{1}{2}}\nn\\
&\ \ \ \  +c(\Delta t)^{\frac{1}{2}}[(\int_{t_{n}}^{t_{n+1}}(t_{n+1}-t)^{2}\|\textbf{f}_{tt}(t)\|_{0}^{2}dt)^{\frac{1}{2}}+(\int_{t_{n-1}}^{t_{n}}(t-t_{n-1})^{2}\|\textbf{f}_{tt}(t)\|_{0}^{2}dt)^{\frac{1}{2}}]\nn\\
&\ \ \ \  +c(\Delta t)^{-\frac{1}{2}}(\int_{t_{n-1}}^{t_{n+1}}(t_{n+1}-t)^{2}(t-t_{n-1})^{2}\|\textbf{u}_{htt}(t)\|_{\alpha}^{2}\|\textbf{u}_{h}(t)\|_{2}^{2}dt)^{\frac{1}{2}}\nn\\
&\ \ \ \ +c(\Delta t)^{-\frac{1}{2}}(\int_{t_{n-1}}^{t_{n+1}}(t_{n+1}-t)^{2}(t-t_{n-1})^{2}\|\textbf{u}_{ht}(t)\|_{\alpha+1}^{2}\|\textbf{u}_{ht}(t)\|_{1}^{2}dt)^{\frac{1}{2}}\nn\\
&\ \ \ \ +c(\Delta t)^{\frac{1}{2}}\|\textbf{u}_{h}(t_{n})\|_{2}[(\int_{t_{n}}^{t_{n+1}}(t_{n+1}-t)^{2}\|\textbf{u}_{htt}(t)\|_{\alpha}^{2}dt)^{\frac{1}{2}}\nn\\&\ \ \ \ +(\int_{t_{n-1}}^{t_{n}}(t-t_{n-1})^{2}\|\textbf{u}_{htt}(t)\|_{\alpha}^{2}dt)^{\frac{1}{2}}]\nn\\
&\ \ \ \  +c(\Delta t)^{\frac{1}{2}}[(\int_{t_{n}}^{t_{n+1}}(t_{n+1}-t)^{2}\|\textbf{u}_{htt}(t)\|_{\alpha}^{2}\|\textbf{u}_{h}(t)\|_{2}^{2}dt)^{\frac{1}{2}}\nn\\&\ \ \ \ +(\int_{t_{n-1}}^{t_{n}}(t-t_{n-1})^{2}\|\textbf{u}_{htt}(t)\|_{\alpha}^{2}\|\textbf{u}_{h}(t)\|_{2}^{2}dt)^{\frac{1}{2}}]\nn\\
&\ \ \ \  +c(\Delta t)^{\frac{1}{2}}[(\int_{t_{n}}^{t_{n+1}}(t_{n+1}-t)^{2}\|\textbf{u}_{ht}(t)\|_{\alpha+1}^{2}\|\textbf{u}_{ht}(t)\|_{1}^{2}dt)^{\frac{1}{2}}\nn\\ & \ \ \ \ +(\int_{t_{n-1}}^{t_{n}}(t-t_{n-1})^{2}\|\textbf{u}_{ht}(t)\|_{\alpha+1}^{2}\|\textbf{u}_{ht}(t)\|_{1}^{2}dt)^{\frac{1}{2}}]\nn\\
&\ \ \ \  +c(\Delta t)^{-\frac{1}{2}}(\int_{t_{n-1}}^{t_{n+1}}(t_{n+1}-t)^{2}(t-t_{n-1})^{2}\|\textbf{H}_{htt}(t)\|_{\alpha}^{2}\|A_{2h}\textbf{H}_{h}(t)\|_{0}^{2}dt)^{\frac{1}{2}}\nn\\
&\ \ \ \  +c(\Delta t)^{-\frac{1}{2}}(\int_{t_{n-1}}^{t_{n+1}}(t_{n+1}-t)^{2}(t-t_{n-1})^{2}\|\textbf{H}_{ht}(t)\|_{\alpha+1}^{2}\|\textbf{H}_{ht}(t)\|_{1}^{2}dt)^{\frac{1}{2}}\nn\\
&\ \ \ \  +c(\Delta t)^{\frac{1}{2}}\|A_{2h}\textbf{H}_{h}(t_{n})\|_{0}[(\int_{t_{n}}^{t_{n+1}}(t_{n+1}-t)^{2}\|\textbf{H}_{htt}(t)\|_{\alpha}^{2}dt)^{\frac{1}{2}}\nn\\
&\ \ \ \ +(\int_{t_{n-1}}^{t_{n}}(t-t_{n-1})^{2}\|\textbf{H}_{htt}(t)\|_{\alpha}^{2}dt)^{\frac{1}{2}}]\nn\\
&\ \ \ \  +c(\Delta t)^{\frac{1}{2}}[(\int_{t_{n}}^{t_{n+1}}(t_{n+1}-t)^{2}\|\textbf{H}_{htt}(t)\|_{\alpha}^{2}\|A_{2h}\textbf{H}_{h}(t)\|_{0}^{2}dt)^{\frac{1}{2}}\nn\\&\ \ \ \ +(\int_{t_{n-1}}^{t_{n}}(t-t_{n-1})^{2}\|\textbf{H}_{htt}(t)\|_{\alpha}^{2}\|A_{2h}\textbf{H}_{h}(t)\|_{0}^{2}dt)^{\frac{1}{2}}]\nn\\
&\ \ \ \  +c(\Delta t)^{\frac{1}{2}}[(\int_{t_{n}}^{t_{n+1}}(t_{n+1}-t)^{2}\|\textbf{H}_{ht}(t)\|_{\alpha+1}^{2}\|\textbf{H}_{ht}(t)\|_{1}^{2}dt)^{\frac{1}{2}}\nn\\ & \ \ \ \ +(\int_{t_{n-1}}^{t_{n}}(t-t_{n-1})^{2}\|\textbf{H}_{ht}(t)\|_{\alpha+1}^{2}\|\textbf{H}_{ht}(t)\|_{1}^{2}dt)^{\frac{1}{2}}],\\
\label{3.25}
\|A_{2h}^{\frac{\alpha-1}{2}}R_{0h}\textbf{F}^{n}\|_{0}&=\sup_{\textbf{B}_{h}\in W_{h}}\frac{(\textbf{F}^{n},\textbf{B}_{h})}{\|A_{2h}^{\frac{1-\alpha}{2}}\textbf{B}_{h}\|_{0}}\leq
c(\Delta t)^{-\frac{1}{2}}(\int_{t_{n-1}}^{t_{n+1}}(t_{n+1}-t)^{2}(t-t_{n-1})^{2}\|A_{2h}^{\frac{\alpha+1}{2}}\textbf{H}_{htt}(t)\|_{0}^{2}dt)^{\frac{1}{2}}\nn\\
&\ \ \ \  +c(\Delta t)^{-\frac{1}{2}}(\int_{t_{n-1}}^{t_{n+1}}(t_{n+1}-t)^{2}(t-t_{n-1})^{2}\|\textbf{g}_{tt}(t)\|_{0}^{2}dt)^{\frac{1}{2}}\nn\\
&\ \ \ \  +c(\Delta t)^{\frac{1}{2}}[(\int_{t_{n}}^{t_{n+1}}(t_{n+1}-t)^{2}\|\textbf{g}_{tt}(t)\|_{0}^{2}dt)^{\frac{1}{2}}+(\int_{t_{n-1}}^{t_{n}}(t-t_{n-1})^{2}\|\textbf{g}_{tt}(t)\|_{0}^{2}dt)^{\frac{1}{2}}]\nn\\
&\ \ \ \  +c(\Delta t)^{-\frac{1}{2}}[(\int_{t_{n-1}}^{t_{n+1}}(t_{n+1}-t)^{2}(t-t_{n-1})^{2}\|\textbf{u}_{htt}(t)\|_{\alpha}^{2}\|A_{2h}\textbf{H}_{h}(t)\|_{0}^{2}dt)^{\frac{1}{2}}\nn\\
&\ \ \ \ +(\int_{t_{n-1}}^{t_{n+1}}(t_{n+1}-t)^{2}(t-t_{n-1})^{2}\|A_{1h}\textbf{u}_{h}(t)\|_{0}^{2}\|\textbf{H}_{htt}(t)\|_{\alpha}^{2}dt)^{\frac{1}{2}}\nn\\
&\ \ \ \  +(\int_{t_{n-1}}^{t_{n+1}}(t_{n+1}-t)^{2}(t-t_{n-1})^{2}\|\textbf{u}_{ht}(t)\|_{\alpha+1}^{2}\|\textbf{H}_{ht}(t)\|_{1}^{2}dt)^{\frac{1}{2}}]\nn\\
&\ \ \ \  +c(\Delta t)^{\frac{1}{2}}\|A_{2h}\textbf{H}_{h}(t_{n})\|_{0}[(\int_{t_{n}}^{t_{n+1}}(t_{n+1}-t)^{2}\|\textbf{u}_{htt}(t)\|_{\alpha}^{2}dt)^{\frac{1}{2}}\nn\\
&\ \ \ \  +(\int_{t_{n-1}}^{t_{n}}(t-t_{n-1})^{2}\|\textbf{u}_{htt}(t)\|_{\alpha}^{2}dt)^{\frac{1}{2}}]
\nn\\
&\ \ \ \  +c(\Delta t)^{\frac{1}{2}}[(\int_{t_{n}}^{t_{n+1}}(t_{n+1}-t)^{2}\|\textbf{u}_{htt}(t)\|_{\alpha}^{2}\|A_{2h}\textbf{H}_{h}(t)\|_{0}^{2}dt)^{\frac{1}{2}}\nn\\
&\ \ \ \  +(\int_{t_{n}}^{t_{n+1}}(t_{n+1}-t)^{2}\|A_{1h}\textbf{u}_{h}(t)\|_{0}^{2}\|\textbf{H}_{htt}(t)\|_{\alpha}^{2}dt)^{\frac{1}{2}}]
\nn\\
&\ \ \ \  +c(\Delta t)^{\frac{1}{2}}[(\int_{t_{n-1}}^{t_{n}}(t-t_{n-1})^{2}\|\textbf{u}_{htt}(t)\|_{\alpha}^{2}\|A_{2h}\textbf{H}_{h}(t)\|_{0}^{2}dt)^{\frac{1}{2}}\nn\\
&\ \ \ \  +(\int_{t_{n-1}}^{t_{n}}(t-t_{n-1})^{2}\|A_{1h}\textbf{u}_{h}(t)\|_{0}^{2}\|\textbf{H}_{htt}(t)\|_{\alpha}^{2}dt)^{\frac{1}{2}}]
\nn\\
&\ \ \ \  +c(\Delta t)^{\frac{1}{2}}[(\int_{t_{n}}^{t_{n+1}}(t_{n+1}-t)^{2}\|\textbf{u}_{ht}(t)\|_{\alpha+1}^{2}\|\textbf{H}_{ht}(t)\|_{1}^{2}dt)^{\frac{1}{2}}\nn\\
&\ \ \ \  +(\int_{t_{n-1}}^{t_{n}}(t-t_{n-1})^{2}\|\textbf{u}_{ht}(t)\|_{\alpha+1}^{2}\|\textbf{H}_{ht}(t)\|_{1}^{2}dt)^{\frac{1}{2}}],
\end{align}
for $\alpha=-1,$ 0, 1. Similarly, we  derive that\begin{align}
\label{3.26}
\|P_{h}\textbf{E}^{n}\|_{1}&=\sup_{\textbf{v}_{h}\in X_{0h}}\frac{(\textbf{E}^{n},\textbf{v}_{h})}{\|A_{1h}^{-\frac{1}{2}}\textbf{v}_{h}\|_{0}}\leq c(\Delta t)^{-\frac{1}{2}}(\int_{t_{n-1}}^{t_{n+1}}(t_{n+1}-t)^{2}(t-t_{n-1})^{2}\|A_{1h}^{\frac{3}{2}}\textbf{u}_{htt}(t)\|_{0}^{2}dt)^{\frac{1}{2}}\nn\\
&\ \ \ \  +c(\Delta t)^{-\frac{1}{2}}(\int_{t_{n-1}}^{t_{n+1}}(t_{n+1}-t)^{2}(t-t_{n-1})^{2}\|\textbf{f}_{tt}(t)\|_{1}^{2}dt)^{\frac{1}{2}}\nn\\
&\ \ \ \  +c(\Delta t)^{\frac{1}{2}}[(\int_{t_{n}}^{t_{n+1}}(t_{n+1}-t)^{2}\|\textbf{f}_{tt}(t)\|_{1}^{2}dt)^{\frac{1}{2}}+(\int_{t_{n-1}}^{t_{n}}(t-t_{n-1})^{2}\|\textbf{f}_{tt}(t)\|_{1}^{2}dt)^{\frac{1}{2}}]\nn\\
&\ \ \ \  +c(\Delta t)^{-\frac{1}{2}}(\int_{t_{n-1}}^{t_{n+1}}(t_{n+1}-t)^{2}(t-t_{n-1})^{2}\|A_{1h}\textbf{u}_{htt}(t)\|_{0}^{2}\|A_{1h}\textbf{u}_{h}(t)\|_{0}^{2}dt)^{\frac{1}{2}}\nn\\
&\ \ \ \ +c(\Delta t)^{-\frac{1}{2}}(\int_{t_{n-1}}^{t_{n+1}}(t_{n+1}-t)^{2}(t-t_{n-1})^{2}\|\textbf{u}_{ht}(t)\|_{2}^{2}\|\textbf{u}_{ht}(t)\|_{2}^{2}dt)^{\frac{1}{2}}\nn\\
&\ \ \ \ +c(\Delta t)^{\frac{1}{2}}\|\textbf{u}_{h}(t_{n})\|_{2}[(\int_{t_{n}}^{t_{n+1}}(t_{n+1}-t)^{2}\|\textbf{u}_{htt}(t)\|_{2}^{2}dt)^{\frac{1}{2}}\nn\\&\ \ \ \ +(\int_{t_{n-1}}^{t_{n}}(t-t_{n-1})^{2}\|\textbf{u}_{htt}(t)\|_{2}^{2}dt)^{\frac{1}{2}}]\nn\\
&\ \ \ \  +c(\Delta t)^{\frac{1}{2}}[(\int_{t_{n}}^{t_{n+1}}(t_{n+1}-t)^{2}\|\textbf{u}_{htt}(t)\|_{2}^{2}\|\textbf{u}_{h}(t)\|_{2}^{2}dt)^{\frac{1}{2}}\nn\\&\ \ \ \ +(\int_{t_{n-1}}^{t_{n}}(t-t_{n-1})^{2}\|\textbf{u}_{htt}(t)\|_{2}^{2}\|\textbf{u}_{h}(t)\|_{2}^{2}dt)^{\frac{1}{2}}]\nn\\
&\ \ \ \  +c(\Delta t)^{\frac{1}{2}}[(\int_{t_{n}}^{t_{n+1}}(t_{n+1}-t)^{2}\|\textbf{u}_{ht}(t)\|_{2}^{2}\|\textbf{u}_{ht}(t)\|_{2}^{2}dt)^{\frac{1}{2}}\nn\\ & \ \ \ \ +(\int_{t_{n-1}}^{t_{n}}(t-t_{n-1})^{2}\|\textbf{u}_{ht}(t)\|_{2}^{2}\|\textbf{u}_{ht}(t)\|_{2}^{2}dt)^{\frac{1}{2}}]\nn\\
&\ \ \ \  +c(\Delta t)^{-\frac{1}{2}}(\int_{t_{n-1}}^{t_{n+1}}(t_{n+1}-t)^{2}(t-t_{n-1})^{2}\|A_{2h}\textbf{H}_{htt}(t)\|_{0}^{2}\|A_{2h}\textbf{H}_{h}(t)\|_{0}^{2}dt)^{\frac{1}{2}}\nn\\
&\ \ \ \  +c(\Delta t)^{-\frac{1}{2}}(\int_{t_{n-1}}^{t_{n+1}}(t_{n+1}-t)^{2}(t-t_{n-1})^{2}\|A_{2h}\textbf{H}_{ht}(t)\|_{0}^{2}\|A_{2h}\textbf{H}_{ht}(t)\|_{0}^{2}dt)^{\frac{1}{2}}\nn\\
&\ \ \ \  +c(\Delta t)^{\frac{1}{2}}\|A_{2h}\textbf{H}_{h}(t_{n})\|_{0}[(\int_{t_{n}}^{t_{n+1}}(t_{n+1}-t)^{2}\|A_{2h}\textbf{H}_{htt}(t)\|_{0}^{2}dt)^{\frac{1}{2}}\nn\\
&\ \ \ \ +(\int_{t_{n-1}}^{t_{n}}(t-t_{n-1})^{2}\|A_{2h}\textbf{H}_{htt}(t)\|_{0}^{2}dt)^{\frac{1}{2}}]\nn\\
&\ \ \ \  +c(\Delta t)^{\frac{1}{2}}[(\int_{t_{n}}^{t_{n+1}}(t_{n+1}-t)^{2}\|A_{2h}\textbf{H}_{htt}(t)\|_{0}^{2}\|A_{2h}\textbf{H}_{h}(t)\|_{0}^{2}dt)^{\frac{1}{2}}\nn\\&\ \ \ \ +(\int_{t_{n-1}}^{t_{n}}(t-t_{n-1})^{2}\|A_{2h}\textbf{H}_{htt}(t)\|_{0}^{2}\|A_{2h}\textbf{H}_{h}(t)\|_{0}^{2}dt)^{\frac{1}{2}}]\nn\\
&\ \ \ \  +c(\Delta t)^{\frac{1}{2}}[(\int_{t_{n}}^{t_{n+1}}(t_{n+1}-t)^{2}\|A_{2h}\textbf{H}_{ht}(t)\|_{0}^{2}\|A_{2h}\textbf{H}_{ht}(t)\|_{0}^{2}dt)^{\frac{1}{2}}\nn\\ & \ \ \ \ +(\int_{t_{n-1}}^{t_{n}}(t-t_{n-1})^{2}\|A_{2h}\textbf{H}_{ht}(t)\|_{0}^{2}\|A_{2h}\textbf{H}_{ht}(t)\|_{0}^{2}dt)^{\frac{1}{2}}],\\
\label{3.27}
\|R_{0h}\textbf{F}^{n}\|_{1}&=\sup_{\textbf{B}_{h}\in W_{h}}\frac{(\textbf{F}^{n},\textbf{B}_{h})}{\|A_{2h}^{-\frac{1}{2}}\textbf{B}_{h}\|_{0}}\leq
c(\Delta t)^{-\frac{1}{2}}(\int_{t_{n-1}}^{t_{n+1}}(t_{n+1}-t)^{2}(t-t_{n-1})^{2}\|A_{2h}^{\frac{3}{2}}\textbf{H}_{htt}(t)\|_{0}^{2}dt)^{\frac{1}{2}}\nn\\
&\ \ \ \  +c(\Delta t)^{-\frac{1}{2}}(\int_{t_{n-1}}^{t_{n+1}}(t_{n+1}-t)^{2}(t-t_{n-1})^{2}\|\textbf{g}_{tt}(t)\|_{1}^{2}dt)^{\frac{1}{2}}\nn\\
&\ \ \ \  +c(\Delta t)^{\frac{1}{2}}[(\int_{t_{n}}^{t_{n+1}}(t_{n+1}-t)^{2}\|\textbf{g}_{tt}(t)\|_{1}^{2}dt)^{\frac{1}{2}}+(\int_{t_{n-1}}^{t_{n}}(t-t_{n-1})^{2}\|\textbf{g}_{tt}(t)\|_{1}^{2}dt)^{\frac{1}{2}}]\nn\\
&\ \ \ \  +c(\Delta t)^{-\frac{1}{2}}[(\int_{t_{n-1}}^{t_{n+1}}(t_{n+1}-t)^{2}(t-t_{n-1})^{2}\|A_{1h}\textbf{u}_{htt}(t)\|_{0}^{2}\|A_{2h}\textbf{H}_{h}(t)\|_{0}^{2}dt)^{\frac{1}{2}}\nn\\
&\ \ \ \ +(\int_{t_{n-1}}^{t_{n+1}}(t_{n+1}-t)^{2}(t-t_{n-1})^{2}\|A_{1h}\textbf{u}_{h}(t)\|_{0}^{2}\|A_{2h}\textbf{H}_{htt}(t)\|_{0}^{2}dt)^{\frac{1}{2}}\nn\\
&\ \ \ \  +(\int_{t_{n-1}}^{t_{n+1}}(t_{n+1}-t)^{2}(t-t_{n-1})^{2}\|A_{1h}\textbf{u}_{ht}(t)\|_{0}^{2}\|A_{2h}\textbf{H}_{ht}(t)\|_{0}^{2}dt)^{\frac{1}{2}}]\nn\\
&\ \ \ \  +c(\Delta t)^{\frac{1}{2}}\|A_{2h}\textbf{H}_{h}(t_{n})\|_{0}[(\int_{t_{n}}^{t_{n+1}}(t_{n+1}-t)^{2}\|A_{1h}\textbf{u}_{htt}(t)\|_{0}^{2}dt)^{\frac{1}{2}}\nn\\
&\ \ \ \  +(\int_{t_{n-1}}^{t_{n}}(t-t_{n-1})^{2}\|A_{1h}\textbf{u}_{htt}(t)\|_{0}^{2}dt)^{\frac{1}{2}}]
\nn\\
&\ \ \ \  +c(\Delta t)^{\frac{1}{2}}[(\int_{t_{n}}^{t_{n+1}}(t_{n+1}-t)^{2}\|A_{1h}\textbf{u}_{htt}(t)\|_{0}^{2}\|A_{2h}\textbf{H}_{h}(t)\|_{0}^{2}dt)^{\frac{1}{2}}\nn\\
&\ \ \ \  +(\int_{t_{n}}^{t_{n+1}}(t_{n+1}-t)^{2}\|A_{1h}\textbf{u}_{h}(t)\|_{0}^{2}\|A_{2h}\textbf{H}_{htt}(t)\|_{0}^{2}dt)^{\frac{1}{2}}]
\nn\\
&\ \ \ \  +c(\Delta t)^{\frac{1}{2}}[(\int_{t_{n-1}}^{t_{n}}(t-t_{n-1})^{2}\|A_{1h}\textbf{u}_{htt}(t)\|_{0}^{2}\|A_{2h}\textbf{H}_{h}(t)\|_{0}^{2}dt)^{\frac{1}{2}}\nn\\
&\ \ \ \  +(\int_{t_{n-1}}^{t_{n}}(t-t_{n-1})^{2}\|A_{1h}\textbf{u}_{h}(t)\|_{0}^{2}\|A_{2h}\textbf{H}_{htt}(t)\|_{0}^{2}dt)^{\frac{1}{2}}]
\nn\\
&\ \ \ \  +c(\Delta t)^{\frac{1}{2}}[(\int_{t_{n}}^{t_{n+1}}(t_{n+1}-t)^{2}\|A_{1h}\textbf{u}_{ht}(t)\|_{0}^{2}\|A_{2h}\textbf{H}_{ht}(t)\|_{0}^{2}dt)^{\frac{1}{2}}\nn\\
&\ \ \ \  +(\int_{t_{n-1}}^{t_{n}}(t-t_{n-1})^{2}\|A_{1h}\textbf{u}_{ht}(t)\|_{0}^{2}\|A_{2h}\textbf{H}_{ht}(t)\|_{0}^{2}dt)^{\frac{1}{2}}].
\end{align}
Note that
\begin{align}\label{3.27'}
&(t_{n+1}-t)\leq2\tau(t),\ \ \  \ (t-t_{n-1})\leq2\tau(t), \ \ \   \tau (t_{n+1})\leq3\tau(t),\ \ \forall t\in[t_{n-1},t_{n+1}],\nn\\
&(t_{n+1}-t)\leq2\tau(t),\ \ \ \ (t-t_{n-1})\leq2\tau(t),\ \ \ \tau (t_{n+1})\leq3\tau(t),\ \ \forall t\in[t_{n-1},t_{n}],\nn\\
&(t_{n+1}-t)\leq2\tau(t),\ \ \ \ (t-t_{n-1})\leq2\tau(t),\ \ \ \tau(t_{n+1})\leq3\tau(t),\ \ \forall t\in[t_{n},t_{n+1}],
\end{align}
hold when $\Delta t\leq\frac{1}{2}$ and $n\geq2$. Summing \refe{3.22}\--\refe{3.27} from $n=2$ to $n=m$, respectively, and making use of Lemma \ref{L2.1} as well as Theorem \ref{T2.1}, we can derive \refe{3.17}\--\refe{3.21}.
\end{proof}

Next, we are ready to analyse the stability for the scheme \refe{3.1}\--\refe{3.6}.
\begin{theorem}\label{T3.2}
Suppose that the conditions of Lemma \ref{T3.1} and
$\Delta t\leq\min\{\frac{1}{2}, \frac{1}{\lambda_{1}}(1-\frac{1}{\delta_{0}})\}$ are satisfied, then
there holds
\begin{align}\label{3.28}
\|\textbf{u}_{h}^{m}\|_{1}^{2}+\mu\|\textbf{H}_{h}^{m}\|_{1}^{2}+\nu\|\textbf{u}_{h}^{m}\|_{2}^{2}+\sigma^{-1}\|\textbf{H}_{h}^{m}\|_{2}^{2}\leq \lambda_{3},
\end{align}
where $\delta_{0}$ is a positive constant more than $1$ , $\lambda_{1}$ is defined by \refe{3.40} and
\begin{align*}
&\lambda_{2}=e^{\delta_{0}\lambda_{1}T}(\lambda_{4}+\lambda_{5}),\ \ \lambda_{3}=2\kappa_{2}(1+\mu+\nu+\sigma^{-1})+2\lambda_{2},\\
&\lambda_{4}=(128c_{2}^{4}\nu^{-3}\kappa_{6}^{2}+\frac{98}{9}c_{4}^{4}\nu^{-1}\kappa_{6}^{2}+\frac{7}{3}c_{4}^{2}\nu^{-1}\kappa_{6}(\Delta t)^{\frac{1}{2}}+2c_{1}\nu^{-1}(\frac{7}{3}c_{4}^{2}+4c_{2}^{2})\kappa_{2}\nn\\
&\ \ \  \ \ \ +\frac{80}{3}c_{3}^{4}\mu^{3}\nu^{-2}\sigma(2+\frac{5}{6}\mu\sigma\nu)\kappa_{6}^{2}+\frac{1}{9}c_{5}^{4}\mu(\frac{245}{6}\mu^{2}\sigma+32\mu\nu^{-1})\kappa_{6}^{2}\nn\\
&\ \ \ \ \ \  +\frac{1}{3}c_{5}^{2}\mu(7\mu\sigma+4\nu^{-1})\kappa_{6}+2\mu(4c_{3}^{2}c_{1}\nu^{-1}+\frac{7}{3}c_{5}^{2}c_{1}\mu \sigma+\frac{16}{3}c_{5}^{2}c_{1}\nu^{-1})\kappa_{2})\nn\\
&\ \ \ \ \ \ \times2(1+\mu+\nu+\sigma^{-1})\kappa_{6},\\
&\lambda_{5}=2(1+\mu+\nu+\sigma^{-1})\kappa_{6}+(16\nu^{-1}+\frac{20}{3}\sigma+\frac{28}{3}+\frac{16}{3}\mu^{-1})\lambda_{0}.
\end{align*}

\begin{proof}
Due to \refe{2.12}, \refe{2.13}, \refe{2.20} and
\begin{align}\label{3.29}
&\|\textbf{u}_{h}^{m}\|_{1}^{2}+\mu\|\textbf{H}_{h}^{m}\|_{1}^{2}+\nu\|A_{1h}\textbf{u}_{h}^{m}\|_{0}^{2}+\sigma^{-1}\|A_{2h}\textbf{H}_{h}^{m}\|_{0}^{2}\leq 2(\|\textbf{u}_{h}(t_m)\|_{1}^{2}+\mu\|\textbf{H}_{h}(t_{m})\|_{1}^{2}\nn\\ &\ \ +\nu\|A_{1h}\textbf{u}_{h}(t_{m})\|_{0}^{2}+\sigma^{-1}\|A_{2h}\textbf{H}_{h}(t_{m})\|_{0}^{2})
 +2(\|e^{m}\|_{1}^{2}+\mu\|\xi^{m}\|_{1}^{2}+\nu\|e^{m}\|_{2}^{2}+\sigma^{-1}\|\xi^{m}\|_{2}^{2}),
\end{align}
 we only need to testify
\begin{align}\label{3.30}
\|e^{m}\|_{1}^{2}+\mu\|\xi^{m}\|_{1}^{2}+\nu\|e^{m}\|_{2}^{2}\Delta t+\sigma^{-1}\|\xi^{m}\|_{2}^{2}\Delta t\leq\lambda_{2}\Delta t.
\end{align}
We prove \refe{3.30} by mathematical induction. Apparently, \refe{3.30} is valid when $m=0$ according to \refe{2.11} and $(\textbf{u}_{h}^{0},\textbf{H}_{h}^{0})=(P_{h}\textbf{u}_{0},R_{0h}\textbf{H}_{0})$. From Lemma \ref{L3.1}, we can know that \refe{3.30} holds for $m=1, 2$. Let us assume that \refe{3.30} holds for $0\leq m \leq J-1$ $(J\geq 4)$. Next, we only need to confirm that \refe{3.30} is valid for $m=J$.
Noting
\begin{align}\label{3.31}
&(\textbf{a}-\textbf{b},\textbf{a})=\frac{1}{2}[(\textbf{a},\textbf{a})-(\textbf{b},\textbf{b})+(\textbf{a}-\textbf{b},\textbf{a}-\textbf{b})],\\
\label{3.32}
&(\mathbf{a}+\textbf{b},\textbf{a})=\frac{1}{2}[(\textbf{a},\textbf{a})-(\textbf{b},\textbf{b})+(\textbf{a}+\textbf{b},\textbf{a}+\textbf{b})].
\end{align}
Taking $\textbf{v}_{h}=4\Delta t A_{1h}e^{n+1}\in X_{0h}$, $q_{h}=0$ and $\textbf{B}_{h}=4\Delta t A_{2h}\xi^{n+1}\in W_{h}$ in \refe{3.9} and \refe{3.10}, respectively, and employing \refe{3.31} and \refe{3.32}, it leads to
\begin{align}\label{3.33}
&\|e^{n+1}\|_{1}^{2}-\|e^{n-1}\|_{1}^{2}+4\|d_{t}e^{n}\|_{1}^{2}(\Delta t)^{2}+\nu(\|e^{n+1}\|_{2}^{2}-\|e^{n-1}\|_{2}^{2}+4\|\tilde{e}^{n}\|_{2}^{2})\Delta t\nn\\&
\ \
=-4b(\textbf{u}_{h}^{n},\tilde{e}^{n},A_{1h}e^{n+1})\Delta t -4b(e^{n},\tilde{\textbf{u}}_{h}(t_{n}),A_{1h}e^{n+1})\Delta t -4\mu(\textbf{H}_{h}^{n}\times\operatorname{curl}\tilde{\xi}^{n},A_{1h}e^{n+1})\Delta t\nn\\
&\ \ \ \ \  -4\mu(\xi^{n}\times\operatorname{curl}\tilde{\textbf{H}}_{h}(t_{n}),A_{1h}e^{n+1})\Delta t +4(\textbf{E}^{n},A_{1h}e^{n+1})\Delta t\equiv \sum_{i=1}^{5}M_{i},\\
\label{3.34}
&\mu(\|\xi^{n+1}\|_{1}^{2}-\|\xi^{n-1}\|_{1}^{2}+4\|d_{t}\xi^{n}\|_{1}^{2}(\Delta t)^{2})+\sigma^{-1}(\|\xi^{n+1}\|_{2}^{2}-\|\xi^{n-1}\|_{2}^{2}+4\|\tilde{\xi}^{n}\|_{2}^{2})\Delta t\nn\\
&\ \ =4\mu(\tilde{\textbf{u}}_{h}(t_{n})\times\xi^{n},\operatorname{curl}A_{2h}\xi^{n+1})\Delta t+4\mu(\tilde{e}^{n}\times\textbf{H}_{h}^{n},\operatorname{curl}A_{2h}\xi^{n+1})\Delta t+4(\textbf{F}^{n},A_{2h}\xi^{n+1})\Delta t\nn\\
&\ \ \equiv\sum_{i=6}^{8}M_{i}.
\end{align}
Applying \refe{2.12}\--\refe{2.19}, $e^{n+1}=\tilde{e}^{n}+\Delta t d_{t}e^{n}$, $\xi^{n+1}=\tilde{\xi}^{n}+\Delta t d_{t}\xi^{n}$ and Young's inequality, it follows that
\begin{align*}
M_{1}&=-4b(\textbf{u}_{h}^{n},\tilde{e}^{n},A_{1h}e^{n+1})\Delta t\nn\\
&=4b(e^{n},\tilde{e}^{n},A_{1h}(\tilde{e}^{n}+\Delta t d_{t}e^{n}))\Delta t-4b(\textbf{u}_{h}(t_{n}),\tilde{e}^{n},A_{1h}(\tilde{e}^{n}+\Delta t d_{t}e^{n}))\Delta t\nn\\
&\leq 4 c_{2}\|e^{n}\|_{1}\|\tilde{e}^{n}\|_{1}^{\frac{1}{2}}\|\tilde{e}^{n}\|_{2}^{\frac{1}{2}}\|\tilde{e}^{n}\|_{2}\Delta t+4c_{2}c_{1}^{\frac{1}{2}}\|A_{1h}\textbf{u}_{h}(t_{n})\|_{0}\|\tilde{e}^{n}\|_{1}\|\tilde{e}^{n}\|_{2}\Delta t
\nn\\
&\ \ \ +2c_{4}(\|e^{n}\|_{2}\|\tilde{e}^{n}\|_{1}^{\frac{1}{2}}\|\tilde{e}^{n}\|_{2}^{\frac{1}{2}}+\|\tilde{e}^{n}\|_{2}\|e^{n}\|_{1}^{\frac{1}{2}}\|e^{n}\|_{2}^{\frac{1}{2}})\|d_{t}e^{n}\|_{1}(\Delta t)^{2}\nn\\
&\ \ \   +4c_{4}c_{1}^{\frac{1}{2}}\|A_{1h}\textbf{u}_{h}(t_{n})\|_{0}\|\tilde{e}^{n}\|_{2}\|d_{t}e^{n}\|_{1}(\Delta t)^{2}\nn\\
&\leq \frac{\nu}{2}\|\tilde{e}^{n}\|_{2}^{2}\Delta t +16c_{2}^{2}\nu^{-1}(\|e^{n}\|_{1}^{2}\|\tilde{e}^{n}\|_{1}\|\tilde{e}^{n}\|_{2}+c_{1}\|A_{1h}\textbf{u}_{h}(t_{n})\|_{0}^{2}\|\tilde{e}^{n}\|_{1}^{2})\Delta t+\frac{6}{7}\|d_{t}e^{n}\|_{1}^{2}(\Delta t)^{2}\nn\\
&\ \ \ +\frac{14}{3}c_{4}^{2}(\|e^{n}\|_{2}^{2}\|\tilde{e}^{n}\|_{1}\|\tilde{e}^{n}\|_{2}+\|\tilde{e}^{n}\|_{2}^{2}\|e^{n}\|_{1}\|e^{n}\|_{2}+2c_{1}\|A_{1h}\textbf{u}_{h}(t_{n})\|_{0}^{2}\|\tilde{e}^{n}\|_{2}^{2})(\Delta t)^{2}\nn\\
&\leq \nu\|\tilde{e}^{n}\|_{2}^{2}\Delta t+256c_{2}^{4}\nu^{-3}\|e^{n}\|_{1}^{4}\|\tilde{e}^{n}\|_{1}^{2}\Delta t+16c_{2}^{2}c_{1}\nu^{-1}\|A_{1h}\textbf{u}_{h}(t_{n})\|_{0}^{2}\|\tilde{e}^{n}\|_{1}^{2}\Delta t+\frac{6}{7}\|d_{t}e^{n}\|_{1}^{2}(\Delta t)^{2}\nn\\
&\ \ \   +\frac{196}{9}c_{4}^{4}\nu^{-1}\|e^{n}\|_{2}^{4}\|\tilde{e}^{n}\|_{1}^{2}(\Delta t)^{3}+\frac{14}{3}c_{4}^{2}\|\tilde{e}^{n}\|_{2}^{2}\|e^{n}\|_{1}\|e^{n}\|_{2}(\Delta t)^{2}+\frac{28}{3}c_{4}^{2}c_{1}\|A_{1h}\textbf{u}_{h}(t_{n})\|_{0}^{2}\|\tilde{e}^{n}\|_{2}^{2}(\Delta t)^{2},\\
M_{2}&=-4b(e^{n},\tilde{\textbf{u}}_{h}(t_{n}),A_{1h}e^{n+1})\Delta t\leq4|b(e^{n},\tilde{\textbf{u}}_{h}(t_{n}),A_{1h}(\tilde{e}^{n}+\Delta t d_{t}e^{n}))|\Delta t\nn\\
&\leq 4c_{2}c_{1}^{\frac{1}{2}}\|A_{1h}\tilde{\textbf{u}}_{h}(t_{n})\|_{0}\|e^{n}\|_{1}\|\tilde{e}^{n}\|_{2}\Delta t+4c_{4}c_{1}^{\frac{1}{2}}\|e^{n}\|_{2}\|A_{1h}\tilde{\textbf{u}}_{h}(t_{n})\|_{0}\|d_{t}e^{n}\|_{1}(\Delta t)^{2}\nn\\
&\leq \frac{\nu}{4}\|\tilde{e}^{n}\|_{2}^{2}\Delta t+16c_{2}^{2}c_{1}\nu^{-1}\|A_{1h}\tilde{\textbf{u}}_{h}(t_{n})\|_{0}^{2}\|e^{n}\|_{1}^{2}\Delta t+\frac{3}{7}\|d_{t}e^{n}\|_{1}^{2}(\Delta t)^{2}\nn\\
&\ \ \  +\frac{28}{3}c_{4}^{2}c_{1}\|e^{n}\|_{2}^{2}\|A_{1h}\tilde{\textbf{u}}_{h}(t_{n})\|_{0}^{2}(\Delta t)^{2},\\
M_{3}&=-4\mu(\textbf{H}_{h}^{n}\times\operatorname{curl}\tilde{\xi}^{n},A_{1h}e^{n+1})\Delta t \nn\\ &\leq4\mu|(\xi^{n}\times\operatorname{curl}\tilde{\xi}^{n},A_{1h}(\tilde{e}^{n}+\Delta t d_{t}e^{n}))-(\textbf{H}_{h}(t_{n})\times\operatorname{curl}\tilde{\xi}^{n},A_{1h}(\tilde{e}^{n}+\Delta t d_{t}e^{n}))|\Delta t\nn\\
&\leq 4c_{3}\mu\|\xi^{n}\|_{1}\|\tilde{\xi}^{n}\|_{1}^{\frac{1}{2}}\|\tilde{\xi}^{n}\|_{2}^{\frac{1}{2}}\|\tilde{e}^{n}\|_{2}\Delta t+4c_{3}c_{1}^{\frac{1}{2}}\mu\|A_{2h}\textbf{H}_{h}(t_{n})\|_{0}\|\tilde{\xi}^{n}\|_{1}\|\tilde{e}^{n}\|_{2}\Delta t\nn\\
&\ \ \ +2c_{5}\mu(\|\xi^{n}\|_{2}\|\tilde{\xi}^{n}\|_{1}^{\frac{1}{2}}\|\tilde{\xi}^{n}\|_{2}^{\frac{1}{2}}+\|\xi^{n}\|_{1}^{\frac{1}{2}}\|\xi^{n}\|_{2}^{\frac{1}{2}}\|\tilde{\xi}^{n}\|_{2})\|d_{t}e^{n}\|_{1}(\Delta t)^{2}\nn\\
&\ \ \  +4c_{5}c_{1}^{\frac{1}{2}}\mu\|A_{2h}\textbf{H}_{h}(t_{n})\|_{0}\|\tilde{\xi}^{n}\|_{2}\|d_{t}e^{n}\|_{1}(\Delta t)^{2}\nn\\
&\leq \frac{\nu}{2}\|\tilde{e}^{n}\|_{2}^{2}\Delta t+16c_{3}^{2}\mu^{2}\nu^{-1}\|\xi^{n}\|_{1}^{2}\|\tilde{\xi}^{n}\|_{1}\|\tilde{\xi}^{n}\|_{2}\Delta t+16c_{3}^{2}c_{1}\mu^{2}\nu^{-1}\|A_{2h}\textbf{H}_{h}(t_{n})\|_{0}^{2}\|\tilde{\xi}^{n}\|_{1}^{2}\Delta t\nn\\
&\ \ \  +\frac{14}{3}c_{5}^{2}\mu^{2}(\|\xi^{n}\|_{2}^{2}\|\tilde{\xi}^{n}\|_{1}\|\tilde{\xi}^{n}\|_{2} +\|\xi^{n}\|_{1}\|\xi^{n}\|_{2}\|\tilde{\xi}^{n}\|_{2}^{2}+2c_{1}\|A_{2h}\textbf{H}_{h}(t_{n})\|_{0}^{2}\|\tilde{\xi}^{n}\|_{2}^{2})(\Delta t)^{2}\nn\\
&\ \ \ +\frac{6}{7}\|d_{t}e^{n}\|_{1}^{2}(\Delta t)^{2}  \nn\\
&\leq \frac{\nu}{2}\|\tilde{e}^{n}\|_{2}^{2}\Delta t+\frac{6}{5}\sigma^{-1}\|\tilde{\xi}^{n}\|_{2}^{2}\Delta t+\frac{320}{3}c_{3}^{4}\mu^{4}\nu^{-2}\sigma\|\xi^{n}\|_{1}^{4}\|\tilde{\xi}^{n}\|_{1}^{2}\Delta t+16c_{3}^{2}c_{1}\mu^{2}\nu^{-1}\|A_{2h}\textbf{H}_{h}(t_{n})\|_{0}^{2}\|\tilde{\xi}^{n}\|_{1}^{2}\Delta t\nn\\
&\ \ \  +\frac{6}{7}\|d_{t}e^{n}\|_{1}^{2}(\Delta t)^{2}+\frac{245}{27}c_{5}^{4}\mu^{4}\sigma\|\xi^{n}\|_{2}^{4}\|\tilde{\xi}^{n}\|_{1}^{2}(\Delta t)^{3}+\frac{14}{3}c_{5}^{2}\mu^{2}\|\xi^{n}\|_{1}\|\xi^{n}\|_{2}\|\tilde{\xi}^{n}\|_{2}^{2}(\Delta t)^{2}\nn\\
&\ \ \ +\frac{28}{3}c_{5}^{2}c_{1}\mu^{2}\|A_{2h}\textbf{H}_{h}(t_{n})\|_{0}^{2}\|\tilde{\xi}^{n}\|_{2}^{2}(\Delta t)^{2},\\
M_{4}&=-4\mu(\xi^{n}\times\operatorname{curl}\tilde{\textbf{H}}_{h}(t_{n}),A_{1h}e^{n+1})\Delta t\leq4\mu|(\xi^{n}\times\operatorname{curl}\tilde{\textbf{H}}_{h}(t_{n}),A_{1h}(\tilde{e}^{n}+\Delta t d_{t}e^{n}))|\Delta t\nn\\
&\leq4c_{3}c_{1}^{\frac{1}{2}}\mu\|A_{2h}\tilde{\textbf{H}}_{h}(t_{n})\|_{0}\|\xi^{n}\|_{1}\|\tilde{e}^{n}\|_{2}\Delta t+4c_{5}c_{1}^{\frac{1}{2}}\mu\|\xi^{n}\|_{2}\|A_{2h}\tilde{\textbf{H}}_{h}(t_{n})\|_{0}\|d_{t}e^{n}\|_{1}(\Delta t)^{2}\nn\\
&\leq \frac{\nu}{4}\|\tilde{e}^{n}\|_{2}^{2}\Delta t+16c_{3}^{2}c_{1}\mu^{2}\nu^{-1}\|A_{2h}\tilde{\textbf{H}}_{h}(t_{n})\|_{0}^{2}\|\xi^{n}\|_{1}^{2}\Delta t
+\frac{3}{7}\|d_{t}e^{n}\|_{1}^{2}(\Delta t)^{2}\nn\\
&\ \ \ +\frac{28}{3}c_{5}^{2}c_{1}\mu^{2}\|\xi^{n}\|_{2}^{2}\|A_{2h}\tilde{\textbf{H}}_{h}(t_{n})\|_{0}^{2}(\Delta t)^{2},\\
M_{5}&=4(\textbf{E}^{n},A_{1h}e^{n+1})\Delta t\leq4|(\textbf{E}^{n},A_{1h}(\tilde{e}^{n}+\Delta t d_{t}e^{n}))|\Delta t\nn\\
&\leq \frac{\nu}{4}\|\tilde{e}^{n}\|_{2}^{2}\Delta t+16\nu^{-1}\|P_{h}\textbf{E}^{n}\|_{0}^{2}\Delta t+\frac{3}{7}\|d_{t}e^{n}\|_{1}^{2}(\Delta t)^{2}+\frac{28}{3}\|P_{h}\textbf{E}^{n}\|_{1}^{2}(\Delta t)^{2},\\
M_{6}&=4\mu(\tilde{\textbf{u}}_{h}(t_{n})\times\xi^{n},\operatorname{curl}A_{2h}\xi^{n+1})\Delta t\leq4\mu|(\tilde{\textbf{u}}_{h}(t_{n})\times\xi^{n},\operatorname{curl}A_{2h}(\tilde{\xi}^{n}+\Delta t d_{t}\xi^{n}))|\Delta t\nn\\
&\leq 4c_{3}c_{1}^{\frac{1}{2}}\mu\|A_{1h}\tilde{\textbf{u}}_{h}(t_{n})\|_{0}\|\xi^{n}\|_{1}\|\tilde{\xi}^{n}\|_{2}\Delta t+4c_{5}c_{1}^{\frac{1}{2}}\mu\|A_{1h}\tilde{\textbf{u}}_{h}(t_{n})\|_{0}\|\xi^{n}\|_{2}\|d_{t}\xi^{n}\|_{1}(\Delta t)^{2}\nn\\
&\leq\frac{3}{5}\sigma^{-1}\|\tilde{\xi}^{n}\|_{2}^{2}\Delta t+\frac{20}{3}c_{3}^{2}c_{1}\mu^{2}\sigma\|A_{1h}\tilde{\textbf{u}}_{h}(t_{n})\|_{0}^{2}\|\xi^{n}\|_{1}^{2}\Delta t+\frac{3}{4}\mu\|d_{t}\xi^{n}\|_{1}^{2}(\Delta t)^{2}\nn\\
&\ \ \ +\frac{16}{3}c_{5}^{2}c_{1}\mu\|A_{1h}\tilde{\textbf{u}}_{h}(t_{n})\|_{0}^{2}\|\xi^{n}\|_{2}^{2}(\Delta t)^{2},\\
M_{7}&=4\mu(\tilde{e}^{n}\times\textbf{H}_{h}^{n},\operatorname{curl}A_{2h}\xi^{n+1})\Delta t\nn\\
&\leq4\mu(|(\tilde{e}^{n}\times\xi^{n},\operatorname{curl}A_{2h}(\tilde{\xi}^{n}+\Delta t d_{t}\xi^{n}))|+|(\tilde{e}^{n}\times\textbf{H}_{h}(t_{n}),\operatorname{curl}A_{2h}(\tilde{\xi}^{n}+\Delta t d_{t}\xi^{n}))|)\Delta t\nn\\
&\leq 4c_{3}\mu\|\tilde{e}^{n}\|_{1}^{\frac{1}{2}}\|\tilde{e}^{n}\|_{2}^{\frac{1}{2}}\|\xi^{n}\|_{1}\|\tilde{\xi}^{n}\|_{2}\Delta t+4c_{5}c_{1}^{\frac{1}{2}}\mu\|A_{2h}\textbf{H}_{h}(t_{n})\|_{0}\|\tilde{\xi}^{n}\|_{1}\|\tilde{e}^{n}\|_{2}\Delta t\nn\\
&\ \ \ +2c_{5}\mu(\|\tilde{e}^{n}\|_{2}\|\xi^{n}\|_{1}^{\frac{1}{2}}\|\xi^{n}\|_{2}^{\frac{1}{2}}+\|\tilde{e}^{n}\|_{1}^{\frac{1}{2}}\|\tilde{e}^{n}\|_{2}^{\frac{1}{2}}\|\xi^{n}\|_{2})\|d_{t}\xi^{n}\|_{1}(\Delta t)^{2}\nn\\
&\ \ \ +4c_{5}c_{1}^{\frac{1}{2}}\mu\|\tilde{e}^{n}\|_{2}\|A_{2h}\textbf{H}_{h}(t_{n})\|_{0}\|d_{t}\xi^{n}\|_{1}(\Delta t)^{2}\nn\\
&\leq \frac{3}{5}\sigma^{-1}\|\tilde{\xi}^{n}\|_{2}^{2}\Delta t+\frac{20}{3}c_{3}^{2}\mu^{2}\sigma\|\tilde{e}^{n}\|_{1}\|\tilde{e}^{n}\|_{2}\|\xi^{n}\|_{1}^{2}\Delta t+16c_{5}^{2}c_{1}\mu^{2}\nu^{-1}\|A_{2h}\textbf{H}_{h}(t_{n})\|_{0}^{2}\|\tilde{\xi}^{n}\|_{1}^{2}\Delta t\nn\\
&\ \ \ +\frac{\nu}{4}\|\tilde{e}^{n}\|_{2}^{2}\Delta t+\frac{3}{2}\mu\|d_{t}\xi^{n}\|_{1}^{2}(\Delta t)^{2}+\frac{8}{3}c_{5}^{2}\mu(\|\tilde{e}^{n}\|_{2}^{2}\|\xi^{n}\|_{1}\|\xi^{n}\|_{2}+\|\tilde{e}^{n}\|_{1}\|\tilde{e}^{n}\|_{2}\|\xi^{n}\|_{2}^{2}\nn\\
&\ \ \ +2c_{1}\|\tilde{e}^{n}\|_{2}^{2}\|A_{2h}\textbf{H}_{h}(t_{n})\|_{0}^{2})(\Delta t)^{2} \nn\\
&\leq \frac{3}{5}\sigma^{-1}\|\tilde{\xi}^{n}\|_{2}^{2}\Delta t+\frac{400}{9}c_{3}^{4}\mu^{4}\sigma^{2}\nu^{-1}\|\tilde{e}^{n}\|_{1}^{2}\|\xi^{n}\|_{1}^{4}\Delta t+\frac{64}{9}c_{5}^{4}\mu^{2}\nu^{-1}\|\tilde{e}^{n}\|_{1}^{2}\|\xi^{n}\|_{2}^{4}(\Delta t)^{3}\nn\\
&\ \ \ +\frac{3}{4}\nu\|\tilde{e}^{n}\|_{2}^{2}\Delta t +16c_{5}^{2}c_{1}\mu^{2}\nu^{-1}\|A_{2h}\textbf{H}_{h}(t_{n})\|_{0}^{2}\|\tilde{\xi}^{n}\|_{1}^{2}\Delta t+\frac{8}{3}c_{5}^{2}\mu\|\tilde{e}^{n}\|_{2}^{2}\|\xi^{n}\|_{1}\|\xi^{n}\|_{2}(\Delta t)^{2}\nn\\
&\ \ \ +\frac{16}{3}c_{5}^{2}c_{1}\mu \|\tilde{e}^{n}\|_{2}^{2}\|A_{2h}\textbf{H}_{h}(t_{n})\|_{0}^{2}(\Delta t)^{2}+\frac{3}{2}\mu\|d_{t}\xi^{n}\|_{1}^{2}(\Delta t)^{2},\\
M_{8}&=4(\textbf{F}^{n},A_{2h}\xi^{n+1})\Delta t=4(\textbf{F}^{n},A_{2h}(\tilde{\xi}^{n}+\Delta t d_{t}\xi^{n}))\Delta t\nn\\
&\leq\frac{3}{5}\sigma^{-1}\|\tilde{\xi}^{n}\|_{2}^{2}\Delta t+\frac{20}{3}\sigma\|R_{0h}\textbf{F}^{n}\|_{0}^{2}\Delta t+\frac{3}{4}\mu\|d_{t}\xi^{n}\|_{1}^{2}(\Delta t)^{2}+\frac{16}{3}\mu^{-1}\|R_{0h}\textbf{F}^{n}\|_{1}^{2}(\Delta t)^{2}.
\end{align*}
Define
\begin{align*}
&a_{n}=\|e^{n+1}\|_{1}^{2}+\|e^{n}\|_{1}^{2}+\mu(\|\xi^{n+1}\|_{1}^{2}+\|\xi^{n}\|_{1}^{2})+\nu(\|e^{n+1}\|_{2}^{2}+\|e^{n}\|_{2}^{2})\Delta t+\sigma^{-1}(\|\xi^{n+1}\|_{2}^{2}+\|\xi^{n}\|_{2}^{2})\Delta t,\\
&b_{n}=\nu\|\tilde{e}^{n}\|_{2}^{2}+\|d_{t}e^{n}\|_{1}^{2}\Delta t+\sigma^{-1}\|\tilde{\xi}^{n}\|_{2}^{2}+\mu\|d_{t}\xi^{n}\|_{1}^{2}\Delta t,\\
&s_{n}=16\nu^{-1}\|P_{h}\textbf{E}^{n}\|_{0}^{2}\Delta t+\frac{28}{3}\|P_{h}\textbf{E}^{n}\|_{1}^{2}(\Delta t)^{2}+\frac{20}{3}\sigma\|R_{0h}\textbf{F}^{n}\|_{0}^{2}\Delta t+\frac{16}{3}\mu^{-1}\|R_{0h}\textbf{F}^{n}\|_{1}^{2}(\Delta t)^{2},\\
&\rho_{1}(n)=4c_{1}\nu^{-1}(4c_{2}^{2}+\frac{7}{3}c_{4}^{2}+\frac{5}{3}c_{3}^{2}\mu\nu\sigma+\frac{4}{3}c_{5}^{2}\mu\nu\sigma)\|A_{1h}\tilde{\textbf{u}}_{h}(t_{n})\|_{0}^{2}\nn\\
&\ \ \ \ \ \ \ \ \  \  +4c_{1}\mu(4c_{3}^{2}\nu^{-1}+\frac{7}{3}c_{5}^{2}\mu\sigma)\|A_{2h}\tilde{\textbf{H}}_{h}(t_n)\|_{0}^{2},\\
&\rho_{2}(n)=\frac{98}{9}c_{4}^{4}\nu^{-1}\|e^{n}\|_{2}^{4}(\Delta t)^{2}+2c_{1}\nu^{-1}(\frac{7}{3}c_{4}^{2}+4c_{2}^{2})\|A_{1h}\textbf{u}_{h}(t_{n})\|_{0}^{2}+\frac{7}{3}c_{4}^{2}\nu^{-1}\|e^{n}\|_{1}\|e^{n}\|_{2}\nn\\
&\ \ \ \ \ \ \ \ \  \  +128c_{2}^{4}\nu^{-3}\|e^{n}\|_{1}^{4}+\frac{80}{3}c_{3}^{4}\mu^{3}\nu^{-2}\sigma(2+\frac{5}{6}\mu\sigma\nu)\|\xi^{n}\|_{1}^{4}+\frac{1}{9}c_{5}^{4}\mu(\frac{245}{6}\mu^{2}\sigma+32\mu\nu^{-1})\|\xi^{n}\|_{2}^{4}(\Delta t)^{2}\nn\\
&\ \ \ \ \ \ \ \ \  \ +\frac{1}{3}c_{5}^{2}\mu(7\mu\sigma+4\nu^{-1})\|\xi^{n}\|_{1}\|\xi^{n}\|_{2}+2\mu c_{1}(4c_{3}^{2}\nu^{-1}+\frac{7}{3}c_{5}^{2}\mu \sigma+\frac{16}{3}c_{5}^{2}\nu^{-1})\|A_{2h}\textbf{H}_{h}(t_{n})\|_{0}^{2},
\end{align*}
taking sum of \refe{3.33} and \refe{3.34}, and combining $M_{1}$-$M_{8}$ with the equality, we conclude that
\begin{align}\label{3.35}
a_{n}-a_{n-1}+b_{n}\Delta t\leq (\rho_{1}(n)+\rho_{2}(n))a_{n}\Delta t+\rho_{2}(n)a_{n-1}\Delta t+s_{n}.
\end{align}
Summing \refe{3.35} from $n=2$ to $J-1$, it follows that
\begin{align}\label{3.36}
a_{J-1}+\Delta t\sum_{i=2}^{J-1}b_{i}\leq\Delta t\sum_{i=2}^{J-1}d_{i}a_{i}+\Delta t d_{1}a_{1}+a_{1}+\sum_{i=2}^{J-1}s_{i},
\end{align}
where
\begin{align}\label{3.37}
d_{i}=\left\{\begin{aligned}
&\rho_{2}(2), \ \ \ \ \ \ \ \ \ \ \ \ \ \ \ \ \ \ \ \ \ \ \ \ \ \ \ \ \  \ i=1,\\ &\rho_{1}(i)+\rho_{2}(i)+\rho_{2}(i+1),\ \ \ \ \ 2\leq i\leq J-2,\\ &\rho_{1}(J-1)+\rho_{2}(J-1),\ \ \ \  \ \ \ \ i=J-1.
\end{aligned}
\right.
\end{align}
Using the induction assumption, $\Delta t\leq\min\{\frac{1}{2}, \frac{1}{\lambda_{1}}(1-\frac{1}{\delta_{0}})\}, $ \refe{2.20}, Lemma \ref{L3.1}, \refe{3.20} and \refe{3.21}, we can get
\begin{align}\label{3.38}
\Delta t d_{1}a_{1}&\leq(128c_{2}^{4}\nu^{-3}\kappa_{6}^{2}+\frac{98}{9}c_{4}^{4}\nu^{-1}\kappa_{6}^{2}+\frac{7}{3}c_{4}^{2}\nu^{-1}\kappa_{6}(\Delta t)^{\frac{1}{2}}+2c_{1}\nu^{-1}(\frac{7}{3}c_{4}^{2}+4c_{2}^{2})\kappa_{2}\nn\\
&\ \ \ \  +\frac{80}{3}c_{3}^{4}\mu^{3}\nu^{-2}\sigma(2+\frac{5}{6}\mu\sigma\nu)\kappa_{6}^{2}+\frac{1}{9}c_{5}^{4}\mu(\frac{245}{6}\mu^{2}\sigma+32\mu\nu^{-1})\kappa_{6}^{2}\nn\\
&\ \ \ \ +\frac{1}{3}c_{5}^{2}\mu(7\mu\sigma+4\nu^{-1})\kappa_{6}+2\mu(4c_{3}^{2}c_{1}\nu^{-1}+\frac{7}{3}c_{5}^{2}c_{1}\mu \sigma+\frac{16}{3}c_{5}^{2}c_{1}\nu^{-1})\kappa_{2})\nn\\
&\ \ \ \ \times2(1+\mu+\nu+\sigma^{-1})\kappa_{6}(\Delta t)^{2}\equiv\lambda_{4}(\Delta t)^{2},\\\label{3.39}
a_{1}+\sum_{i=2}^{J-1}s_{i}&\leq
2(1+\mu+\nu+\sigma^{-1})\kappa_{6}\Delta t+(16\nu^{-1}+\frac{20}{3}\sigma+\frac{28}{3}+\frac{16}{3}\mu^{-1})\lambda_{0}\Delta t \equiv\lambda_{5}\Delta t,
\end{align}
where
\begin{align}\label{3.40}
    \lambda_{1}=\max_{2\leq i\leq N-1}\{d_{i}\}.
\end{align}
Employing Lemma \ref{L2.3}, plugging \refe{3.38} and \refe{3.39} into \refe{3.36}, there holds
\begin{align}\label{3.41}
a_{J-1}+\Delta t\sum_{i=2}^{J-1}b_{i}\leq e^{\delta_{0}\lambda_{1}T}(\lambda_{4}+\lambda_{5})\Delta t,
\end{align}
which demonstrates that \refe{3.30} is valid for $m=J$, where $\lambda_{2}=e^{\delta_{0}\lambda_{1}T}(\lambda_{4}+\lambda_{5})$.
Furthermore, an application of \refe{2.20}, \refe{3.29} and \refe{3.30} yields \refe{3.28}.
\end{proof}
\end{theorem}
	\section{Error estimation for the fully discrete CNLF scheme}
	\setcounter{equation}{0}
To derive the optimal error estimate, it is necessary firstly to prepare the following Lemma and Theorems.
\begin{lemma}\label{L4.1}
Assume that the conditions of Theorem \ref{T3.2} are satisfied and $\lambda_{6}\Delta t\leq1-\frac{1}{\delta_{0}}$, then there holds
\begin{align}\label{4.1}
&\nu(\|e^{n+1}\|_{0}^{2}+\|e^{n}\|_{0}^{2})+\sigma^{-1}(\|\xi^{n+1}\|_{0}^{2}+\|\xi^{n}\|_{0}^{2}) +\Delta t\sum_{n=2}^{m}\|d_{t}e^{n}\|_{-1}^{2}+\mu\Delta t\sum_{n=2}^{m}\|d_{t}\xi^{n}\|_{-1}^{2}\leq\lambda_{8}(\Delta t)^{2},
\end{align}
where
\begin{align*}
&\lambda_{7}=(1+\frac{10}{3}c_{2}^{2}c_{1}\nu^{-1}\kappa_{6}+2c_{1}(\frac{5}{3}c_{2}^{2}\nu^{-1}+ c_{3}^{2}\mu\sigma)\kappa_{2} +2c_{3}^{2}c_{1}\mu(\frac{5}{3}\mu\sigma+\nu^{-1})\kappa_{6}+\frac{10}{3}c_{3}^{2}c_{1}\mu^{2}\sigma\kappa_{2})\nn\\
&\ \ \ \ \ \ \times2(\nu+\sigma^{-1})\kappa_{6}+4(\frac{5}{3}+\mu^{-1})\lambda_{0},\ \ \ \lambda_{8}=e^{\delta_{0}\lambda_{6}T}\lambda_{7},
\end{align*}
 $\lambda_{6}$ is defined by \refe{4.6}.

\begin{proof}
Taking $\textbf{v}_{h}=4\Delta t A_{1h}^{\alpha}d_{t}e^{n}\in X_{0h}$, $q_{h}=0$ and $\textbf{B}_{h}=4\Delta t A_{2h}^{\alpha}d_{t}\xi^{n}\in W_{h}$ in \refe{3.9} and \refe{3.10} for $\alpha=-1, 0$, we deduce from adding them together that
\begin{align}\label{4.2}
&\nu(\|e^{n+1}\|_{\alpha+1}^{2}-\|e^{n-1}\|_{\alpha+1}^{2})+\sigma^{-1}(\|\xi^{n+1}\|_{\alpha+1}^{2}-\|\xi^{n-1}\|_{\alpha+1}^{2})+4\|d_{t}e^{n}\|_{\alpha}^{2}\Delta t+4\mu\|d_{t}\xi^{n}\|_{\alpha}^{2}\Delta t\nn\\
&\ \ =-4 b(\textbf{u}_{h}^{n},\tilde{e}^{n},A_{1h}^{\alpha}d_{t}e^{n})\Delta t-4 b(e^{n},\tilde{\textbf{u}}_{h}(t_{n}),A_{1h}^{\alpha}d_{t}e^{n})\Delta t-4\mu(\textbf{H}_{h}^{n}\times\operatorname{curl}\tilde{\xi}^{n},A_{1h}^{\alpha}d_{t}e^{n})\Delta t\nn\\
&\ \ \ \ \ -4\mu(\xi^{n}\times\operatorname{curl}\tilde{\textbf{H}}_{h}(t_{n}),A_{1h}^{\alpha}d_{t}e^{n})\Delta t+4(\textbf{E}^{n},A_{1h}^{\alpha}d_{t}e^{n})\Delta t+4\mu(\tilde{\textbf{u}}_{h}(t_{n})\times\xi^{n},\operatorname{curl}A_{2h}^{\alpha}d_{t}\xi^{n})\Delta t\nn\\
&\ \ \ \ \
+4\mu(\tilde{e}^{n}\times\textbf{H}_{h}^{n},\operatorname{curl}A_{2h}^{\alpha}d_{t}\xi^{n})\Delta t+4(\textbf{F}^{n},A_{2h}^{\alpha}d_{t}\xi^{n})\Delta t\equiv\sum_{i=9}^{12}M_{i}.
\end{align}
Employing \refe{2.14'}\--\refe{2.16} and Young's inequality, we have
\begin{align*}
M_{9}&=-4 b(\textbf{u}_{h}^{n},\tilde{e}^{n},A_{1h}^{\alpha}d_{t}e^{n})\Delta t-4 b(e^{n},\tilde{\textbf{u}}_{h}(t_{n}),A_{1h}^{\alpha}d_{t}e^{n})\Delta t\nn\\
&\leq4c_{2}c_{1}^{\frac{1}{2}}(\|\textbf{u}_{h}^{n}\|_{2}\|\tilde{e}^{n}\|_{\alpha+1}+\|e^{n}\|_{\alpha+1}\|\tilde{\textbf{u}}_{h}(t_{n})\|_{2})\|d_{t}e^{n}\|_{\alpha}\Delta t\nn\\
&\leq \frac{6}{5}\|d_{t}e^{n}\|_{\alpha}^{2}\Delta t+\frac{20}{3}c_{2}^{2}c_{1}(\|\textbf{u}_{h}^{n}\|_{2}^{2}\|\tilde{e}^{n}\|_{\alpha+1}^{2}+\|e^{n}\|_{\alpha+1}^{2}\|\tilde{\textbf{u}}_{h}(t_{n})\|_{2}^{2})\Delta t,\\
M_{10}&=-4\mu(\textbf{H}_{h}^{n}\times\operatorname{curl}\tilde{\xi}^{n},A_{1h}^{\alpha}d_{t}e^{n})\Delta t-4\mu(\xi^{n}\times\operatorname{curl}\tilde{\textbf{H}}_{h}(t_{n}),A_{1h}^{\alpha}d_{t}e^{n})\Delta t\nn\\
&\leq 4 c_{3}c_{1}^{\frac{1}{2}}\mu(\|\textbf{H}_{h}^{n}\|_{2}\|\tilde{\xi}^{n}\|_{\alpha+1}+\|\xi^{n}\|_{\alpha+1}\|\tilde{\textbf{H}}_{h}(t_{n})\|_{2})\|d_{t}e^{n}\|_{\alpha}\Delta t\nn\\
&\leq \frac{6}{5}\|d_{t}e^{n}\|_{\alpha}^{2}\Delta t+\frac{20}{3}c_{3}^{2}c_{1}\mu^{2}(\|\textbf{H}_{h}^{n}\|_{2}^{2}\|\tilde{\xi}^{n}\|_{\alpha+1}^{2}+\|\xi^{n}\|_{\alpha+1}^{2}\|\tilde{\textbf{H}}_{h}(t_{n})\|_{2}^{2})\Delta t,\\
M_{11}&=4\mu(\tilde{\textbf{u}}_{h}(t_{n})\times\xi^{n},\operatorname{curl}A_{2h}^{\alpha}d_{t}\xi^{n})\Delta t+4\mu(\tilde{e}^{n}\times\textbf{H}_{h}^{n},\operatorname{curl}A_{2h}^{\alpha}d_{t}\xi^{n})\Delta t\nn\\
&\leq4 c_{3}c_{1}^{\frac{1}{2}}\mu(\|\tilde{\textbf{u}}_{h}(t_{n})\|_{2}\|\xi^{n}\|_{\alpha+1}+\|\tilde{e}^{n}\|_{\alpha+1}\|\textbf{H}_{h}^{n}\|_{2})\|d_{t}\xi^{n}\|_{\alpha}\Delta t\nn\\
&\leq 2\mu\|d_{t}\xi^{n}\|_{\alpha}^{2}\Delta t+4 c_{3}^{2}c_{1}\mu(\|\tilde{\textbf{u}}_{h}(t_{n})\|_{2}^{2}\|\xi^{n}\|_{\alpha+1}^{2}+\|\tilde{e}^{n}\|_{\alpha+1}^{2}\|\textbf{H}_{h}^{n}\|_{2}^{2})\Delta t,\\
M_{12}&=4(\textbf{E}^{n},A_{1h}^{\alpha}d_{t}e^{n})\Delta t+4(\textbf{F}^{n},A_{2h}^{\alpha}d_{t}\xi^{n})\Delta t\nn\\
&\leq\frac{3}{5}\|d_{t}e^{n}\|_{\alpha}^{2}\Delta t+\mu\|d_{t}\xi^{n}\|_{\alpha}^{2}\Delta t+\frac{20}{3}\|P_{h}\textbf{E}^{n}\|_{\alpha}^{2}\Delta t+4\mu^{-1}\|R_{0h}\textbf{F}^{n}\|_{\alpha}^{2}\Delta t.
\end{align*}
Let
\begin{align*}
&a_{n}=\nu(\|e^{n+1}\|_{\alpha+1}^{2}+\|e^{n}\|_{\alpha+1}^{2})+\sigma^{-1}(\|\xi^{n+1}\|_{\alpha+1}^{2}+\|\xi^{n}\|_{\alpha+1}^{2}),\\
&b_{n}=\|d_{t}e^{n}\|_{\alpha}^{2}+\mu\|d_{t}\xi^{n}\|_{\alpha}^{2},\\
&\rho_{3}(n)=\frac{10}{3}c_{2}^{2}c_{1}\nu^{-1}\|\textbf{u}_{h}^{n}\|_{2}^{2}+2c_{1}(\frac{5}{3}c_{2}^{2}\nu^{-1}+ c_{3}^{2}\mu\sigma)\|\tilde{\textbf{u}}_{h}(t_{n})\|_{2}^{2}\nn\\
&\ \ \ \ \ \ \ \ \ \ +2c_{3}^{2}c_{1}\mu(\frac{5}{3}\mu\sigma+\nu^{-1})\|\textbf{H}_{h}^{n}\|_{2}^{2}+\frac{10}{3}c_{3}^{2}c_{1}\mu^{2}\sigma\|\tilde{\textbf{H}}_{h}(t_{n})\|_{2}^{2}.
\end{align*}
Combining $M_{9}$-$M_{12}$ with \refe{4.2} leads to
\begin{align}\label{4.3}
a_{n}-a_{n-1}+b_{n}\Delta t\leq\rho_{3}(n)(a_{n}+a_{n-1})\Delta t+\frac{20}{3}\|P_{h}\textbf{E}^{n}\|_{\alpha}^{2}\Delta t+4\mu^{-1}\|R_{0h}\textbf{F}^{n}\|_{\alpha}^{2}\Delta t.
\end{align}
It follows from summing up \refe{4.3} from $n=2$ to $m$ that
\begin{align}\label{4.4}
a_{m}+\Delta t\sum_{i=2}^{m}b_{i}\leq \Delta t\sum_{i=2}^{m}d_{i}a_{i}+(1+d_{1}\Delta t)a_{1}+\frac{20}{3}\Delta t\sum_{i=2}^{m}\|P_{h}\textbf{E}^{i}\|_{\alpha}^{2}+4\mu^{-1}\Delta t\sum_{i=2}^{m}\|R_{0h}\textbf{F}^{i}\|_{\alpha}^{2},
\end{align}
where
\begin{align}\label{4.5}
d_{i}=\left\{\begin{aligned}
&\rho_{3}(2), \ \ \ \ \ \ \ \ \ \ \ \ \ \ \ \ \ \ \ i=1,\\ &\rho_{3}(i)+\rho_{3}(i+1),\  \ \ \  2\leq i\leq m-1,\\ &\rho_{3}(m),\ \ \  \ \  \ \ \ \ \ \ \ \ \ \ \ \ \ i=m.
\end{aligned}
\right.
\end{align}
Due to Lemma \ref{L2.1} and Theorem \ref{T3.2}, there is a positive constant $\lambda_{6}$ such that
\begin{align}\label{4.6}
\lambda_{6}=\max_{2\leq i \leq N-1}\{d_{i}\}.
\end{align}
For $\alpha=-1$, using Lemma \ref{L2.1}, Lemma \ref{L3.1} and Lemma \ref{T3.1}, we have
\begin{align}\label{4.7}
&(1+d_{1}\Delta t)a_{1}+\frac{20}{3}\Delta t\sum_{i=2}^{m}\|P_{h}\textbf{E}^{i}\|_{\alpha}^{2}+4\mu^{-1}\Delta t\sum_{i=2}^{m}\|R_{0h}\textbf{F}^{i}\|_{\alpha}^{2}\nn\\
&\ \ \leq(1+\frac{10}{3}c_{2}^{2}c_{1}\nu^{-1}\kappa_{6}+2c_{1}(\frac{5}{3}c_{2}^{2}\nu^{-1}+ c_{3}^{2}\mu\sigma)\kappa_{2} +2c_{3}^{2}c_{1}\mu(\frac{5}{3}\mu\sigma+\nu^{-1})\kappa_{6}+\frac{10}{3}c_{3}^{2}c_{1}\mu^{2}\sigma\kappa_{2})\nn\\
&\ \ \ \ \ \times2(\nu+\sigma^{-1})\kappa_{6}(\Delta t)^{2}+4(\frac{5}{3}+\mu^{-1})\lambda_{0}(\Delta t)^{2}\equiv \lambda_{7}(\Delta t)^{2}.
\end{align}
Applying $\lambda_{6}\Delta t\leq1-\frac{1}{\delta_{0}}$, Lemma \ref{L2.3},  and \refe{4.7}  to \refe{4.4} emerges
\begin{align}\label{4.8}
a_{m}+\Delta t\sum_{i=2}^{m}b_{i}\leq e^{\delta_{0}\lambda_{6}T}\lambda_{7}(\Delta t)^{2}.
\end{align}
Therefore, \refe{4.1} is verified.
\end{proof}
\end{lemma}
\begin{theorem}\label{T4.1}
Suppose the conditions of Lemma \ref{L4.1} are established and $\lambda_{9}\Delta t\leq1-\frac{1}{\delta_{0}}$, then the following estimate exists
\begin{align}\label{4.9}
&\|e^{m+1}\|_{-2}^{2}+\|e^{m}\|_{-2}^{2}+\mu(\|\xi^{m+1}\|_{-2}^{2}+\|\xi^{m}\|_{-2}^{2})+\nu\Delta t\sum_{n=2}^{m}\|\tilde{e}^{n}\|_{-1}^{2} +\sigma^{-1}\Delta t\sum_{n=2}^{m}\|\tilde{\xi}^{n}\|_{-1}^{2}\leq\lambda_{12}(\Delta t)^{4},
\end{align}
where $\lambda_{9}$ is defined by \refe{4.14} and
\begin{align*}
\lambda_{10}&=2\kappa_{6}(1+\mu)(1+3c_{4}^{2}c_{1}\nu^{-1}\kappa_{6}+2c_{1}(c_{4}^{2}(\nu^{-1}+1)+2c_{5}^{2}(\mu\sigma+1))\kappa_{2} \nn\\
&\ \ \   +c_{5}^{2}c_{1}\mu(2\mu\sigma+3\nu^{-1})\kappa_{6}+4c_{4}^{2}c_{1}\mu(\mu\sigma+1)\kappa_{2}),\\
\lambda_{11}&=(1+\mu)\kappa_{6}+(6\nu^{-1}+4\sigma)\lambda_{0} +\lambda_{8},\\
\lambda_{12}&=e^{\delta_{0}\lambda_{9}T}(\lambda_{10}+\lambda_{11}).
\end{align*}
\begin{proof}
Taking $\textbf{v}_{h}=4\Delta t A_{1h}^{-2}\tilde{e}^{n}\in X_{0h}$, $q_{h}=0$ and $\textbf{B}_{h}=4\Delta t A_{2h}^{-2}\tilde{\xi}^{n}\in W_{h}$ in \refe{3.9} and \refe{3.10}, respectively, and summing up the two equalities, we can obtain
\begin{align}\label{4.10}
&\|e^{n+1}\|_{-2}^{2}-\|e^{n-1}\|_{-2}^{2}+\mu(\|\xi^{n+1}\|_{-2}^{2}-\|\xi^{n-1}\|_{-2}^{2})+4\nu\|\tilde{e}^{n}\|_{-1}^{2}\Delta t+4\sigma^{-1}\|\tilde{\xi}^{n}\|_{-1}^{2}\Delta t\nn\\
&=-4b(\textbf{u}_{h}^{n},\tilde{e}^{n},A_{1h}^{-2}\tilde{e}^{n})\Delta t-4b(e^{n},\tilde{\textbf{u}}_{h}(t_{n}),A_{1h}^{-2}\tilde{e}^{n})\Delta t-4\mu(\textbf{H}_{h}^{n}\times\operatorname{curl}\tilde{\xi}^{n},A_{1h}^{-2}\tilde{e}^{n})\Delta t\nn\\
&\ \ \ -4\mu(\xi^{n}\times\operatorname{curl}\tilde{\textbf{H}}_{h}(t_{n}),A_{1h}^{-2}\tilde{e}^{n})\Delta t+4\mu(\tilde{\textbf{u}}_{h}(t_{n})\times\xi^{n},\operatorname{curl}A_{2h}^{-2}\tilde{\xi}^{n})\Delta t+4\mu(\tilde{e}^{n}\times\textbf{H}_{h}^{n},\operatorname{curl}A_{2h}^{-2}\tilde{\xi}^{n})\Delta t\nn\\
&\ \ \   +4(\textbf{E}^{n},A_{1h}^{-2}\tilde{e}^{n})\Delta t+4(\textbf{F}^{n},A_{2h}^{-2}\tilde{\xi}^{n})\Delta t\equiv\sum_{i=13}^{19}M_{i}.
\end{align}
Using \refe{2.12}, \refe{2.13}, \refe{2.17}\--\refe{2.19}, $e^{n+1}=\tilde{e}^{n}+\Delta t d_{t}e^{n}$, $\xi^{n+1}=\tilde{\xi}^{n}+\Delta t d_{t}\xi^{n}$ and Young's inequality, we deduce that
\begin{align*}
M_{13}&=-4b(\textbf{u}_{h}^{n},\tilde{e}^{n},A_{1h}^{-2}\tilde{e}^{n})\Delta t\leq 4c_{4}c_{1}^{\frac{1}{2}}\|\textbf{u}_{h}^{n}\|_{2}\|\tilde{e}^{n}\|_{-2}\|\tilde{e}^{n}\|_{-1}\Delta t\nn\\
&\leq\frac{2}{3}\nu\|\tilde{e}^{n}\|_{-1}^{2}\Delta t+6c_{4}^{2}c_{1}\nu^{-1}
\|\textbf{u}_{h}^{n}\|_{2}^{2}\|\tilde{e}^{n}\|_{-2}^{2}\Delta t,\\
M_{14}&=-4b(e^{n},\tilde{\textbf{u}}_{h}(t_{n}),A_{1h}^{-2}\tilde{e}^{n})\Delta t\leq4|b(\tilde{e}^{n-1},\tilde{\textbf{u}}_{h}(t_{n}),A_{1h}^{-2}\tilde{e}^{n})+b(\Delta t d_{t}{e^{n-1}},\tilde{\textbf{u}}_{h}(t_{n}),A_{1h}^{-2}\tilde{e}^{n})|\Delta t\nn\\
&\leq 4c_{4}c_{1}^{\frac{1}{2}}(\|\tilde{e}^{n-1}\|_{-1}+\Delta t\|d_{t}e^{n-1}\|_{-1})\|A_{1h}\tilde{\textbf{u}}_{h}(t_{n})\|_{0}\|\tilde{e}^{n}\|_{-2}\Delta t\nn\\
&\leq\nu\|\tilde{e}^{n-1}\|_{-1}^{2}\Delta t+(\Delta t)^{3}\|d_{t}e^{n-1}\|_{-1}^{2}+4c_{4}^{2}c_{1}(\nu^{-1}+1)\|A_{1h}\tilde{\textbf{u}}_{h}(t_{n})\|_{0}^{2}\|\tilde{e}^{n}\|_{-2}^{2}\Delta t,\\
M_{15}&=-4\mu(\textbf{H}_{h}^{n}\times\operatorname{curl}\tilde{\xi}^{n},A_{1h}^{-2}\tilde{e}^{n})\Delta t\leq4c_{5}c_{1}^{\frac{1}{2}}\mu\|\textbf{H}_{h}^{n}\|_{2}\|\tilde{e}^{n}\|_{-2}\|\tilde{\xi}^{n}\|_{-1}\Delta t\nn\\
&\leq\sigma^{-1}\|\tilde{\xi}^{n}\|_{-1}^{2}\Delta t+4c_{5}^{2}c_{1}\mu^{2}\sigma\|\textbf{H}_{h}^{n}\|_{2}^{2}\|\tilde{e}^{n}\|_{-2}^{2}\Delta t,\\
M_{16}&=-4\mu(\xi^{n}\times\operatorname{curl}\tilde{\textbf{H}}_{h}(t_{n}),A_{1h}^{-2}\tilde{e}^{n})\Delta t\leq4\mu|((\tilde{\xi}^{n-1}+\Delta t d_{t}\xi^{n-1})\times\operatorname{curl}\tilde{\textbf{H}}_{h}(t_{n}),A_{1h}^{-2}\tilde{e}^{n})|\Delta t\nn\\
&\leq4c_{4}c_{1}^{\frac{1}{2}}\mu(\|\tilde{\xi}^{n-1}\|_{-1}+\Delta t\|d_{t}\xi^{n-1}\|_{-1})\|A_{2h}\tilde{\textbf{H}}_{h}(t_{n})\|_{0}\|\tilde{e}^{n}\|_{-2}\Delta t\nn\\
&\leq\frac{1}{2}\sigma^{-1}\|\tilde{\xi}^{n-1}\|_{-1}^{2}\Delta t+\frac{1}{2}\mu\|d_{t}\xi^{n-1}\|_{-1}^{2}(\Delta t)^{3}+8c_{4}^{2}c_{1}\mu(\mu\sigma+1)\|A_{2h}\tilde{\textbf{H}}_{h}(t_{n})\|_{0}^{2}\|\tilde{e}^{n}\|_{-2}^{2}\Delta t,\\
M_{17}&=4\mu(\tilde{\textbf{u}}_{h}(t_{n})\times\xi^{n},\operatorname{curl}A_{2h}^{-2}\tilde{\xi}^{n})\Delta t\leq4\mu|(\tilde{\textbf{u}}_{h}(t_{n})\times(\tilde{\xi}^{n-1}+\Delta t d_{t}\xi^{n-1}),\operatorname{curl}A_{2h}^{-2}\tilde{\xi}^{n})|\Delta t\nn\\
&\leq4c_{5}c_{1}^{\frac{1}{2}}\mu(\|\tilde{\xi}^{n-1}\|_{-1}+\Delta t\|d_{t}\xi^{n-1}\|_{-1})\|A_{1h}\tilde{\textbf{u}}_{h}(t_{n})\|_{0}\|\tilde{\xi}^{n}\|_{-2}\Delta t\nn\\
&\leq \frac{1}{2}\sigma^{-1}\|\tilde{\xi}^{n-1}\|_{-1}^{2}\Delta t+\frac{1}{2}\mu\|d_{t}\xi^{n-1}\|_{-1}^{2}(\Delta t)^{3}+8c_{5}^{2}c_{1}\mu(\mu\sigma+1)\|A_{1h}\tilde{\textbf{u}}_{h}(t_{n})\|_{0}^{2}\|\tilde{\xi}^{n}\|_{-2}^{2}\Delta t,\\
M_{18}&=4\mu(\tilde{e}^{n}\times\textbf{H}_{h}^{n},\operatorname{curl}A_{2h}^{-2}\tilde{\xi^{n}})\Delta t\leq4c_{5}c_{1}^{\frac{1}{2}}\mu\|\textbf{H}_{h}^{n}\|_{2}\|\tilde{e}^{n}\|_{-1}\|\tilde{\xi}^{n}\|_{-2}\Delta t\nn\\
&\leq\frac{2}{3}\nu\|\tilde{e}^{n}\|_{-1}^{2}\Delta t+6c_{5}^{2}c_{1}\mu^{2}\nu^{-1}\|\textbf{H}_{h}^{n}\|_{2}^{2}\|\tilde{\xi}^{n}\|_{-2}^{2}\Delta t,\\
M_{19}&=4(\textbf{E}^{n},A_{1h}^{-2}\tilde{e}^{n})\Delta t+4(\textbf{F}^{n},A_{2h}^{-2}\tilde{\xi}^{n})\Delta t\nn\\
&\leq\frac{2}{3}\nu\|\tilde{e}^{n}\|_{-1}^{2}\Delta t+6\nu^{-1}\|A_{1h}^{-\frac{3}{2}}P_{h}\textbf{E}^{n}\|_{0}^{2}\Delta t+\sigma^{-1}\|\tilde{\xi}^{n}\|_{-1}^{2}\Delta t+4\sigma\|A_{2h}^{-\frac{3}{2}}R_{0h}\textbf{F}^{n}\|_{0}^{2}\Delta t.
\end{align*}
Noting 
\begin{align*}
&a_{n}=\|e^{n+1}\|_{-2}^{2}+\|e^{n}\|_{-2}^{2}+\mu(\|\xi^{n+1}\|_{-2}^{2}+\|\xi^{n}\|_{-2}^{2}),\\
&b_{n}=\nu\|\tilde{e}^{n}\|_{-1}^{2}+\sigma^{-1}\|\tilde{\xi}^{n}\|_{-1}^{2},\\
&\rho_{4}(n)=3c_{4}^{2}c_{1}\nu^{-1}\|\textbf{u}_{h}^{n}\|_{2}^{2}+2c_{1}(c_{4}^{2}(\nu^{-1}+1)+2c_{5}^{2}(\mu\sigma+1))\|A_{1h}\tilde{\textbf{u}}_{h}(t_{n})\|_{0}^{2}\nn\\
&\ \ \ \ \ \ \ \ \ \ +c_{5}^{2}c_{1}\mu(2\mu\sigma+3\nu^{-1})\|\textbf{H}_{h}^{n}\|_{2}^{2}+4c_{4}^{2}c_{1}\mu(\mu\sigma+1)\|A_{2h}\tilde{\textbf{H}}_{h}(t_{n})\|_{0}^{2},\\
&s_{n}=6\nu^{-1}\|A_{1h}^{-\frac{3}{2}}P_{h}\textbf{E}^{n}\|_{0}^{2}+4\sigma\|A_{2h}^{-\frac{3}{2}}R_{0h}\textbf{F}^{n}\|_{0}^{2},
\end{align*}
and substituting $M_{13}$-$M_{19}$ into \refe{4.10}, it leads to
\begin{align}\label{4.11}
&a_{n}-a_{n-1}+2b_{n}\Delta t\leq \rho_{4}(n)(a_{n}+a_{n-1})\Delta t+b_{n-1}\Delta t +\|d_{t}e^{n-1}\|_{-1}^{2}(\Delta t)^{3}\nn\\
&\ \ \ +\mu\|d_{t}\xi^{n-1}\|_{-1}^{2}(\Delta t)^{3}+s_{n}\Delta t.
\end{align}
Summing up \refe{4.11} from $n=2$ to $m$, we have
\begin{align}\label{4.12}
a_{m}+\Delta t\sum_{i=2}^{m}b_{i}&\leq \Delta t\sum_{i=2}^{m}d_{i}a_{i}+(\Delta t)^{3}\sum_{i=2}^{m-1}(\|d_{t}e^{i}\|_{-1}^{2}+\mu\|d_{t}\xi^{i}\|_{-1}^{2})+\Delta t\sum_{i=2}^{m}s_{i}\nn\\
&\ \ \  +(1+d_{1}\Delta t)a_{1}+\|d_{t}e^{1}\|_{-1}^{2}(\Delta t)^{3}+\mu\|d_{t}\xi^{1}\|_{-1}^{2}(\Delta t)^{3},
\end{align}
where
\begin{align}\label{4.13}
d_{i}=\left\{\begin{aligned}
&\rho_{4}(2), \ \ \ \ \ \ \ \ \ \ \ \ \ \ \ \ \ \ \ i=1,\\ &\rho_{4}(i)+\rho_{4}(i+1),\ \ \ \ 2\leq i\leq m-1,\\ &\rho_{4}(m),\ \ \ \ \  \ \ \ \ \ \ \ \ \ \ \ \ \ i=m.
\end{aligned}
\right.
\end{align}
In view of Lemma \ref{L2.1} and Theorem \ref{T3.2}, there exists a positive constant $\lambda_{9}$ such that
\begin{align}\label{4.14}
\lambda_{9}=\max_{2\leq i\leq N-1}\{d_{i}\}.
\end{align}
 Since $e^{0}=\textbf{0}$, $\xi^{0}=\textbf{0}$, taking advantage of Lemma \ref{L2.1}, Lemma \ref{L3.1}, Lemma \ref{T3.1} and Lemma \ref{L4.1}, we can derive that
\begin{align}\label{4.15}
&(1+d_{1}\Delta t)a_{1}\leq2\kappa_{6}(1+\mu)(1+3c_{4}^{2}c_{1}\nu^{-1}\kappa_{6}+2c_{1}(c_{4}^{2}(\nu^{-1}+1)+2c_{5}^{2}(\mu\sigma+1))\kappa_{2} \nn\\
&\ \ \   +c_{5}^{2}c_{1}\mu(2\mu\sigma+3\nu^{-1})\kappa_{6}+4c_{4}^{2}c_{1}\mu(\mu\sigma+1)\kappa_{2})(\Delta t)^{4}\equiv\lambda_{10}(\Delta t)^{4},\\
\label{4.16}
&\|d_{t}e^{1}\|_{-1}^{2}(\Delta t)^{3}+\mu\|d_{t}\xi^{1}\|_{-1}^{2}(\Delta t)^{3} +\Delta t\sum_{i=2}^{m}s_{i}+(\Delta t)^{3}\sum_{i=2}^{m-1}(\|d_{t}e^{i}\|_{-1}^{2}+\mu\|d_{t}\xi^{i}\|_{-1}^{2})\nn\\ &\ \ \ \leq(1+\mu)\kappa_{6}(\Delta t)^{4}+(6\nu^{-1}+4\sigma)\lambda_{0}(\Delta t)^{4} +\lambda_{8}(\Delta t)^{4}\equiv\lambda_{11}(\Delta t)^{4}.
\end{align}
Combining \refe{4.12} with \refe{4.15} and \refe{4.16}, employing $\lambda_{9}\Delta t\leq1-\frac{1}{\delta_{0}}$ and Lemma \ref{L2.3}, we can get
\begin{align}\label{4.17}
a_{m}+\Delta t\sum_{i=2}^{m}b_{i}\leq e^{\delta_{0}\lambda_{9}T}(\lambda_{10}+\lambda_{11})(\Delta t)^{4},
\end{align}
which finishes the proof.
\end{proof}
\end{theorem}
\begin{theorem}\label{T4.2}
Under the conditions of Theorem \ref{T4.1} and $\lambda_{13}\Delta t\leq 1-\frac{1}{\delta_{0}}$, then there holds
\begin{align}\label{4.18}
&\tau(t_{m+1})(\|e^{m+1}\|_{-1}^{2}+\|e^{m}\|_{-1}^{2}+\mu(\|\xi^{m+1}\|_{-1}^{2}+\|\xi^{m}\|_{-1}^{2}))\nn\\
&\ \ +\Delta t\sum_{n=2}^{m}(\tau(t_{n+1})(\nu\|\tilde{e}^{n}\|_{0}^{2}+\sigma^{-1}\|\tilde{\xi}^{n}\|_{0}^{2}))\leq \lambda_{17}(\Delta t)^{4},
\end{align}
where $\lambda_{13}$ is defined by \refe{4.26} and
\begin{align*}
&\lambda_{14}=4\kappa_{6}(1+\mu)(4c_{4}^{2}c_{1}\nu^{-1}\kappa_{6}+4c_{1}(c_{4}^{2}\nu^{-1}+\frac{1}{3}c_{5}^{2}\mu\sigma)\kappa_{2}   +4c_{1}\mu\nu^{-1}(c_{5}^{2}+c_{3}^{2})\kappa_{6}+4c_{5}^{2}c_{1}\mu\nu^{-1}\kappa_{2}),\nn\\
&
\lambda_{15}=4\kappa_{6}(1+\mu)+8\lambda_{0}(\nu^{-1}+\frac{1}{3}\sigma),\ \lambda_{16} =4(\kappa_{6}(1+\mu)+\lambda_{8}+\lambda_{12}(\nu^{-1}+\mu\sigma)),\\
& \lambda_{17}=e^{\delta_{0}\lambda_{13}T}(2\lambda_{14}+\lambda_{15}+2\lambda_{16}).
\end{align*}

\begin{proof}
Taking $\textbf{v}_{h}=4\Delta t A_{1h}^{-1}\tilde{e}^{n}\in X_{0h}$, $q_{h}=0$ and $\textbf{B}_{h}=4\Delta t A_{2h}^{-1}\tilde{\xi}^{n}$ in \refe{3.9} and \refe{3.10}, respectively, and taking the sum of the two equalities, we can acquire
\begin{align}\label{4.19}
&\|e^{n+1}\|_{-1}^{2}-\|e^{n-1}\|_{-1}^{2}+\mu(\|\xi^{n+1}\|_{-1}^{2}-\|\xi^{n-1}\|_{-1}^{2})+4\nu\|\tilde{e}^{n}\|_{0}^{2}\Delta t+4\sigma^{-1}\|\tilde{\xi}^{n}\|_{0}^{2}\Delta t\nn\\
&\ \ =-4b(\textbf{u}_{h}^{n},\tilde{e}^{n},A_{1h}^{-1}\tilde{e}^{n})\Delta t-4b(e^{n},\tilde{\textbf{u}}_{h}(t_{n}),A_{1h}^{-1}\tilde{e}^{n})\Delta t-4\mu(\textbf{H}_{h}^{n}\times\operatorname{curl}\tilde{\xi}^{n},A_{1h}^{-1}\tilde{e}^{n})\Delta t\nn\\
&\ \ \ \ \ -4\mu(\xi^{n}\times\operatorname{curl}\tilde{\textbf{H}}_{h}(t_{n}),A_{1h}^{-1}\tilde{e}^{n})\Delta t+4\mu(\tilde{\textbf{u}}_{h}(t_{n})\times\xi^{n},\operatorname{curl}A_{2h}^{-1}\tilde{\xi}^{n})\Delta t\nn\\
&\ \ \ \ \   +4\mu(\tilde{e}^{n}\times\textbf{H}_{h}^{n},\operatorname{curl}A_{2h}^{-1}\tilde{\xi}^{n})\Delta t+4(\textbf{E}^{n},A_{1h}^{-1}\tilde{e}^{n})\Delta t+4(\textbf{F}^{n},A_{2h}^{-1}\tilde{\xi}^{n})\Delta t\equiv\sum_{i=20}^{23}M_{i}.
\end{align}
Applying \refe{2.12}, \refe{2.13}, \refe{2.15}, \refe{2.17}\--\refe{2.19}, and Young's inequality yields
\begin{align*}
M_{20}&=-4b(\textbf{u}_{h}^{n},\tilde{e}^{n},A_{1h}^{-1}\tilde{e}^{n})\Delta t-4b(e^{n},\tilde{\textbf{u}}_{h}(t_{n}),A_{1h}^{-1}\tilde{e}^{n})\Delta t \nn\\ &\leq4c_{4}c_{1}^{\frac{1}{2}}(\|\textbf{u}_{h}^{n}\|_{2}\|\tilde{e}^{n}\|_{-1}+\|A_{1h}\tilde{\textbf{u}}_{h}(t_{n})\|_{0}\|e^{n}\|_{-1})\|\tilde{e}^{n}\|_{0}\Delta t\nn\\
&\leq\nu\|\tilde{e}^{n}\|_{0}^{2}\Delta t+8c_{4}^{2}c_{1}\nu^{-1}(\|\textbf{u}_{h}^{n}\|_{2}^{2}\|\tilde{e}^{n}\|_{-1}^{2}+\|A_{1h}\tilde{\textbf{u}}_{h}(t_{n})\|_{0}^{2}\|e^{n}\|_{-1}^{2})\Delta t,\\
M_{21}&=-4\mu(\textbf{H}_{h}^{n}\times\operatorname{curl}\tilde{\xi}^{n},A_{1h}^{-1}\tilde{e}^{n})\Delta t-4\mu(\xi^{n}\times\operatorname{curl}\tilde{\textbf{H}}_{h}(t_{n}),A_{1h}^{-1}\tilde{e}^{n})\Delta t\nn\\
&\leq4c_{5}c_{1}^{\frac{1}{2}}\mu(\|\textbf{H}_{h}^{n}\|_{2}\|\tilde{\xi}^{n}\|_{-1}+\|A_{2h}\tilde{\textbf{H}}_{h}(t_{n})\|_{0}\|\xi^{n}\|_{-1})\|\tilde{e}^{n}\|_{0}\Delta t\nn\\
&\leq \nu\|\tilde{e}^{n}\|_{0}^{2}\Delta t+8c_{5}^{2}c_{1}\mu^{2}\nu^{-1}(\|\textbf{H}_{h}^{n}\|_{2}^{2}\|\tilde{\xi}^{n}\|_{-1}^{2}+\|A_{2h}\tilde{\textbf{H}}_{h}(t_{n})\|_{0}^{2}\|\xi^{n}\|_{-1}^{2})\Delta t,\\
M_{22}&=4\mu(\tilde{\textbf{u}}_{h}(t_{n})\times\xi^{n},\operatorname{curl}A_{2h}^{-1}\tilde{\xi}^{n})\Delta t+4\mu(\tilde{e}^{n}\times\textbf{H}_{h}^{n},\operatorname{curl}A_{2h}^{-1}\tilde{\xi}^{n})\Delta t\nn\\
&\leq4c_{5}c_{1}^{\frac{1}{2}}\mu\|A_{1h}\tilde{\textbf{u}}_{h}(t_{n})\|_{0}\|\xi^{n}\|_{-1}\|\tilde{\xi}^{n}\|_{0}\Delta t+4c_{3}c_{1}^{\frac{1}{2}}\mu\|\textbf{H}_{h}^{n}\|_{2}\|\tilde{\xi}^{n}\|_{-1}\|\tilde{e}^{n}\|_{0}\Delta t\nn\\
&\leq\frac{3}{2}\sigma^{-1}\|\tilde{\xi}^{n}\|_{0}^{2}\Delta t+\frac{\nu}{2}\|\tilde{e}^{n}\|_{0}^{2}\Delta t+\frac{8}{3}c_{5}^{2}c_{1}\mu^{2}\sigma\|A_{1h}\tilde{\textbf{u}}_{h}(t_{n})\|_{0}^{2}\|{\xi^{n}}\|_{-1}^{2}\Delta t\nn\\
&\ \ \ +8c_{3}^{2}c_{1}\mu^{2}\nu^{-1}
\|\textbf{H}_{h}^{n}\|_{2}^{2}\|\tilde{\xi}^{n}\|_{-1}^{2}\Delta t,\\
M_{23}&=4(\textbf{E}^{n},A_{1h}^{-1}\tilde{e}^{n})\Delta t+4(\textbf{F}^{n},A_{2h}^{-1}\tilde{\xi}^{n})\Delta t
\nn\\&\leq\frac{\nu}{2}\|\tilde{e}^{n}\|_{0}^{2}\Delta t+8\nu^{-1}\|P_{h}\textbf{E}^{n}\|_{-2}^{2}\Delta t+\frac{3}{2}\sigma^{-1}\|\tilde{\xi}^{n}\|_{0}^{2}\Delta t+\frac{8}{3}\sigma\|R_{0h}\textbf{F}^{n}\|_{-2}^{2}\Delta t.
\end{align*}
Denoting
\begin{align*}
&a_{n}=\|e^{n+1}\|_{-1}^{2}+\|e^{n}\|_{-1}^{2}+\mu(\|\xi^{n+1}\|_{-1}^{2}+\|\xi^{n}\|_{-1}^{2}),\\
&b_{n}=\nu\|\tilde{e}^{n}\|_{0}^{2}+\sigma^{-1}\|\tilde{\xi}^{n}\|_{0}^{2},\\
&\rho_{5}(n)=4c_{4}^{2}c_{1}\nu^{-1}\|\textbf{u}_{h}^{n}\|_{2}^{2}+4c_{1}(c_{4}^{2}\nu^{-1}+\frac{1}{3}c_{5}^{2}\mu\sigma)\|A_{1h}\tilde{\textbf{u}}_{h}(t_{n})\|_{0}^{2} \nn\\
&\ \ \ \ \ \ \ \ \ \  +4c_{1}\mu\nu^{-1}(c_{5}^{2}+c_{3}^{2})\|\textbf{H}_{h}^{n}\|_{2}^{2}+4c_{5}^{2}c_{1}\mu\nu^{-1}\|A_{2h}\tilde{\textbf{H}}_{h}(t_{n})\|_{0}^{2},
\end{align*}
substituting $M_{20}$-$M_{23}$ into \refe{4.19}, we can see
\begin{align}\label{4.20}
a_{n}-a_{n-1}+b_{n}\Delta t\leq \rho_{5}(n)(a_{n}+a_{n-1})\Delta t+8\nu^{-1}\|P_{h}\textbf{E}^{n}\|_{-2}^{2}\Delta t+\frac{8}{3}\sigma\|R_{0h}\textbf{F}^{n}\|_{-2}^{2}\Delta t.
\end{align}
Multiplying \refe{4.20} by $\tau(t_{n+1})$, noting that $\tau(t_{n+1})\leq\tau(t_{n})+\Delta t$, using $\lambda_{13}\Delta t\leq 1-\frac{1}{\delta_{0}}$ and taking the sum of the inequality from $n=2$ to $m$, we have
\begin{align}\label{4.21}
&\tau(t_{m+1})(\|e^{m+1}\|_{-1}^{2}+\|e^{m}\|_{-1}^{2}+\mu(\|\xi^{m+1}\|_{-1}^{2}+\|\xi^{m}\|_{-1}^{2}))+\Delta t\sum_{i=2}^{m}\tau(t_{i+1})b_{i}
\nn\\&\ \ \leq\tau(t_{2})a_{1}+2\Delta t\sum_{i=1}^{m}a_{i}+2\tau(t_{2})d_{1}a_{1}\Delta t+\Delta t\sum_{i=2}^{m}\tau(t_{i+1})d_{i}a_{i}\nn\\
&\ \ \ \ \ +8\nu^{-1}\Delta t\sum_{i=2}^{m}\tau(t_{i+1})\|P_{h}\textbf{E}^{i}\|_{-2}^{2} +\frac{8}{3}\sigma\Delta t\sum_{i=2}^{m}\tau(t_{i+1})\|R_{0h}\textbf{F}^{i}\|_{-2}^{2},
\end{align}
where exists a positive constant $\lambda_{13}$ from Lemma \ref{L2.1} and Theorem \ref{T3.2}, such that
\begin{align}\label{4.26}
\lambda_{13}=\max_{2\leq i\leq N-1}\{d_{i}\},
\end{align}
and
\begin{align}\label{4.22}
d_{i}=\left\{\begin{aligned}
&\rho_{5}(2), \ \ \ \ \ \ \ \ \ \ \ \ \ \ \ \ \ \ \ i=1,\\ &\rho_{5}(i)+\rho_{5}(i+1),\ \ \ \ 2\leq i\leq m-1,\\ &\rho_{5}(m),\ \ \ \ \  \ \ \ \ \ \ \ \ \ \ \ \ \ i=m.
\end{aligned}
\right.
\end{align}
Using Lemma \ref{L2.1}, Lemma \ref{L3.1} and Lemma \ref{T3.1}, we can get
\begin{align}\label{4.23}
&\tau(t_{2})d_{1}a_{1}\Delta t\leq4\kappa_{6}(1+\mu)(4c_{4}^{2}c_{1}\nu^{-1}\kappa_{6}+4c_{1}(c_{4}^{2}\nu^{-1}+\frac{1}{3}c_{5}^{2}\mu\sigma)\kappa_{2} \nn\\
&\ \ \   +4c_{1}\mu\nu^{-1}(c_{5}^{2}+c_{3}^{2})\kappa_{6}+4c_{5}^{2}c_{1}\mu\nu^{-1}\kappa_{2})(\Delta t)^{4}\equiv\lambda_{14}(\Delta t)^{4},\\
\label{4.24}
&\tau(t_{2})a_{1}+8\nu^{-1}\Delta t\sum_{i=2}^{m}\tau(t_{i+1})\|P_{h}\textbf{E}^{i}\|_{-2}^{2}+\frac{8}{3}\sigma\Delta t\sum_{i=2}^{m}\tau(t_{i+1})\|R_{0h}\textbf{F}^{i}\|_{-2}^{2}\nn\\
&\ \ \  \leq4\kappa_{6}(1+\mu)(\Delta t)^{4}+8\lambda_{0}(\nu^{-1}+\frac{1}{3}\sigma)(\Delta t)^{4}\equiv\lambda_{15}(\Delta t)^{4}.
\end{align}
Thanks to Lemma \ref{L3.1}, Lemma \ref{L4.1}, Theorem \ref{T4.1}, $e^{i+1}=\tilde{e}^{i}+\Delta t d_{t}e^{i}$, $\xi^{i+1}=\tilde{\xi}^{i}+\Delta t d_{t}\xi^{i}$, $e^{0}=\textbf{0}$ and $\xi^{0}=\textbf{0}$, there holds
\begin{align}\label{4.25}
\Delta t\sum_{i=1}^{m}a_{i}&\leq 2\Delta t\sum_{i=2}^{m}(\|\tilde{e}^{i}\|_{-1}^{2}+\|\tilde{e}^{i-1}\|_{-1}^{2})+2\mu\Delta t\sum_{i=2}^{m}(\|\tilde{\xi}^{i}\|_{-1}^{2}+\|\tilde{\xi}^{i-1}\|_{-1}^{2})+a_{1}\Delta t\nn\\
&\ \ \  +2(\Delta t)^{3}\sum_{i=2}^{m}(\|d_{t}e^{i}\|_{-1}^{2}+\|d_{t}e^{i-1}\|_{-1}^{2})+2\mu(\Delta t)^{3}\sum_{i=2}^{m}(\|d_{t}\xi^{i}\|_{-1}^{2}+\|d_{t}\xi^{i-1}\|_{-1}^{2})\nn\\
& \leq4\Delta t\sum_{i=1}^{m}(\|\tilde{e}^{i}\|_{-1}^{2}+\mu\|\tilde{\xi}^{i}\|_{-1}^{2})+4(\Delta t)^{3}\sum_{i=1}^{m}(\|d_{t}e^{i}\|_{-1}^{2}+\mu\|d_{t}\xi^{i}\|_{-1}^{2})+a_{1}\Delta t\nn\\
&\leq4(\kappa_{6}(1+\mu)+\lambda_{8}+\lambda_{12}(\nu^{-1}+\mu\sigma))(\Delta t)^{4}\equiv\lambda_{16}(\Delta t)^{4}.
\end{align}
 Combining \refe{4.21} with \refe{4.23}\--\refe{4.25}, applying $\lambda_{13}\Delta t\leq1-\frac{1}{\delta_{0}}$ and Lemma \ref{L2.3}, we obtain
\begin{align}\label{4.27}
&\tau(t_{m+1})(\|e^{m+1}\|_{-1}^{2}+\|e^{m}\|_{-1}^{2}+\mu(\|\xi^{m+1}\|_{-1}^{2}+\|\xi^{m}\|_{-1}^{2}))+\Delta t\sum_{i=2}^{m}\tau(t_{i+1})b_{i}\nn\\
&\ \ \leq e^{\delta_{0}\lambda_{13}T}(2\lambda_{14}+\lambda_{15}+2\lambda_{16})(\Delta t)^{4},
\end{align}
which achieves \refe{4.18}.
\end{proof}
\end{theorem}
\begin{theorem}\label{T4.3}
Assume that the conditions of Theorem \ref{T4.2} are satisfied, then there holds
\begin{align}\label{4.28}
&\tau^{2}(t_{m+1})(\|e^{m+1}\|_{0}^{2}+\|e^{m}\|_{0}^{2}+\mu(\|\xi^{m+1}\|_{0}^{2}+\|\xi^{m}\|_{0}^{2}))\nn\\
&\ \ +\nu\Delta t\sum_{n=2}^{m}\tau^{2}(t_{n+1})\|\tilde{e}^{n}\|_{1}^{2}+\sigma^{-1}\Delta t\sum_{n=2}^{m}\tau^{2}(t_{n+1})\|\tilde{\xi}^{n}\|_{1}^{2}\leq\lambda_{23}(\Delta t)^{4},
\end{align}
where $\lambda_{18}$ is defined by \refe{4.32'} and
\begin{align*}
\lambda_{19}&=e^{\lambda_{18}T}((8c_{1}\kappa_{2}\kappa_{6}(c_{2}^{2}\nu^{-1}+\frac{2}{3}c_{3}^{2}\mu\sigma+c_{3}^{2}\mu\nu^{-1})+2\kappa_{6})(1+\mu)+4(\nu^{-1}+\frac{2}{3}\sigma)\lambda_{0}),\\
\lambda_{20}&=12\kappa_{6}(\nu+\sigma^{-1})+4(\frac{5}{3}+\mu^{-1})\lambda_{0},\\
\lambda_{21}&=e^{\delta_{0}\lambda_{6}T}((3\kappa_{6}(\nu+\sigma^{-1})+4\lambda_{19}+4\lambda_{2}(\nu+\sigma^{-1}\mu^{-1}))+\lambda_{20}),\\
\lambda_{22}&=21(1+\lambda_{18})(1+\mu)\kappa_{6}+18 ((\nu^{-1}+\mu\sigma)\lambda_{17}+\lambda_{21})(1+\lambda_{18}),\\
\lambda_{23}&= e^{\lambda_{18}T}(8\kappa_{6}(1+\lambda_{18})(1+\mu)+\lambda_{22}+4(\nu^{-1}+\frac{2}{3}\sigma)\lambda_{0}).
\end{align*}

\begin{proof}
Taking $\textbf{v}_{h}=4\Delta t\tilde{e}^{n}\in X_{0h}$, $q_{h}=0$ and $\textbf{B}_{h}=4\Delta t\tilde{\xi}^{n}$ in \refe{3.9} and \refe{3.10}, respectively, and summing up the two equalities, we can deduce that
\begin{align}\label{4.29}
&\|e^{n+1}\|_{0}^{2}-\|e^{n-1}\|_{0}^{2}+\mu(\|\xi^{n+1}\|_{0}^{2}-\|\xi^{n-1}\|_{0}^{2})+4\nu\|\tilde{e}^{n}\|_{1}^{2}\Delta t+4\sigma^{-1}\|\tilde{\xi}^{n}\|_{1}^{2}\Delta t\nn\\
&\ \ =-4 b(\textbf{u}_{h}^{n},\tilde{e}^{n},\tilde{e}^{n})\Delta t-4b(e^{n},\tilde{\textbf{u}}_{h}(t_{n}),\tilde{e}^{n})\Delta t-4\mu(\textbf{H}_{h}^{n}\times\operatorname{curl}\tilde{\xi}^{n},\tilde{e}^{n})\Delta t\nn\\
&\ \ \ \ \ \ +4\mu(\tilde{e}^{n}\times \textbf{H}_{h}^{n},\operatorname{curl}\tilde{\xi}^{n})\Delta t-4\mu(\xi^{n}\times\operatorname{curl}\tilde{\textbf{H}}_{h}(t_{n}),\tilde{e}^{n})\Delta t \nn\\ &\ \ \ \ \ \ +4\mu(\tilde{\textbf{u}}_{h}(t_{n})\times\xi^{n},\operatorname{curl}\tilde{\xi}^{n})\Delta t+4(\textbf{E}^{n},\tilde{e}^{n})\Delta t+4(\textbf{F}^{n},\tilde{\xi}^{n})\Delta t\equiv\sum_{i=24}^{27}M_{i}.
\end{align}
By \refe{2.12}\--\refe{2.16} and Young's inequality, we can get
\begin{align*}
M_{24}&=-4 b(\textbf{u}_{h}^{n},\tilde{e}^{n},\tilde{e}^{n})\Delta t-4b(e^{n},\tilde{\textbf{u}}_{h}(t_{n}),\tilde{e}^{n})\Delta t\nn\\
&\leq 4c_{2}c_{1}^{\frac{1}{2}}\|e^{n}\|_{0}\|\tilde{\textbf{u}}_{h}(t_{n})\|_{2}\|\tilde{e}^{n}\|_{1}\Delta t\leq\nu\|\tilde{e}^{n}\|_{1}^{2}\Delta t+4c_{2}^{2}c_{1}\nu^{-1}\|e^{n}\|_{0}^{2}\|\tilde{\textbf{u}}_{h}(t_{n})\|_{2}^{2}\Delta t,\\
M_{25}&=-4\mu(\textbf{H}_{h}^{n}\times\operatorname{curl}\tilde{\xi}^{n},\tilde{e}^{n})\Delta t+4\mu(\tilde{e}^{n}\times \textbf{H}_{h}^{n},\operatorname{curl}\tilde{\xi}^{n})\Delta t=0,\\
M_{26}&=-4\mu(\xi^{n}\times\operatorname{curl}\tilde{\textbf{H}}_{h}(t_{n}),\tilde{e}^{n})\Delta t+4\mu(\tilde{\textbf{u}}_{h}(t_{n})\times\xi^{n},\operatorname{curl}\tilde{\xi}^{n})\Delta t\nn\\
&\leq 4c_{3}c_{1}^{\frac{1}{2}}\mu\|\xi^{n}\|_{0}\|\tilde{\textbf{H}}_{h}(t_{n})\|_{2}\|\tilde{e}^{n}\|_{1}\Delta t+4c_{3}c_{1}^{\frac{1}{2}}\mu\|\tilde{\textbf{u}}_{h}(t_{n})\|_{2}\|\xi^{n}\|_{0}\|\tilde{\xi}^{n}\|_{1}\Delta t\nn\\
&\leq \nu\|\tilde{e}^{n}\|_{1}^{2}\Delta t+4c_{3}^{2}c_{1}\mu^{2}\nu^{-1}\|\xi^{n}\|_{0}^{2}\|\tilde{\textbf{H}}_{h}(t_{n})\|_{2}^{2}\Delta t+\frac{3}{2}\sigma^{-1}\|\tilde{\xi}^{n}\|_{1}^{2}\Delta t+\frac{8}{3}c_{3}^{2}c_{1}\mu^{2}\sigma\|\tilde{\textbf{u}}_{h}(t_{n})\|_{2}^{2}\|\xi^{n}\|_{0}^{2}\Delta t,\\
M_{27}&=4(\textbf{E}^{n},\tilde{e}^{n})\Delta t+4(\textbf{F}^{n},\tilde{\xi}^{n})\Delta t\nn\\
&\leq \nu\|\tilde{e}^{n}\|_{1}^{2}\Delta t+4\nu^{-1}\|P_{h}\textbf{E}^{n}\|_{-1}^{2}\Delta t+\frac{3}{2}\sigma^{-1}\|\tilde{\xi}^{n}\|_{1}^{2}\Delta t+\frac{8}{3}\sigma\|R_{0h}\textbf{F}^{n}\|_{-1}^{2}\Delta t.
\end{align*}
Define
\begin{align*}
&a_{n}=\|e^{n+1}\|_{0}^{2}+\|e^{n}\|_{0}^{2}+\mu(\|\xi^{n+1}\|_{0}^{2}+\|\xi^{n}\|_{0}^{2}),\\
&b_{n}=\nu\|\tilde{e}^{n}\|_{1}^{2}+\sigma^{-1}\|\tilde{\xi}^{n}\|_{1}^{2},\\
&\rho_{6}(n)=4c_{1}(c_{2}^{2}\nu^{-1}+\frac{2}{3}c_{3}^{2}\mu\sigma)\|\tilde{\textbf{u}}_{h}(t_{n})\|_{2}^{2}+4c_{3}^{2}c_{1}\mu\nu^{-1}\|\tilde{\textbf{H}}_{h}(t_{n})\|_{2}^{2}.
\end{align*}
From $M_{24}$-$M_{27}$ and \refe{4.29}, we conclude that
\begin{align}\label{4.30}
a_{n}-a_{n-1}+ b_{n}\Delta t\leq \rho_{6}(n)a_{n-1}\Delta t+4\nu^{-1}\|P_{h}\textbf{E}^{n}\|_{-1}^{2}\Delta t+\frac{8}{3}\sigma\|R_{0h}\textbf{F}^{n}\|_{-1}^{2}\Delta t.
\end{align}
It follows from summing up \refe{4.30} from $n=2$ to $m$ that
\begin{align}\label{4.31}
a_{m}+\Delta t\sum_{n=2}^{m}b_{n}&\leq \rho_{6}(2)a_{1}\Delta t+a_{1}+\Delta t\sum_{n=2}^{m-1}\rho_{6}(n+1)a_{n}+4\nu^{-1}\Delta t\sum_{n=2}^{m}\|P_{h}\textbf{E}^{n}\|_{-1}^{2}\nn\\
&\ \ \ +\frac{8}{3}\sigma\Delta t\sum_{n=2}^{m}\|R_{0h}\textbf{F}^{n}\|_{-1}^{2}.
\end{align}
Combining Lemma \ref{L2.1}, Lemma \ref{L2.3}, Lemma \ref{L3.1} and Lemma \ref{T3.1} with \refe{4.31}, we can obtain
\begin{align}\label{4.32}
a_{m}+\Delta t\sum_{n=2}^{m}b_{n}&\leq e^{\lambda_{18}T}((8c_{1}\kappa_{2}\kappa_{6}(c_{2}^{2}\nu^{-1}+\frac{2}{3}c_{3}^{2}\mu\sigma+c_{3}^{2}\mu\nu^{-1})+2\kappa_{6})(1+\mu)(\Delta t)^{2}\nn\\
&\ \ \ +4(\nu^{-1}+\frac{2}{3}\sigma)\lambda_{0}(\Delta t)^{2})\equiv\lambda_{19}(\Delta t)^{2},
\end{align}
where exists a positive constant $\lambda_{18}$ due to Lemma \ref{L2.1} that
\begin{align}\label{4.32'}
\lambda_{18}=\max_{2\leq n\leq N}\{\rho_{6}(n)\}.
\end{align}
Next, the arithmetic formula
\begin{align}\label{4.33}
\Delta t\sum_{n=2}^{m}{\tau(t_{n+1})(\|d_{t}e^{n}\|_{0}^{2}+\mu\|d_{t}\xi^{n}\|_{0}^{2})},
\end{align}
 will be estimated.

Multiplying \refe{4.3} by $\tau(t_{n+1})$ on its both sides for $\alpha=0$, and summing up it from $n=2$ to $m$, we acquire
\begin{align}\label{4.34}
&\tau(t_{m+1})(\nu(\|e^{m+1}\|_{1}^{2}+\|e^{m}\|_{1}^{2})+\sigma^{-1}(\|\xi^{m+1}\|_{1}^{2}+\|\xi^{m}\|_{1}^{2}))+\Delta t\sum_{n=2}^{m} \tau(t_{n+1})(\|d_{t}e^{n}\|_{0}^{2}+\mu\|d_{t}\xi^{n}\|_{0}^{2})\nn\\
&\ \ \leq(\tau(t_{2})+2\tau(t_{2})\rho_{3}(2)\Delta t)(\nu(\|e^{2}\|_{1}^{2}+\|e^{1}\|_{1}^{2})+\sigma^{-1}(\|\xi^{2}\|_{1}^{2}+\|\xi^{1}\|_{1}^{2}))+\frac{20}{3}\Delta t\sum_{n=2}^{m}\tau(t_{n+1})\|P_{h}\textbf{E}^{n}\|_{0}^{2}\nn\\
&\ \ \ \ \
+4\mu^{-1}\Delta t\sum_{n=2}^{m}\tau(t_{n+1})\|R_{0h}\textbf{F}^{n}\|_{0}^{2}+2\Delta t\sum_{n=2}^{m+1}(\nu(\|e^{n}\|_{1}^{2}+\|e^{n-1}\|_{1}^{2})+\sigma^{-1}(\|\xi^{n}\|_{1}^{2}+\|\xi^{n-1}\|_{1}^{2}))\nn\\
&\ \ \ \ \
+\Delta t\sum_{n=2}^{m}\tau(t_{n+1})(\rho_{3}(n)+\rho_{3}(n+1))(\nu(\|e^{n+1}\|_{1}^{2}+\|e^{n}\|_{1}^{2})+\sigma^{-1}(\|\xi^{n+1}\|_{1}^{2}+\|\xi^{n}\|_{1}^{2})).
\end{align}
Applying Lemma \ref{L3.1}, Lemma \ref{T3.1}, \refe{3.41}, \refe{4.32} and $\lambda_{6}\Delta t\leq1-\frac{1}{\delta_{0}}$, we can derive that
\begin{align}\label{4.35}
&(\tau(t_{2})+2\tau(t_{2})\rho_{3}(2)\Delta t)(\nu(\|e^{2}\|_{1}^{2}+\|e^{1}\|_{1}^{2})+\sigma^{-1}(\|\xi^{2}\|_{1}^{2}+\|\xi^{1}\|_{1}^{2}))+\frac{20}{3}\Delta t\sum_{n=2}^{m}\tau(t_{n+1})\|P_{h}\textbf{E}^{n}\|_{0}^{2}\nn\\
&\ \
+4\mu^{-1}\Delta t\sum_{n=2}^{m}\tau(t_{n+1})\|R_{0h}\textbf{F}^{n}\|_{0}^{2}\leq12\kappa_{6}(\nu+\sigma^{-1})(\Delta t)^{2}+4(\frac{5}{3}+\mu^{-1})\lambda_{0}(\Delta t)^{2}
\equiv\lambda_{20}(\Delta t)^{2},\\
\label{4.36}
&\Delta t\sum_{n=2}^{m+1}(\nu(\|e^{n}\|_{1}^{2}+\|e^{n-1}\|_{1}^{2})+\sigma^{-1}(\|\xi^{n}\|_{1}^{2}+\|\xi^{n-1}\|_{1}^{2}))\nn\\
&\ \ \leq \nu\Delta t(2\|e^{2}\|_{1}^{2}+\|e^{1}\|_{1}^{2})+\sigma^{-1}\Delta t(2\|\xi^{2}\|_{1}^{2}+\|\xi^{1}\|_{1}^{2})+2\Delta t\sum_{n=3}^{m+1}(\nu\|e^{n}\|_{1}^{2}+\sigma^{-1}\|\xi^{n}\|_{1}^{2})\nn\\
&\ \ \leq \nu\Delta t(2\|e^{2}\|_{1}^{2}+\|e^{1}\|_{1}^{2})+\sigma^{-1}\Delta t(2\|\xi^{2}\|_{1}^{2}+\|\xi^{1}\|_{1}^{2})+4\Delta t\sum_{n=2}^{m}(\nu\|\tilde{e}^{n}\|_{1}^{2}+\sigma^{-1}\|\tilde{\xi}^{n}\|_{1}^{2})\nn\\
&\ \ \ \ \ +4(\Delta t)^{3}\sum_{n=2}^{m}(\nu\|d_{t}e^{n}\|_{1}^{2}+\sigma^{-1}\|d_{t}\xi^{n}\|_{1}^{2})\nn\\
&\ \ \leq (3\kappa_{6}(\nu+\sigma^{-1})+4\lambda_{19}+4\lambda_{2}(\nu+\sigma^{-1}\mu^{-1}))(\Delta t)^{2}.
\end{align}
Combining Lemma \ref{L2.3}, \refe{4.6}, \refe{4.35}, and \refe{4.36} with \refe{4.34} leads to
\begin{align}\label{4.37}
&\tau(t_{m+1})(\nu(\|e^{m+1}\|_{1}^{2}+\|e^{m}\|_{1}^{2})+\sigma^{-1}(\|\xi^{n+1}\|_{1}^{2}+\|\xi^{n}\|_{1}^{2}))+\Delta t\sum_{n=2}^{m} \tau(t_{n+1})(\|d_{t}e^{n}\|_{0}^{2}+\mu\|d_{t}\xi^{n}\|_{0}^{2})\nn\\
&\ \  \leq e^{\delta_{0}\lambda_{6}T}((3\kappa_{6}(\nu+\sigma^{-1})+4\lambda_{19}+4\lambda_{2}(\nu+\sigma^{-1}\mu^{-1}))+\lambda_{20})(\Delta t)^{2}\equiv\lambda_{21}(\Delta t)^{2}.
\end{align}
Multiplying \refe{4.30} by $\tau^{2}(t_{n+1})$ on its both sides, and summing up the inequality from $n=2$ to $m$, we can obtain
\begin{align}\label{4.38}
&\tau^{2}(t_{m+1})(\|e^{m+1}\|_{0}^{2}+\|e^{m}\|_{0}^{2}+\mu(\|\xi^{m+1}\|_{0}^{2}+\|\xi^{m}\|_{0}^{2}))+\Delta t\sum_{n=2}^{m}\tau^{2}(t_{n+1})(\nu\|\tilde{e}^{n}\|_{1}^{2}+\sigma^{-1}\|\tilde{\xi}^{n}\|_{1}^{2})\nn\\
&\ \ \leq\tau^{2}(t_{2})(\lambda_{18}+1)(\|e^{2}\|_{0}^{2}+\|e^{1}\|_{0}^{2}+\mu(\|\xi^{2}\|_{0}^{2}+\|\xi^{1}\|_{0}^{2}))\nn\\
&\ \ \ \ \ + (1+\lambda_{18}\Delta t)(\Delta t)^{2}\sum_{n=2}^{m}(\|e^{n}\|_{0}^{2}+\|e^{n-1}\|_{0}^{2}+\mu(\|\xi^{n}\|_{0}^{2}+\|\xi^{n-1}\|_{0}^{2}))\nn\\
&\ \ \ \ \ +2\Delta t\sum_{n=2}^{m}(1+\lambda_{18})\tau(t_{n})(\|e^{n}\|_{0}^{2}+\|e^{n-1}\|_{0}^{2}+\mu(\|\xi^{n}\|_{0}^{2}+\|\xi^{n-1}\|_{0}^{2}))\nn\\
&\ \ \ \ \ +\Delta t\sum_{n=2}^{m-1}\tau^{2}(t_{n+1})\lambda_{18}(\|e^{n+1}\|_{0}^{2}+\|e^{n}\|_{0}^{2}+\mu(\|\xi^{n+1}\|_{0}^{2}+\|\xi^{n}\|_{0}^{2}))\nn\\
&\ \ \ \ \ +4\nu^{-1}\Delta t\sum_{n=2}^{m}\tau^{2}(t_{n+1})\|P_{h}\textbf{E}^{n}\|_{-1}^{2}+\frac{8}{3}\sigma\Delta t\sum_{n=2}^{m}\tau^{2}(t_{n+1})\|R_{0h}\textbf{F}^{n}\|_{-1}^{2}.
\end{align}
By Lemma \ref{L3.1}, \refe{4.18}, \refe{4.37}, $\tau(t_{n+1})\leq\tau(t_{n})+\Delta t$, $e^{n+1}=\tilde{e}^{n}+\Delta t d_{t}e^{n}$ and $\xi^{n+1}=\tilde{\xi}^{n}+\Delta t d_{t}\xi^{n}$, we have
\begin{align}\label{4.39}
& (1+\lambda_{18}\Delta t) (\Delta t)^{2}\sum_{n=2}^{m}(\|e^{n}\|_{0}^{2}+\|e^{n-1}\|_{0}^{2}+\mu(\|\xi^{n}\|_{0}^{2}+\|\xi^{n-1}\|_{0}^{2}))\nn\\
&\ \  +2\Delta t\sum_{n=2}^{m}(1+\lambda_{18})\tau(t_{n})(\|e^{n}\|_{0}^{2}+\|e^{n-1}\|_{0}^{2}+\mu(\|\xi^{n}\|_{0}^{2}+\|\xi^{n-1}\|_{0}^{2}))\nn\\
&\ \ \leq  3\Delta t(1+\lambda_{18})\tau(t_{2})(\|e^{1}\|_{0}^{2}+\mu\|\xi^{1}\|_{0}^{2})+ 3\Delta t(1+\lambda_{18}) (\tau(t_{2})+\tau(t_{3}))
(\|e^{2}\|_{0}^{2}+\mu\|\xi^{2}\|_{0}^{2})\nn\\
&\ \ \ \ \  
+\sum_{n=3}^{m-1} 3\Delta t(1+\lambda_{18})(\tau(t_{n})+\tau(t_{n+1}))
(\|e^{n}\|_{0}^{2}+\mu\|\xi^{n}\|_{0}^{2})+ 3 \Delta t(1+\lambda_{18})\tau(t_{m})(\|e^{m}\|_{0}^{2}+\mu\|\xi^{m}\|_{0}^{2})\nn\\
&\ \ \leq 3 \Delta t(1+\lambda_{18})\tau(t_{2})(\|e^{1}\|_{0}^{2}+\mu\|\xi^{1}\|_{0}^{2})+3 \Delta t(1+\lambda_{18})(\tau(t_{2})+\tau(t_{3}))(\|e^{2}\|_{0}^{2}+\mu\|\xi^{2}\|_{0}^{2})\nn\\
&\ \ \ \ \  + 9(1+\lambda_{18}) \Delta t\sum_{n=3}^{m}\tau(t_{n})(\|e^{n}\|_{0}^{2}+\mu\|\xi^{n}\|_{0}^{2}) \nn\\
&\ \ \leq  21(1+\lambda_{18}) (1+\mu)\kappa_{6}(\Delta t)^{4}+ 9(1+\lambda_{18}) \Delta t\sum_{n=3}^{m}\tau(t_{n})(\|e^{n}\|_{0}^{2}+\mu\|\xi^{n}\|_{0}^{2}) \nn\\
&\ \ \leq  21(1+\lambda_{18}) (1+\mu)\kappa_{6}(\Delta t)^{4}+ 18(1+\lambda_{18}) \Delta t\sum_{n=2}^{m-1}\tau(t_{n+1})(\|\tilde{e}^{n}\|_{0}^{2}+(\Delta t)^{2}\|d_{t}e^{n}\|_{0}^{2}\nn\\
&\ \ \ \ \  +\mu\|\tilde{\xi}^{n}\|_{0}^{2}+\mu(\Delta t)^{2}\|d_{t}\xi^{n}\|_{0}^{2})\nn\\
&\ \ \leq 21(1+\lambda_{18}) (1+\mu)\kappa_{6}(\Delta t)^{4}+ 18 ((\nu^{-1}+\mu\sigma)\lambda_{17}+\lambda_{21}) (1+\lambda_{18}) (\Delta t)^{4}\equiv\lambda_{22}(\Delta t)^{4}.
\end{align}
Combining Lemma \ref{L2.3}, Lemma \ref{L3.1}, Lemma \ref{T3.1} and \refe{4.39} with \refe{4.38}, \refe{4.28} can be verified.
\end{proof}
\end{theorem}

\begin{theorem}\label{T4.4}
Suppose that the conditions of Theorem \ref{T4.3} are satisfied, then the following estimation holds
\begin{align}
\Delta t \sum_{n=2}^{m}\tau^{3}(t_{n+1})\|p_{h}(t_{n})-p_{h}^{n}\|_{0}^{2}\leq \lambda_{27}(\Delta t)^{4},
\end{align}
where
 \begin{align*}
&\lambda_{24}=8\kappa_{6}(\nu+\sigma^{-1})+6\lambda_{17}+6\lambda_{21}(\nu+\sigma^{-1}\mu^{-1}),
\\
&\lambda_{25}= e^{\delta_{0}\lambda_{6}T} (32\kappa_{6}(\nu+\sigma^{-1})+6\lambda_{24}+(\frac{20}{3}+\mu^{-1})\lambda_{0}),\\
&\lambda_{26}=\frac{1}{\beta_{1}^{2}}\lambda_{1}'(\lambda_{25}+\lambda_0+(c_{2}^{2}c_{1}\nu^{-1}+c_{3}^{2}c_{1}\mu^{2}\sigma)\kappa_{2}\lambda_{24}+(\nu^{-1}+c_{2}^{2}c_{1}\nu^{-2}\lambda_3+c_{3}^{2}c_{1}\mu^{2}\sigma^{2}\lambda_{3})\lambda_{23}),\\
&\lambda_{27}=(4\lambda_{26}+126\kappa_{2}),
\end{align*}
for $2\leq m \leq N-1$.
\end{theorem}
\begin{proof}
Multiplying \refe{4.3} by $\tau^{2}(t_{n+1})$ on both sides with $\alpha=-1$, and taking the sum of the inequality from $n=2$ to $m$, we obtain
\begin{align}\label{4.40'}
&\tau^{2}(t_{m+1})(\nu(\|e^{m+1}\|_{0}^{2}+\|e^{m}\|_{0}^{2})+\sigma^{-1}(\|\xi^{m+1}\|_{0}^{2}+\|\xi^{m}\|_{0}^{2})) +\Delta t\sum_{n=2}^{m} \tau^{2}(t_{n+1})(\|d_{t}e^{n}\|_{-1}^{2}+\mu\|d_{t}\xi^{n}\|_{-1}^{2})\nn\\
&\ \ \leq4\tau^{2}(t_{2})(\nu(\|e^{2}\|_{0}^{2}+\|e^{1}\|_{0}^{2})+\sigma^{-1}(\|\xi^{2}\|_{0}^{2}+\|\xi^{1}\|_{0}^{2}))\nn\\
&\ \ \ \ \ +4\Delta t\sum_{n=1}^{m}\tau(t_{n+1})(\nu(\|e^{n+1}\|_{0}^{2}+\|e^{n}\|_{0}^{2})+\sigma^{-1}(\|\xi^{n+1}\|_{0}^{2}+\|\xi^{n}\|_{0}^{2}))\nn\\
&\ \ \ \ \
+\Delta t\sum_{n=2}^{m}\tau^{2}(t_{n+1})(\rho_{3}(n)+\rho_{3}(n+1))(\nu(\|e^{n+1}\|_{0}^{2}+\|e^{n}\|_{0}^{2})+\sigma^{-1}(\|\xi^{n+1}\|_{0}^{2}+\|\xi^{n}\|_{0}^{2}))\nn\\
&\ \ \ \ \ +2(\Delta t)^2 \sum_{n=2}^{m}\tau(t_{n+1})\rho_{3}(n+1)(\nu(\|e^{n+1}\|_{0}^{2}+\|e^{n}\|_{0}^{2})+\sigma^{-1}(\|\xi^{n+1}\|_{0}^{2}+\|\xi^{n}\|_{0}^{2}))\nn\\
&\ \ \ \ \ +\frac{20}{3}\Delta t\sum_{n=2}^{m}\tau^{2}(t_{n+1})\|P_{h}\textbf{E}^{n}\|_{-1}^{2}
+4\mu^{-1}\Delta t\sum_{n=2}^{m}\tau^{2}(t_{n+1})\|R_{0h}\textbf{F}^{n}\|_{-1}^{2}.
\end{align}
Making use of Lemma \ref{L3.1}, Theorem \ref{T4.2}, \refe{4.37}, $\tau(t_{n+1})\leq\tau(t_{n})+\Delta t$, $e^{n+1}=\tilde{e}^{n}+\Delta t d_{t}e^{n}$ and $\xi^{n+1}=\tilde{\xi}^{n}+\Delta t d_{t}\xi^{n}$, we get
\begin{align}\label{0.1}
&\Delta t\sum_{n=1}^{m}\tau(t_{n+1})(\nu(\|e^{n+1}\|_{0}^{2}+\|e^{n}\|_{0}^{2})+\sigma^{-1}(\|\xi^{n+1}\|_{0}^{2}+\|\xi^{n}\|_{0}^{2}))\nn\\
&\ \ \leq\Delta t \tau(t_{2})(\nu(3\|e^{2}\|_{0}^{2}+\|e^{1}\|_{0}^{2})+\sigma^{-1}(3\|\xi^{2}\|_{0}^{2}+\|\xi^{1}\|_{0}^{2}))+3\Delta t\sum_{n=2}^{m}\tau(t_{n+1})(\nu\|e^{n+1}\|_{0}^{2}+\sigma^{-1}\|\xi^{n+1}\|_{0}^{2})
\nn\\ &\  \ \leq8\kappa_{6}(\nu+\sigma^{-1})(\Delta t)^{4}+6\Delta t\sum_{n=2}^{m}\tau(t_{n+1})(\nu\|\tilde{e}^{n}\|_{0}^{2}+\sigma^{-1}\|\tilde{\xi}^{n}\|_{0}^{2})\nn\\
&\ \ \ \ \ +6(\Delta t)^{3}\sum_{n=2}^{m}\tau(t_{n+1})(\nu\|d_{t}{e^{n}}\|_{0}^{2}+\sigma^{-1}\|d_{t}{\xi^{n}}\|_{0}^{2})\nn\\
&\ \ \leq 8\kappa_{6}(\nu+\sigma^{-1})(\Delta t)^{4}+6\lambda_{17}(\Delta t)^{4}+6\lambda_{21}(\nu+\sigma^{-1}\mu^{-1})(\Delta t)^{4}\equiv\lambda_{24}(\Delta t)^{4}.
\end{align}
Combining Lemma \ref{L2.3}, Lemma \ref{L3.1}, Lemma \ref{T3.1}, \refe{0.1} with \refe{4.40'}, it yields
\begin{align}\label{0.2}
&\tau^{2}(t_{m+1})(\nu(\|e^{m+1}\|_{0}^{2}+\|e^{m}\|_{0}^{2})+\sigma^{-1}(\|\xi^{m+1}\|_{0}^{2}+\|\xi^{m}\|_{0}^{2})) +\Delta t\sum_{n=2}^{m} \tau^{2}(t_{n+1})(\|d_{t}e^{n}\|_{-1}^{2}+\mu\|d_{t}\xi^{n}\|_{-1}^{2})\nn\\
& \ \ \leq e^{\delta_{0}\lambda_{6}T} (32\kappa_{6}(\nu+\sigma^{-1})+6\lambda_{24}+(\frac{20}{3}+\mu^{-1})\lambda_{0})(\Delta t)^{4}\equiv\lambda_{25}(\Delta t)^{4}.
\end{align}
It follows from \refe{2.8}, \refe{2.12}\--\refe{2.16} and \refe{3.9} that
\begin{align}\label{4.40}
\beta_{1}\|\eta^{n}\|_{0}&\leq \|d_{t}e^{n}\|_{-1}+\|\tilde{e}^{n}\|_{1}+c_{2}c_{1}^{\frac{1}{2}}\|\textbf{u}_{h}^{n}\|_{2}\|\tilde{e}^{n}\|_{1}+c_{2}c_{1}^{\frac{1}{2}}\|\tilde{\textbf{u}}_{h}(t_{n})\|_{2}\|e^{n}\|_{0}+c_{3}c_{1}^{\frac{1}{2}}\mu\|\textbf{H}_{h}^{n}\|_{2}\|\tilde{\xi}^{n}\|_{1}\nn\\
&\ \ \ +c_{3}c_{1}^{\frac{1}{2}}\mu\|\tilde{\textbf{H}}_{h}(t_{n})\|_{2}\|\xi^{n}\|_{0}+\|P_{h}\textbf{E}^{n}\|_{-1}.
\end{align}
 Taking the square of \refe{4.40} on both sides, multiplying it by $\tau^{3}(t_{n+1})$, it yields that
\begin{align}\label{4.41}
\tau^{3}(t_{n+1})\|\eta^{n}\|_{0}^{2}&\leq \frac{1}{\beta_{1}^{2}}\lambda_{1}'\tau^{2}(t_{n+1}) (\|d_{t}e^{n}\|_{-1}^{2}+\|\tilde{e}^{n}\|_{1}^{2}+c_{2}^{2}c_{1}\|\textbf{u}_{h}^{n}\|_{2}^{2}\|\tilde{e}^{n}\|_{1}^{2}+c_{2}^{2}c_{1}\|\tilde{\textbf{u}}_{h}(t_{n})\|_{2}^{2}\|e^{n}\|_{0}^{2}\nn\\
&\ \ \ +c_{3}^{2}c_{1}\mu^{2}\|\textbf{H}_{h}^{n}\|_{2}^{2}\|\tilde{\xi}^{n}\|_{1}^{2}+c_{3}^{2}c_{1}\mu^{2}\|\tilde{\textbf{H}}_{h}(t_{n})\|_{2}^{2}\|\xi^{n}\|_{0}^{2}+\|P_{h}\textbf{E}^{n}\|_{-1}^{2}),
\end{align}
 for all $2\leq n \leq N-1$.

 Multiplying \refe{4.41} by $\Delta t$, taking the sum of the inequality from $n=2$ to $n=m$, and using Lemma \ref{L2.1}, Lemma \ref{T3.1}, Theorem \ref{T3.2}, Theorem \ref{T4.3}, \refe{0.1} and \refe{0.2}, we arrive at
 \begin{align}\label{4.42'}
 &\Delta t \sum_{n=2}^{m}\tau^{3}(t_{n+1})\|\eta^{n}\|_{0}^{2}\leq \frac{1}{\beta_{1}^{2}}\lambda_{1}'\Delta t\sum_{n=2}^{m}\tau^{2}(t_{n+1}) (\|d_{t}e^{n}\|_{-1}^{2}+\|\tilde{e}^{n}\|_{1}^{2}+c_{2}^{2}c_{1}\|\textbf{u}_{h}^{n}\|_{2}^{2}\|\tilde{e}^{n}\|_{1}^{2}\nn\\
 &\ \ +c_{2}^{2}c_{1}\|\tilde{\textbf{u}}_{h}(t_{n})\|_{2}^{2}\|e^{n}\|_{0}^{2} +c_{3}^{2}c_{1}\mu^{2}\|\textbf{H}_{h}^{n}\|_{2}^{2}\|\tilde{\xi}^{n}\|_{1}^{2}+c_{3}^{2}c_{1}\mu^{2}\|\tilde{\textbf{H}}_{h}(t_{n})\|_{2}^{2}\|\xi^{n}\|_{0}^{2}+\|P_{h}\textbf{E}^{n}\|_{-1}^{2})\nn\\
 &\ \ \leq\frac{1}{\beta_{1}^{2}}\lambda_{1}'(\lambda_{25}+\lambda_0+(c_{2}^{2}c_{1}\nu^{-1}+c_{3}^{2}c_{1}\mu^{2}\sigma)\kappa_{2}\lambda_{24}+(\nu^{-1}+c_{2}^{2}c_{1}\nu^{-2}\lambda_3+c_{3}^{2}c_{1}\mu^{2}\sigma^{2}\lambda_{3})\lambda_{23})(\Delta t)^{4}\nn\\
 &\ \ \equiv\lambda_{26}(\Delta t)^{4},
 \end{align}
 where $\lambda_{1}'$ is a positive constant.

 Thanks to \refe{3.16}, using Cauchy-Schwartz inequality and triangle inequality yields
 \begin{align}\label{4.42}
 &\tau^{3}(t_{n+1})\|p_{h}(t_{n})-p_{h}^{n}\|_{0}^{2}\nn\\
 &\ \ =\tau^{3}(t_{n+1})\|\eta^{n}+p_{h}(t_{n})-\frac{p_{h}(t_{n+1})+p_{h}(t_{n-1})}{2}-\frac{1}{4\Delta t}\int_{t_{n-1}}^{t_{n+1}}{(t-t_{n+1})(t-t_{n-1})p_{htt}}dt\|_{0}^{2}\nn\\
 &\ \ \leq4\tau^{3}(t_{n+1})\|\eta^{n}\|_{0}^{2}+\tau^{3}(t_{n+1})\|p_{h}(t_{n+1})-2p_{h}(t_{n})+p_{h}(t_{n-1})\|_{0}^{2}\nn\\
 &\ \ \ \ \ +\frac{1}{4\Delta t}\tau^{3}(t_{n+1})\int_{t_{n-1}}^{t_{n+1}}(t-t_{n+1})^2(t-t_{n-1})^2\|p_{htt}(t)\|_{0}^{2}dt\nn\\
 &\ \ \leq 4\tau^{3}(t_{n+1})\|\eta^{n}\|_{0}^{2}+\frac{2}{3}(\Delta t)^{3}\tau^{3}(t_{n+1})\int_{t_{n-1}}^{t_{n+1}}\|p_{htt}\|_{0}^{2}dt+4(\Delta t)^{3}\tau^{3}(t_{n+1})\int_{t_{n-1}}^{t_{n+1}}\|p_{htt}\|_{0}^{2}dt,
 \end{align}
 here, we have used that  \cite{sj, 42}
 $$\|p(t_{n+1})-2p(t_{n})+p(t_{n-1})\|_{0}^{2}\leq\frac{2}{3}(\Delta t)^{3}\int_{t_{n-1}}^{t_{n+1}}\|p_{htt}\|_{0}^{2}dt.$$

Multiplying \refe{4.42} by $\Delta t$, taking sum of the inequality from $n=2$ to $n=m$, and using \refe{2.25}, \refe{3.27'}, \refe{4.42'} and \refe{4.42}, we can observe
\begin{align}\label{4.43'}
\Delta t\sum_{n=2}^{m}\tau^{3}(t_{n+1})\|p_{h}(t_{n})-p_{h}^{n}\|_{0}^{2}dt\leq (4\lambda_{26}+126\kappa_{2})(\Delta t)^{4},
\end{align}
which completes the proof.
\end{proof}

 Combining Lemma \ref{L3.1}, Theorem \ref{T3.2}, Theorem \ref{T4.3}, Theorem \ref{T4.4} and \refe{4.37} with Theorem \ref{T2.2}, the Theorem \ref{T4.5} is given as follows.
 \begin{theorem}\label{T4.5}
  Assume  that the conditions of Theorem \ref{T4.3} are satisfied, then the fully discrete CNLF scheme \refe{3.1}\--\refe{3.2} is almost unconditional stable and has the following error estimates
 \begin{align}\label{4.43}
&\|\textbf{u}_{h}^{m}\|_{1}^{2}+\|\textbf{H}_{h}^{m}\|_{1}^{2}+\|\textbf{u}_{h}^{m}\|_{2}^{2}+\|\textbf{H}_{h}^{m}\|_{2}^{2}\leq \frac{\lambda_{3}}{\min\{1,\nu,\mu, \sigma^{-1}\}},\\
\label{4.44}
&\|\textbf{u}(t_{m})-\textbf{u}_{h}^{m}\|_{0}+\|\textbf{H}(t_{m})-\textbf{H}_{h}^{m}\|_{0}\leq\tau(t_{m})^{-1}\left(\left(\frac{\lambda_{23}}{\min\{1,\mu\}}\right)^{\frac{1}{2}}(\Delta t)^{2}+\kappa_{4}^{\frac{1}{2}}h^{2}\right),\\
\label{4.44'}
&\|\nabla(\textbf{u}(t_{m})-\textbf{u}_{h}^{m})\|_{0}+\|\nabla(\textbf{B}(t_{m})-\textbf{B}_{h}^{m})\|_{0}\leq\tau(t_{m})^{-\frac{1}{2}}\left(\left(\frac{\lambda_{21}}{\min\{\nu,\sigma^{-1}\}}\right)^{\frac{1}{2}}\Delta t+\kappa_{4}^{\frac{1}{2}}h\right),
\end{align}
for all $2\leq m \leq N$. Furthermore, there holds
\begin{align}
\label{4.45}
\Delta t \sum_{n=2}^{m}\tau^{3}(t_{n})\|p(t_{n})-p_{h}^{n}\|_{0}^{2}\leq 2\lambda_{27}(\Delta t)^{4}+2\kappa_{4}Th^{4},
\end{align}
for $2\leq m \leq N-1.$
 \end{theorem}
	\section{Numerical examples}
In this section, we show two numerical examples to illustrate the effectiveness of the CNLF algorithm. The finite element spaces are chosen as $P1b$-$P1$-$P1b$ finite element spaces, which satisfy the LBB condition.
	\subsection{2D problems with the smooth solution }
	In the first place, we consider two-dimensional non-stationary incompressible MHD equations with exact solution as follows
	\begin{align*}
	 u_1&=(y+y^4)*\cos(t),\\
	 u_2&=(x+x^2)*\cos(t),\\
	 p&=(2.0*x-1.0)*(2.0*y-1.0)*\cos(t),\\
	 b_1&=(\sin(y)+y) *\cos(t),\\
	 b_2&=(\sin(x)+x^2)*\cos(t).
	\end{align*}
	The computation is carried out by choosing $\tau=1/10h$ and $h = 1/n$, $n = 8,
	16, 24, 32, 40$. The code was implemented by using the software package FreeFEM++ \cite{HPO}.
	
	Table \ref{Ta1} displays the numerical results of the model with $\nu=1.0$, $\mu=1.0$ and $\sigma=1.0$. From the numerical results, we find that the errors change small as the space step go small, and the convergence orders are optimal. Table \ref{Ta2} shows the numerical results of the model with $\nu=0.01$, $\mu=1.0$ and $\sigma=100$.
	Table \ref{Ta3} presents the numerical results of the model with $\nu=0.001$, $\mu=1.0$ and $\sigma=1000$, and Table \ref{Ta4} is the numerical results of the model with $\nu=0.002$, $\mu=1.0$ and $\sigma=2000$. The numerical results are as alike as before, the convergence orders are also optimal, and the numerical results are consistent with our theoretical analysis, and our method is effective.
	\begin{table}[ht!]
		\tabcolsep 0pt \caption{The numerical results of $\nu=1.0$, $\mu=1.0$, $\sigma=1.0$ and $T=1.0$ for different $h$.}\label{Ta1} \vspace*{-18pt}
		\begin{center}
		\def\temptablewidth{1.0\textwidth}
		{\rule{\temptablewidth}{1pt}}
		\begin{tabular*}{\temptablewidth}{@{\extracolsep{\fill}}llllll}
		$\footnotesize1/h$ & $\footnotesize\|\textbf{u}_h-\textbf{u}\|_0$ & $\footnotesize\|\nabla(\textbf{u}_{h}-\textbf{u})\|_0$ &  $\footnotesize\|p_h^1-p\|_0$ & $\footnotesize\|\textbf{B}_h-\textbf{B}\|_0$ & $\footnotesize\|\nabla(\textbf{B}_{h}-\textbf{B})\|_0$ \\
		\hline
		8  &       0.00427603   &    0.103026    &   0.0256856      &      0.00100479  &     0.0295769  \\
		16 &       0.00106969   &    0.0515209    &   0.0076434      &      0.000239404  &     0.0147887    \\
		24 &       0.000475199   &    0.0343395    &   0.0037971      &      0.000102137  &     0.00985905     \\
		32  &       0.000267202   &    0.0257507    &   0.00232872      &      5.56424e-05  &     0.00739424    \\
		40 &       0.00017096   &    0.0205985    &   0.00160141      &      3.48671e-05  &     0.00591537    \\
	  \hline
	  $\footnotesize 1/h$ & $\footnotesize\textbf{u}_{\text{order}L^2}$ &  $\footnotesize\textbf{u}_{\text{order}H^1}$ & $\footnotesize p_{\text{order}L^2}$&  $\footnotesize\textbf{B}_{\text{order}L^2}$ & $\footnotesize\textbf{B}_{\text{order}H^1}$ \\
	  \hline
	  8  &     /  & / & /  &  /  &/   \\
	  16    &   1.99908   &    0.999782    &   1.74867   &    2.06938  &     0.999978    \\
      24    &   2.00113   &    1.00055    &   1.72544   &    2.10088  &     1.00002    \\
      32    &   2.00126   &    1.00053    &   1.69951   &    2.11126  &     1.00002    \\
      40    &   2.00129   &    1.00046    &   1.678   &    2.09462  &     1.00002    \\
			   \end{tabular*}
	 {\rule{\temptablewidth}{1pt}}
	\end{center}
	\end{table}

	
	\begin{table}[ht!]
		\tabcolsep 0pt \caption{The numerical results of $\nu=0.01$, $\mu=1.0$, $\sigma=100$ and $T=1.0$ for different $h$.}\label{Ta2} \vspace*{-18pt}
		\begin{center}
		\def\temptablewidth{1.0\textwidth}
		{\rule{\temptablewidth}{1pt}}
		\begin{tabular*}{\temptablewidth}{@{\extracolsep{\fill}}llllll}
		$\footnotesize1/h$ & $\footnotesize\|\textbf{u}_h-\textbf{u}\|_0$ & $\footnotesize\|\nabla(\textbf{u}_{h}-\textbf{u})\|_0$ &  $\footnotesize\|p_h^1-p\|_0$ & $\footnotesize\|\textbf{B}_h-\textbf{B}\|_0$ & $\footnotesize\|\nabla(\textbf{B}_{h}-\textbf{B})\|_0$ \\
		\hline
		8  &        0.00558679   &    0.202146    &   0.00270053      &      0.00222407  &     0.0663236  \\
		16 &       0.00115566   &    0.0664967    &   0.000686662      &      0.000792706  &     0.0213601    \\
		24 &       0.000492982   &    0.0390932    &   0.000338277      &      0.000497216  &     0.0125546     \\
		32  &        0.000274777   &    0.027839    &   0.000211351      &      0.000352145  &     0.00891291    \\
		40 &       0.00017697   &    0.0217034    &   0.000150492      &      0.000270266  &     0.00692847     \\
	  \hline
	  $\footnotesize 1/h$ & $\footnotesize\textbf{u}_{\text{order}L^2}$ &  $\footnotesize\textbf{u}_{\text{order}H^1}$ & $\footnotesize p_{\text{order}L^2}$&  $\footnotesize\textbf{B}_{\text{order}L^2}$ & $\footnotesize\textbf{B}_{\text{order}H^1}$ \\
	  \hline
	  8  &     /  & / & /  &  /  &/   \\
	  16    &   2.27331   &    1.60404    &   1.97557   &    1.48835  &     1.6346    \\
24    &   2.10117   &    1.31011    &   1.74609   &    1.15035  &     1.31068    \\
32    &   2.03181   &    1.18015    &   1.63495   &    1.19918  &     1.19086    \\
40    &   1.97174   &    1.11574    &   1.52194   &    1.18595  &     1.1287    \\
			   \end{tabular*}
	 {\rule{\temptablewidth}{1pt}}
	\end{center}
	\end{table}
	
	
		\begin{table}[ht!]
		\tabcolsep 0pt \caption{The numerical results of $\nu=0.001$, $\mu=1.0$, $\sigma=1000$ and $T=1.0$ for different $h$.}\label{Ta3} \vspace*{-18pt}
		\begin{center}
		\def\temptablewidth{1.0\textwidth}
		{\rule{\temptablewidth}{1pt}}
		\begin{tabular*}{\temptablewidth}{@{\extracolsep{\fill}}llllll}
		$\footnotesize1/h$ & $\footnotesize\|\textbf{u}_h-\textbf{u}\|_0$ & $\footnotesize\|\nabla(\textbf{u}_{h}-\textbf{u})\|_0$ &  $\footnotesize\|p_h^1-p\|_0$ & $\footnotesize\|\textbf{B}_h-\textbf{B}\|_0$ & $\footnotesize\|\nabla(\textbf{B}_{h}-\textbf{B})\|_0$ \\
		\hline
		8  &        0.0494654   &    2.73996    &   0.0209561      &      0.0359843  &     1.27368   \\
		16 &        0.00503382   &    0.498632    &   0.00224957      &      0.00373482  &     0.189918    \\
		24 &        0.00147128   &    0.205262    &   0.000885435      &      0.00120948  &     0.0812076      \\
		32  &       0.000624692   &    0.111706    &   0.000399888      &      0.000557328  &     0.03984      \\
		40 &        0.000348353   &    0.0713689    &   0.00025601      &      0.000370717  &     0.023794     \\
	  \hline
	  $\footnotesize 1/h$ & $\footnotesize\textbf{u}_{\text{order}L^2}$ &  $\footnotesize\textbf{u}_{\text{order}H^1}$ & $\footnotesize p_{\text{order}L^2}$&  $\footnotesize\textbf{B}_{\text{order}L^2}$ & $\footnotesize\textbf{B}_{\text{order}H^1}$ \\
	  \hline
	  8  &     /  & / & /  &  /  &/   \\
	  16    &   3.29669   &    2.45811    &   3.21965   &    3.26826  &     2.74555    \\
24    &   3.03366   &    2.18905    &   2.29962   &    2.78077  &     2.09534    \\
32    &   2.9777   &    2.11487    &   2.7631   &    2.69323  &     2.47543    \\
40    &   2.61733   &    2.00773    &   1.99856   &    1.82714  &     2.30989    \\
			   \end{tabular*}
	 {\rule{\temptablewidth}{1pt}}
	\end{center}
	\end{table}
		\begin{table}[ht!]
		\tabcolsep 0pt \caption{The numerical results of $\nu=0.002$, $\mu=1.0$, $\sigma=2000$ and $T=1.0$ for different $h$.}\label{Ta4} \vspace*{-18pt}
		\begin{center}
		\def\temptablewidth{1.0\textwidth}
		{\rule{\temptablewidth}{1pt}}
		\begin{tabular*}{\temptablewidth}{@{\extracolsep{\fill}}llllll}
		$\footnotesize1/h$ & $\footnotesize\|\textbf{u}_h-\textbf{u}\|_0$ & $\footnotesize\|\nabla(\textbf{u}_{h}-\textbf{u})\|_0$ &  $\footnotesize\|p_h^1-p\|_0$ & $\footnotesize\|\textbf{B}_h-\textbf{B}\|_0$ & $\footnotesize\|\nabla(\textbf{B}_{h}-\textbf{B})\|_0$ \\
		\hline
		8  &         0.0215627   &    1.11785    &   0.0103966      &      0.0174424  &     0.452495   \\
		16 &         0.00242045   &    0.229028    &   0.000931281      &      0.00134777  &     0.0747986    \\
		24 &        0.000791297   &    0.101216    &   0.000411105      &      0.000583189  &     0.0325257     \\
		32  &       0.000399139   &    0.0589293    &   0.000283206      &      0.000425523  &     0.019387    \\
		40 &       0.000250524   &    0.0396481    &   0.000217605      &      0.000335836  &     0.0134441     \\
	  \hline
	  $\footnotesize 1/h$ & $\footnotesize\textbf{u}_{\text{order}L^2}$ &  $\footnotesize\textbf{u}_{\text{order}H^1}$ & $\footnotesize p_{\text{order}L^2}$&  $\footnotesize\textbf{B}_{\text{order}L^2}$ & $\footnotesize\textbf{B}_{\text{order}H^1}$ \\
	  \hline
	  8  &     /  & / & /  &  /  &/   \\
	 16    &   3.15519   &    2.28713    &   3.48075   &    3.69395  &     2.59682    \\
24    &   2.75741   &    2.01396    &   2.01673   &    2.06601  &     2.05386    \\
32    &   2.37889   &    1.88026    &   1.29544   &    1.09563  &     1.79862    \\
40    &   2.08724   &    1.77596    &   1.18083   &    1.06074  &     1.64048    \\
			   \end{tabular*}
	 {\rule{\temptablewidth}{1pt}}
	\end{center}
	\end{table}

\newpage
\subsection{Hartmann flow }
Next, we introduce the 2D Hartmann flow, which considers the 2D domain $\Omega=[0,10]\times[-1,1]$ and the external magnetic field $\textbf{B}^{d}=(0,1)$. The source forces are taken $\textbf{f}=\textbf{g}=\textbf{0}$ and the exact solution is as follows
\begin{align*}
&\textbf{u}=(u_1(y),0), \ \  \textbf{H}=(H_1(y),1),\nn\\
&p(x,y)=-Gx- \frac{\mu \cdot H_1(y)^2}{2}+p_0,
\end{align*}
with
\begin{align*}
u_1(y)=\frac{G}{\nu Ha \cdot  \tanh(Ha)}\left(1-\frac{\cosh(y Ha)}{\cosh(Ha)}\right), \ \  H_1(y)=\frac{G}{\mu}\left(\frac{\sinh(y Ha)}{\sinh(Ha)}-y \right),
\end{align*}
where the Hartmann number is denoted by $Ha=\sqrt{\frac{\sigma}{\nu}\mu^2}$.

The walls of the duct are either conducting or insulating. We impose no-slip boundary conditions on the wall, the tangential part of the magnetic field on the boundary and Neumann boundary conditions on the inlet and the outlet:
\begin{align*}
&\textbf{u}=0,\ \ \text{on }  y=\pm 1,\nn\\
&(p\textbf{I}-\nu\nabla \textbf{u})\textbf{n}=p_{d}\textbf{n},\ \ \text{on } x=0 \text{ and } x=10,\nn\\
&\textbf{n}\times \textbf{B}=\textbf{n}\times\textbf{B}_{d} \text { on }\partial\Omega.
\end{align*}
Here, $p_{d}(x,y)=p(x,y)$, $p_{0}$ is a constant and $\textbf{I}$ is identity matrix. The initial conditions are $\textbf{u}_{0}=(0,0)$ and $\textbf{H}_{0}=(0,0)$. We assume $G=1$ and choose $h=\frac{1}{24}$ and $\Delta t=\frac{1}{9}h$. It will stop when the error between the solutions of two steps $\|\textbf{u}_{h}^{n+1}-\textbf{u}_{h}^{n}\|_{0}+\|\textbf{u}_{h}^{n+1}-\textbf{u}_{h}^{n}\|_{0}\leq 1e^{-6}$.
The analytical solutions of the first components $u_{1}(y)$ and $H_{1}(y)$ are compared with numerical ones $u_{1}(y_{k})$ and $H_{1}(y_{k})$ $(y_{k}=-1+0.1k,k=0,\ldots,20 )$ in Figure{1}, when we choose $Ha=1$ ($\nu=1.0$, $\mu=1.0$, $\sigma=1.0$).
Figure 2 shows the approximation of the analytical solutions of the first components $u_{1}(y)$ and $H_{1}(y)$ and numerical ones $u_{1}(y_{k})$ and $H_{1}(y_{k})$ ($y_{k}=-1+0.1k$, $k=0,\ldots,20$) with $Ha=10$ ($\nu=0.1$, $\mu=1.0$, $\sigma=10.0$). The exact solutions are consistent with the numerical solutions along the vertical line passing the geometric center.
\begin{center}
\begin{figure}[!ht]
	\subfigure[Velocity]{\includegraphics[width=0.5 \textwidth, height=50mm]{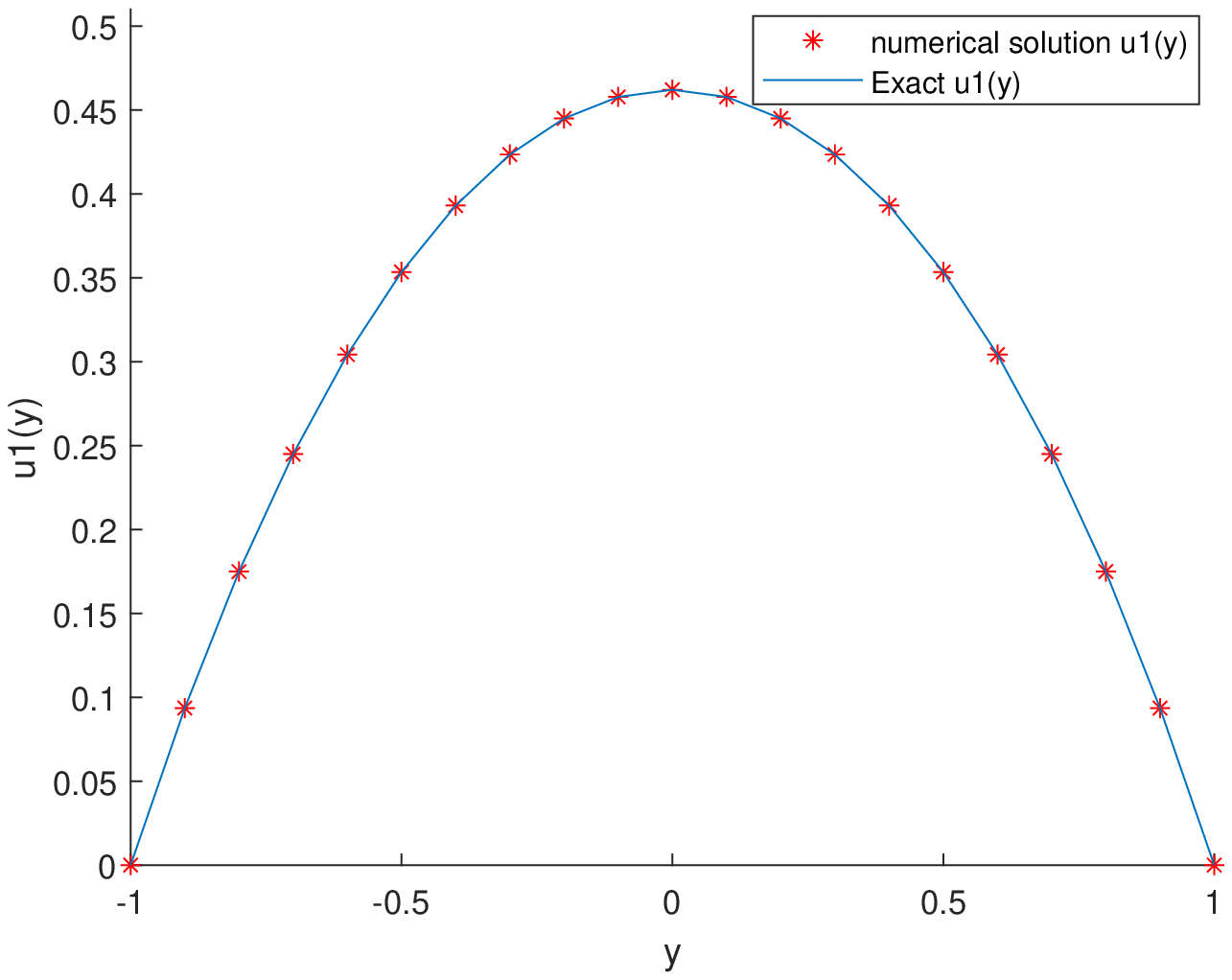}}
	\subfigure[Magnetic]{\includegraphics[width=0.5 \textwidth, height=50mm]{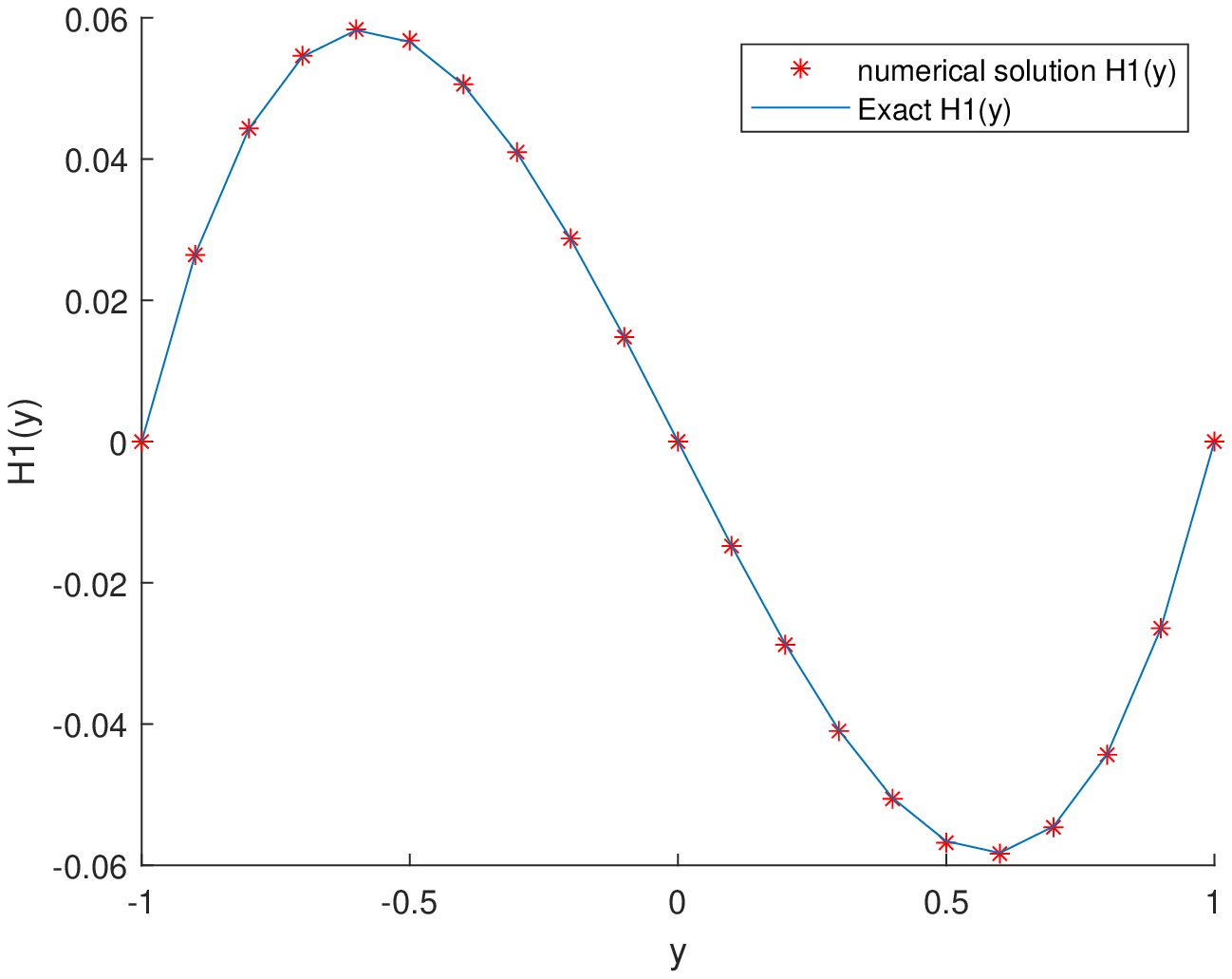}}
	\caption{The Numerical results of slicing along $-1\leq y\leq 1$, $x=5$ with $Ha=1$ ($\nu=1.0$, $\mu=1.0$, $\sigma=1.0$).}\label{F1}
\end{figure}	
\end{center}
\begin{center}
\begin{figure}[!ht]
	\subfigure[Velocity]{\includegraphics[width=0.5 \textwidth, height=50mm]{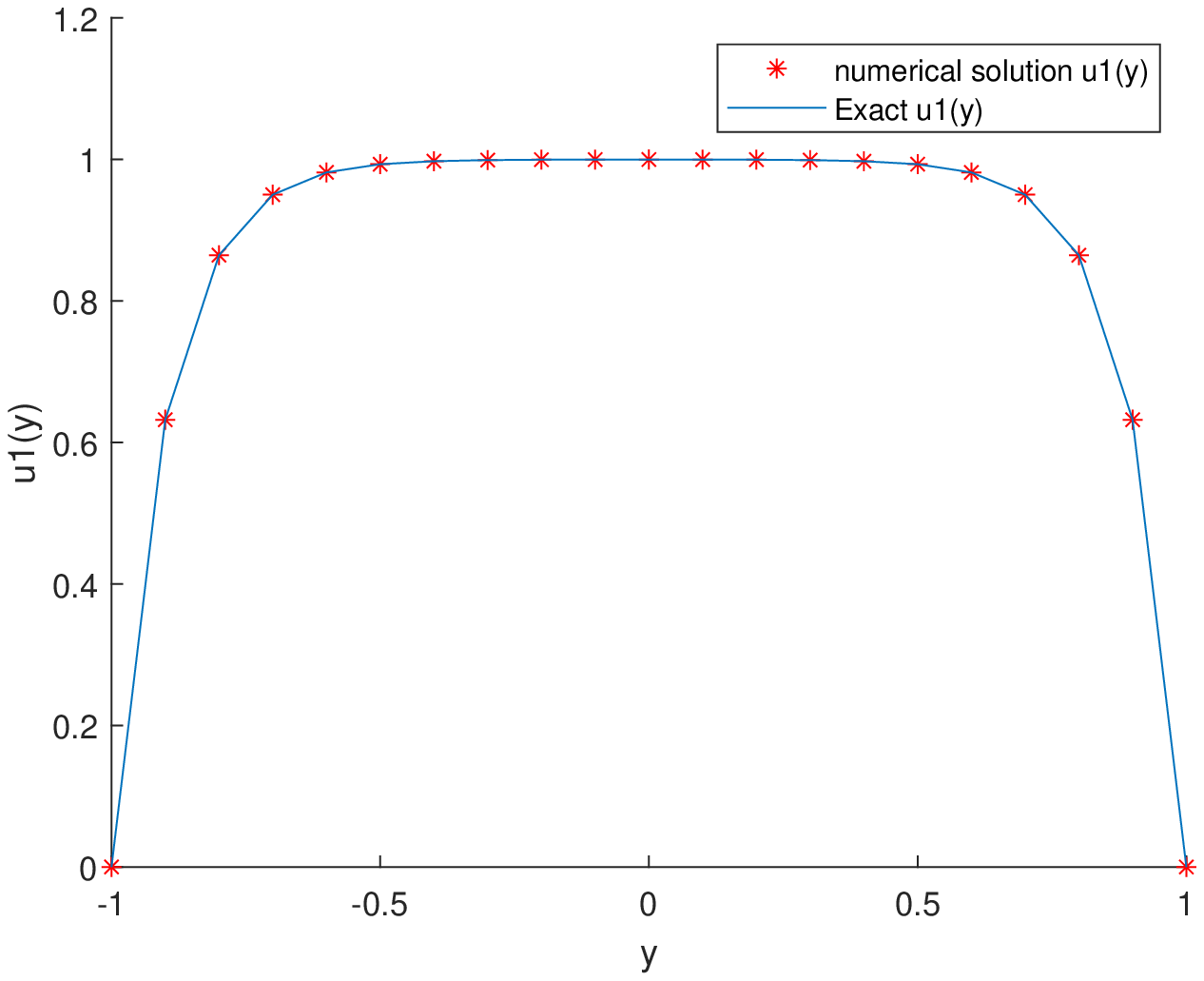}}
	\subfigure[Magnetic]{\includegraphics[width=0.5 \textwidth, height=50mm]{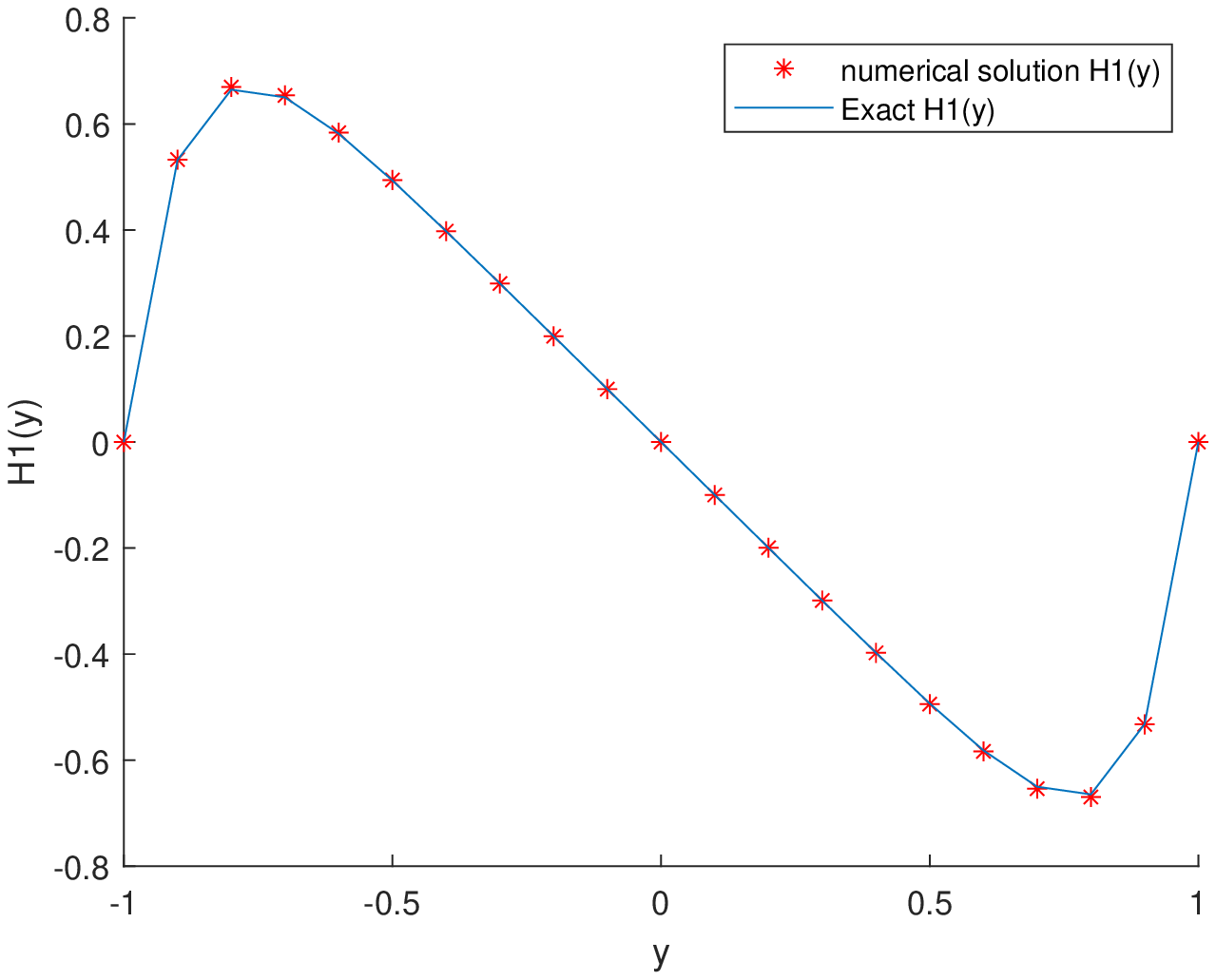}}
	\caption{The Numerical results of slicing along $-1\leq y\leq 1$, $x=5$ with $Ha=10$ ($\nu=0.1$, $\mu=1.0$, $\sigma=10.0$).}\label{F2}
\end{figure}	
\end{center}

	\end{document}